# Visco-Cosserat periporomechanics for dynamic shear bands and crack branching in porous media

Xiaoyu Song[1], Hossein Pashazad

[a]Engineering School of Sustainable Infrastructure & Environment University of Florida, Gainesville, Florida 32611 USA

**Abstract**

Periporomechanmics is a strong nonlocal framework for modeling the mechanics and physics of variably saturated porous media with evolving discontinuities. In periporomechanics, the horizon serves as a mathematical nonlocal parameter that lacks a precise physical meaning. In this article, as a new contribution we formulate a Cosserat periporomechanics paradigm for modeling dynamic shear banding and crack branching in dry porous media incorporating a micro-structure based length scale. In this micro-periporomechanics framework, each material point has both translational and rotational degrees of freedom in line with the Cosserat continuum theory. We formulate a stabilized Cosserat constitutive correspondence principle via the energy method through which classical viscous material models for porous media can be used in the proposed Cosserat periporomechanics. We have numerically implemented the micro-periporomechanics model through an explicit Lagrangian meshfree algorithm for dynamic problems. Two benchmark numerical examples are presented to test the computational micro-periporomechanics paradigm in modeling shear bands and mode-I cracks. We then present two numerical examples to show the efficacy of the proposed micro-periporomechanics for modeling the characteristics of shear banding bifurcation and dynamic crack branching in dry porous media.

## 1. Introduction

Periporomechanmics (e.g., [1–8]) is a strong non-local reformulation of classical poromechanics (e.g., [9–13]) for modeling the mechanics and physics of variably saturated porous media with evolving discontinuities. In line with the peridynamics for solids [14–17], the motion equation and mass balance equations of periporomechanics are in the form of integrodifferential equations (integration in space and difference in time) in lieu of partial differential equations through the peridynamic (PD) state and effective force state concept [3, 15]. This salient feature of periporomechanics makes it a legitimate and robust numerical tool for modeling discontinuities and progressive failure in porous media such as shear bands and cracks. It is noted that the previous periporomechanics was developed for non-polar porous materials in which material points have two kinds of degrees of freedom, i.e., translational displacement and pore fluid pressure. It is known that shear bands in porous media such as soils have a finite thickness on the order of several particle sizes and involve particle rotations in the banded deformation zone [18, 19]. Meanwhile, in periporomechanics the horizon serves as a mathematical non-local parameter that lacks a clear physical meaning and direct relation to the micro-structure of porous media. Therefore, the formulation of a visco-Cosserat periporomechanics paradigm significantly contributes to realistically modeling dynamic shear bands and crack branching in porous media through periporomechanics.

[1] Corresponding author

*Email address:* xysong@ufl.edu (Xiaoyu Song)

In this article, as a new contribution, we formulate a visco-Cosserat periporomechanics paradigm incorporating a micro-structure-based length scale for modeling dynamic shear banding and crack branching in dry porous media. The viscosity is included to study the rate dependency of shear banding and cracking in porous media (e.g, [9, 20–22]). In this micro-polar periporomechanics framework, each material point under the dry condition has both translational and rotational degrees of freedom following the Cosserat continuum. It is noted that the Cosserat continuum theory that incorporates the micro-structure of materials has been successfully applied for modeling shear bands in porous media using the classical finite element method while mitigating its pathological issue (e.g., [23–27], among others). We refer to the distinguished literature on this subject (e.g., [18, 19, 23, 28]). Cosserat continuum theory has also been used in peridynamics for modeling the mechanical behavior of solids (e.g., [29–31], among others). For instance, Gerstle et al. [29] first proposed a micropolar PD model in which the bond-based material model was adopted for modeling concrete with any Poisson's ratio. Moreover, a micropolar lattice model with an anisotropic damage parameter was proposed in [30]. Chowdhury et al. [31] formulated a micropolar state-based PD framework for modeling cracks in which a correspondence material model was proposed. It is noted that the framework in [31] is limited for isotropic elastic materials. It is known that porous media could show a rate dependency (e.g., viscosity) in the process of shear banding and cracking (e.g., fault creep) [32–34]. In this study, as an original contribution, we develop a visco-Cosserat periporomechanics framework for dry porous media. This new framework incorporates a physics-based material length scale (i.e., Cosserat length scale) and rate-dependency for more realistic modeling of progressive localized failure and cracking of porous media under dry conditions.

We note that no PD material model in the periporomechanics framework is available for modeling porous media. Thus, in this new visco-Cosserat periporomechanics framework, the multiphase correspondence principle for non-polar porous media [3, 5] is reformulated to include the rotational degree of freedom. Through the Cosserat multiphase correspondence principle, advanced elastoplastic constitutive models for geomaterials (e.g., [35–37]) can be reformulated for micropolar porous media [28] and then incorporated into the Cosserat periporomechanics paradigm. Nonetheless, the standard Cosserat correspondence principle inherits the zero-energy mode instability from the non-polar constitutive correspondence principle. In this study, we formulate a stabilized Cosserat constitutive correspondence principle to circumvent this issue using the energy method (e.g., [5, 38]). With this stabilized micro-polar correspondence principle, classical viscous constitutive models for micro-polar porous media can be incorporated into the proposed Cosserat periporomechanics. We refer to [7] for a comprehensive review of methods for mitigating zeroenergy mode instability associated with the original constitutive correspondence principle in PD for solids. We have numerically implemented the proposed Cosserat periporomechanics paradigm through an explicit Lagrangian meshfree algorithm for dynamic problems [7, 11, 39–41]. Two numerical examples inspired by experimental works in the literature are presented to test the implemented computational micro-periporomechanics paradigm in modeling shear banding and mode I cracking in porous media. We then present two numerical examples to demonstrate the efficacy and robustness of the Cosserat periporomechanics for modeling the characteristics of shear banding bifurcation (e.g., inclination angle and thickness of the two conjugate shear bands) and dynamic crack branching (e.g., the timing of branching) in porous media under drained conditions (i.e., single-phase porous media).

The reminder of this article is organized as follows. Section 2 deals with the mathematical formulation of the visco-Cosserat periporomechanic paradigm that includes the governing equations, the stabilized Cosserat constitutive correspondence principle, and the Cosserat visco-plasticity and visco-elasticity, and a bilinear damage model. Section 3 presents the numerical implementation of the proposed Cosserat periporomechanmics paradigm through an explicit Newmark scheme with an augmented energy criterion for numerical stability. Section 4 deals with numerical examples to validate the implemented Cosserat periporomechanics model and demonstrate its efficacy and robustness in modeling dynamic shear banding and crack branching in dry porous media. For sign convection, the assumption in continuum mechanics is adopted, i.e., the tensile force and deformation under tension is positive.

## 2. Mathematical formulation

In this section, we present the mathematical formulation of the visco-Cosserat periporomechanics for dry porous media. First, we present the governing equations of Cosserat periporomechanics for dry porous materials. Second, we develop a Cosserat constitutive correspondence principle through which the classical visco-Cosserat material models can be incorporated into the proposed Cosserat periporomechanics. Third, we present an energy-based stabilization scheme in the Cosserat periporomechanmics framework to mitigate the zero-energy mode deformation instability. Finally, the classical Cosserat visco-plasticity, visco-elasticity, and damage models are presented.

### *2.1. Governing equations of Cosserat periporomechanics*

In periporomechanics, a porous material body can be conceptualized as a collection of a finite number of material points with poromechanical and physical interactions at a finite distance. The scope of the present contribution is to formulate visco-Cosserat periporomechnics for a dry porous material. In line with the classical Cosserat continuum theory (e.g., [18, 42–44]) in the proposed Cosserat periporomechanmics (CPPM) each material point is endowed with translational and rotational degrees of freedom. Let $\boldsymbol{X}$ denote a material point and let $\boldsymbol{X}'$ denote its neighboring material point in the initial configuration. For notations, a state or parameter without a prime is associated with material point $\boldsymbol{X}$, and a state or parameter with a prime is associated with material point $\boldsymbol{X}'$. For example, we can define $\boldsymbol{y}$ and $\boldsymbol{y}'$ as the spatial locations of material points $\boldsymbol{X}$ and $\boldsymbol{X}'$ in the deformed configuration. The partial density of a solid skeleton is defined as

$$\rho^s = \phi\rho_s, \tag{1}$$

where $\phi$ is the volume fraction of the solid skeleton (i.e., one minus porosity) and $\rho_s$ is the intrinsic density of the solid skeleton. Here, the porosity is the fraction between the pore space and the porous media's total volume.

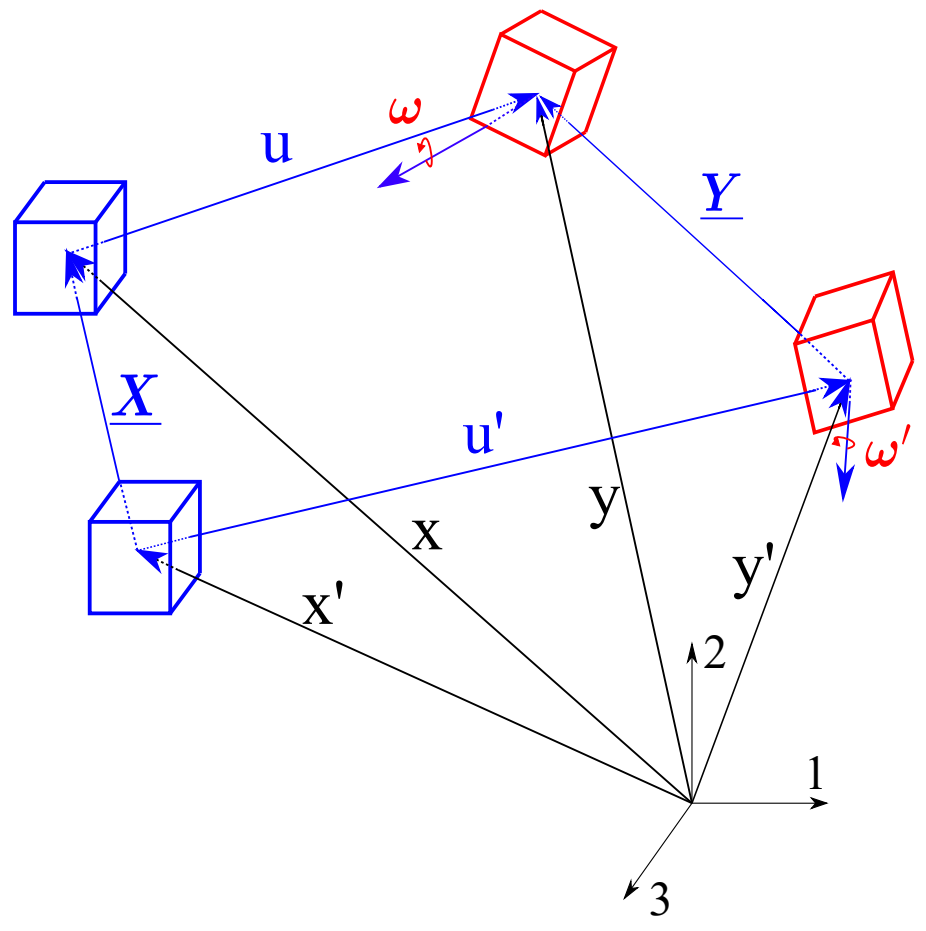

Figure 1: Kinematics of two material points x and x' in Cosserat periporomechanics.

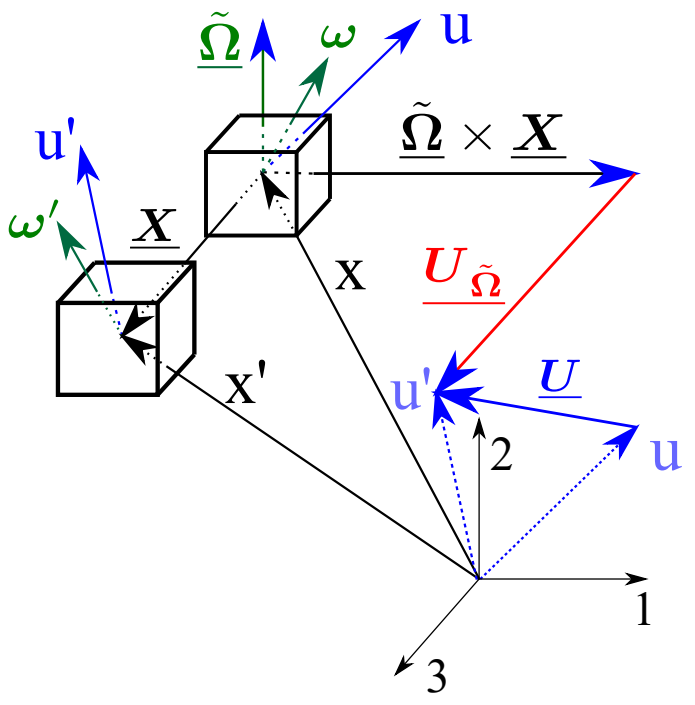

Figure 2: Schematic of composite states in Cosserat periporomechanics.

Figure 1 plots the kinematics of CPPM for material points. In line with the non-polar PPM, the deformation state and the relative displacement state associated with bond $\boldsymbol{\xi} = \boldsymbol{X}' - \boldsymbol{X}$ at material point $\boldsymbol{X}'$ read

$$\underline{\boldsymbol{Y}} = \boldsymbol{y}' - \boldsymbol{y}, \tag{2}$$

$$\underline{\mathbf{U}} = \mathbf{u}' - \mathbf{u}, \tag{3}$$

where $\boldsymbol{u}$ and $\boldsymbol{u}'$ are the displacement vectors of material points $\boldsymbol{X}$ and $\boldsymbol{X}'$, respectively. It is assumed that material points have no macro rotations. Given the micro-rotations $\boldsymbol{\omega}$ and $\boldsymbol{\omega}'$ at $\boldsymbol{X}$ and $\boldsymbol{X}'$, respectively, the micro-rotation state and the the mean micro-rotation state associated with bond $\boldsymbol{\xi}$ at $\boldsymbol{X}$ can be defined as

$$\underline{\boldsymbol{\Omega}} = \boldsymbol{\omega}' - \boldsymbol{\omega}, \tag{4}$$

$$\overline{\underline{\boldsymbol{\Omega}}} = \frac{\mathbf{1}}{\mathbf{2}}(\boldsymbol{\omega}' + \boldsymbol{\omega}). \tag{5}$$

Figure 2 plots the composite state in Cosserat periporomechanics. Referring to Figure 2, the composite deformation state in Cosserat periporomechanics can be defined as

$$\underline{\tilde{\mathbf{Y}}} = \underline{\mathbf{Y}} - \overline{\underline{\boldsymbol{\Omega}}} \times \underline{\mathbf{X}}, \tag{6}$$

$$\underline{\tilde{\mathbf{U}}} = \underline{\mathbf{U}} - \overline{\underline{\boldsymbol{\Omega}}} \times \underline{\mathbf{X}}, \tag{7}$$

where $\times$ is the cross product operator of two vectors.

Next, we derive the governing equations from a free energy of a polar porous body following [3]. Assuming no heat flux, the free energy density of the solid skeleton in Cosserat periporomechanics under dry condition (i.e., single phase) can be defined as

$$\overline{\mathcal{W}} = \overline{\mathcal{W}}(\underline{\tilde{\mathrm{Y}}}, \underline{\Omega}), \tag{8}$$

Let $\overline{\underline{\boldsymbol{\mathcal{T}}}}$ be the effective force state and $\underline{\boldsymbol{\mathcal{M}}}$ be the moment state. Through Fr´echet derivatives [3, 15], the effective force state and moment state on a Cosserat periporomechanics bond can be defined as

$$\overline{\underline{\boldsymbol{\mathcal{T}}}} = \frac{\partial \overline{\mathcal{W}}}{\partial \underline{\tilde{\mathrm{Y}}}}, \tag{9}$$

$$\underline{\boldsymbol{\mathcal{M}}} = \frac{\partial \overline{\mathcal{W}}}{\partial \underline{\Omega}} \quad . \tag{10}$$

Then, the variational form of the free energy of the solid skeleton can be written as

$$\begin{aligned} \Delta\overline{\mathcal{W}} &= \frac{\partial \overline{\mathcal{W}}}{\partial \underline{\tilde{Y}}} \bullet \Delta\underline{\tilde{Y}} + \frac{\partial \overline{\mathcal{W}}}{\partial \underline{\Omega}} \bullet \Delta\underline{\Omega} \\ &= \overline{\underline{\boldsymbol{\mathcal{T}}}} \bullet \Delta\underline{\tilde{Y}} + \underline{\boldsymbol{\mathcal{M}}} \bullet \Delta\underline{\boldsymbol{\Omega}}, \end{aligned} \tag{11}$$

where $\bullet$ is the inner product of vector states [15].

Assuming small deformation and replacing $\underline{\boldsymbol{X}}$ by $\underline{\boldsymbol{Y}}$ in 6, equation (11) can be written as

$$\Delta\overline{\mathcal{W}} = \overline{\underline{\boldsymbol{\mathcal{T}}}} \cdot \left(\boldsymbol{\Delta}\underline{\mathbf{Y}} - \boldsymbol{\Delta}\overline{\underline{\boldsymbol{\Omega}}} \times \underline{\mathbf{Y}}\right) + \underline{\boldsymbol{\mathcal{M}}} \cdot \Delta\underline{\boldsymbol{\Omega}}. \tag{12}$$

The total potential energy $\mathcal{W}$ of a bounded single-phase porous body $\mathcal{B}$ can be written as

$$\mathcal{W} = \int_{\mathcal{B}} \left(\overline{\mathcal{W}} - \rho^s \mathrm{g} \cdot \vec{\mathrm{y}} - \mathrm{l} \cdot \omega\right) \, \mathrm{dV}, \tag{13}$$

where $\rho^s$ is partial density of the solid skeleton as defined in (1), $\boldsymbol{g}$ is gravity acceleration, and $\boldsymbol{l}$ is body couple density.

Next, we derive the equilibrium equation of the Cosserat periporomechanics for a dry porous body from (13). The first-order variational form of $W$ can be written as

$$\begin{aligned} 0 &= \int_{\mathcal{B}} \left[\underline{\boldsymbol{\mathcal{T}}} \cdot \left(\boldsymbol{\Delta}\underline{\mathbf{Y}} - \boldsymbol{\Delta}\underline{\boldsymbol{\Omega}} \times \underline{\mathbf{Y}}\right) + \underline{\boldsymbol{\mathcal{M}}} \cdot \boldsymbol{\Delta}\underline{\boldsymbol{\Omega}} - \rho^s \vec{\mathbf{g}} \cdot \boldsymbol{\Delta}\mathbf{y} - \mathbf{l} \cdot \boldsymbol{\Delta\omega}\right] \, \mathrm{dV} \\ &= \int_{\mathcal{B}} \left[\int_{\mathcal{H}} \underline{\boldsymbol{\mathcal{T}}} \cdot \left(\Delta\vec{y'} - \Delta\vec{y}\right) dV' + \int_{\mathcal{H}} \underline{\boldsymbol{\mathcal{T}}}' \times \underline{\boldsymbol{Y}} \cdot \left(\frac{\Delta\omega' + \Delta\omega}{2}\right) dV' \right. \\ &\quad \left. + \int_{\mathcal{H}} \underline{\boldsymbol{\mathcal{M}}} \cdot (\Delta\omega' - \Delta\omega) \, dV' - \rho^s g \cdot \Delta y - l \cdot \Delta\omega \right] dV \end{aligned} \tag{14}$$

By interchanging the dummy variable $X \longleftrightarrow X'$ [3] in the first, third and fifth terms in the integrand of (14) leads to

$$\begin{aligned} 0 = \int_{\mathcal{B}} &\left[\int_{\mathcal{H}} \left(\underline{\boldsymbol{\mathcal{T}}}' - \underline{\boldsymbol{\mathcal{T}}}\right) \mathrm{dV}' \cdot \Delta \mathrm{y} + \frac{1}{2}\int_{\mathcal{H}} \underline{\mathrm{Y}} \times \left(\underline{\boldsymbol{\mathcal{T}}}' - \underline{\boldsymbol{\mathcal{T}}}\right) \mathrm{dV}' \cdot \Delta\omega \right. \\ &\left. + \int_{\mathcal{H}} \left(\underline{\boldsymbol{\mathcal{M}}}' - \underline{\boldsymbol{\mathcal{M}}}\right) dV' \cdot \Delta\omega - \rho^s g \cdot \Delta y - l \cdot \Delta\omega \right] dV \end{aligned} \tag{15}$$

Note that equation (15) must hold for any $\Delta\boldsymbol{y}$ and $\Delta\boldsymbol{\omega}$. Thus, we have

$$0 = \int_{\mathcal{H}} \left(\underline{\boldsymbol{\mathcal{T}}}' - \underline{\boldsymbol{\mathcal{T}}}\right) dV' + \rho^s g, \tag{16}$$

$$0 = \int_{\mathcal{H}} \left(\underline{\boldsymbol{\mathcal{M}}}' - \underline{\boldsymbol{\mathcal{M}}}\right) dV' + \frac{1}{2}\int_{H} \underline{Y} \times \left(\underline{\boldsymbol{\mathcal{T}}}' - \underline{\boldsymbol{\mathcal{T}}}\right) dV' + l \tag{17}$$

The kinetic energy for a bounded dry porous body $\mathcal{B}$ can be written as

$$\mathcal{K} = \int_{\mathcal{B}} \rho^s \frac{\dot{\mathrm{u}} \cdot \dot{\mathrm{u}}}{2} \, \mathrm{dV} + \int_{\mathcal{B}} \mathcal{J}^s \frac{\dot{\omega} \cdot \dot{\omega}}{2} \mathrm{dV}, \tag{18}$$

where $\dot{\mathbf{u}}$ is velocity, $\mathcal{J}^s$ is the micro-inertia of the solid skeleton. With (18) the equations of motion and rotation in Cosserat periporomechanics under the dry condition can be written as

$$\rho^s \ddot{u} = \int_{\mathcal{H}} \left(\underline{\boldsymbol{\mathcal{T}}}' - \underline{\boldsymbol{\mathcal{T}}}\right) dV' + \rho^s g, \tag{19}$$

$$\mathcal{J}^s \ddot{\omega} = \int_{\mathcal{H}} \left(\underline{\boldsymbol{\mathcal{M}}} - \underline{\boldsymbol{\mathcal{M}}}'\right) dV' + \frac{1}{2}\int_{\mathcal{H}} \underline{Y} \times \left(\underline{\boldsymbol{\mathcal{T}}}' - \underline{\boldsymbol{\mathcal{T}}}\right) dV' + l, \tag{20}$$

where $\ddot{\boldsymbol{u}}$ is acceleration vector and $\ddot{\boldsymbol{\omega}}$ is angular acceleration vector. We note that in (20) $\mathcal{J}^s \ddot{\boldsymbol{\omega}}$ represents the angular momentum of a spinning material point.

To complete the equations (19) and (20), we need to introduce the constitutive models for the effective force state and the moment state. In this article, as a new contribution we develop a stabilized Cosserat correspondence principle through which classical local constitutive models for micro-polar porous media can be incorporated into the proposed Cosserat periporomechanics. In what follows we derive the Cosserat correspondence principle through an equivalence of internal energy between local micro-polar poromechanics and Cosserat periporomechanics.

*2.2. Cosserat periporomechanics correspondence principal*

This part deals with the Cosserat correspondence principle following the lines in non-polar periporomechanics [3]. To achieve this task we first derive an expression of internal energy density for micro-polar porous media in the framework of Cosserat periporomechanics By ignoring the heat source and/or heat flux, let us multiply (19) by $\dot{\boldsymbol{u}}$ and integrate over a finite sub-region $\mathcal{P}$ within bounded micro-polar porous body $\mathcal{B}$, and we have

$$\int_{\mathcal{P}} \rho^s \ddot{u} \cdot \dot{u} dV = \int_{\mathcal{P}} \int_{\mathcal{H}} (\underline{\mathcal{T}}' - \underline{\mathcal{T}}) \cdot \dot{u} \, dV' dV + \int_{\mathcal{P}} \rho^s g \cdot \dot{u} dV. \tag{21}$$

Similarly, multiply (20) by $\dot{\boldsymbol{\omega}}$ and integrate over $\mathcal{P}$, and we have

$$\begin{aligned} \int_{\mathcal{P}} \mathcal{I}^s \ddot{\omega} \cdot \dot{\omega} dV &= \int_{\mathcal{P}} \int_{\mathcal{H}} (\underline{\mathcal{M}}' - \underline{\mathcal{M}}) \cdot \dot{\omega} \, dV' dV + \int_{\mathcal{P}} l \cdot \dot{\omega} dV \\ &+ \frac{1}{2} \int_{\mathcal{P}} \int_{\mathcal{H}} \underline{Y} \times (\underline{\mathcal{T}}' - \underline{\mathcal{T}}) \cdot \dot{\omega} \, dV' dV. \end{aligned} \tag{22}$$

Summation of (21) and (22) gives the balance of energy as

$$\begin{aligned} \int_{\mathcal{P}} (\rho^s \ddot{u} \cdot \dot{u} + \mathcal{I}^s \ddot{\omega} \cdot \dot{\omega}) dV &= \int_{\mathcal{P}} \int_{\mathcal{H}} (\underline{\mathcal{T}}' - \underline{\mathcal{T}}) \cdot \dot{u} \, dV' dV \\ &+ \int_{\mathcal{P}} \int_{\mathcal{H}} (\underline{\mathcal{M}}' - \underline{\mathcal{M}}) \cdot \dot{\omega} \, dV' dV \\ &+ \frac{1}{2} \int_{\mathcal{P}} \int_{\boldsymbol{\mathcal{H}}} \underline{\boldsymbol{Y}} \times (\underline{\mathcal{T}}' - \underline{\mathcal{T}}) \cdot \dot{\omega} \, dV' dV \\ &+ \int_{\mathcal{P}} (\rho^s g \cdot \dot{u} + l \cdot \dot{\omega}) dV. \end{aligned} \tag{23}$$

The first three terms in the right-hand side of (23) can be rewritten as

$$\begin{aligned} \int_{\mathcal{P}} \int_{\mathcal{H}} (\underline{\mathcal{T}}' - \underline{\mathcal{T}}) \cdot \dot{u} \, dV' dV &= \int_{\mathcal{P}} \int_{\mathcal{B}} (\underline{\mathcal{T}} \cdot \dot{u}' - \underline{\mathcal{T}}' \cdot u) \, dV' dV - \int_{\mathcal{P}} \int_{\mathcal{B}} \underline{\mathcal{T}} (\dot{u}' - \dot{u}) \, dV' dV \\ &= \int_{\mathcal{P}} \int_{\mathcal{B} \setminus \mathcal{P}} (\underline{\mathcal{T}} \cdot \dot{u}' - \underline{\mathcal{T}}' \cdot u) \, dV' dV - \int_{\mathcal{P}} \int_{\mathcal{B}} \underline{\mathcal{T}} (\dot{u}' - \dot{u}) \, dV' dV \end{aligned} \tag{24}$$

$$\begin{aligned} & \int_{\mathcal{P}} \int_{\mathcal{H}} (\underline{\mathcal{M}}' - \underline{\mathcal{M}}) \cdot \dot{\boldsymbol{\omega}} \, dV' dV \\ &= \int_{\mathcal{P}} \int_{\mathcal{B}} (\underline{\mathcal{M}} \cdot \dot{\omega}' - \underline{\mathcal{M}}' \cdot \dot{\omega}) \, dV' dV - \int_{\mathcal{P}} \int_{\mathcal{B}} \underline{\mathcal{M}} \cdot (\dot{\omega}' - \dot{\omega}) \, dV' dV \\ &= \int_{\mathcal{P}} \int_{\mathcal{B} \setminus \mathcal{P}} (\underline{\mathcal{M}} \cdot \dot{\omega}' - \underline{\mathcal{M}}' \cdot \dot{\omega}) \, dV' dV - \int_{\mathcal{P}} \int_{\mathcal{B}} \underline{\mathcal{M}} \cdot (\dot{\omega}' - \dot{\omega}) \, dV' dV \end{aligned} \tag{25}$$

$$
\begin{aligned}
&\frac{1}{2}\int_{\mathcal{P}}\int_{\mathcal{H}} \underline{Y} \times \left(\underline{\bar{\mathcal{T}}}' - \underline{\bar{\mathcal{T}}}\right) \cdot \dot{\omega}\, dV' dV \\
&= \int_{\mathcal{P}}\int_{\mathcal{B}} \left(\underline{\bar{\mathcal{T}}} \cdot \frac{\dot{\omega}'}{2} \times \underline{Y}' - \underline{\bar{\mathcal{T}}}' \cdot \frac{\dot{\omega}}{2} \times \underline{\boldsymbol{Y}}\right) dV' dV - \int_{\mathcal{P}}\int_{\mathcal{B}} \underline{\bar{\mathcal{T}}} \cdot \frac{\dot{\omega}' + \dot{\omega}}{2} \times \underline{Y}\, dV' dV \\
&= \int_{\mathcal{P}}\int_{\mathcal{B}\setminus\mathcal{P}} \left(\underline{\bar{\mathcal{T}}} \cdot \frac{\dot{\omega}'}{2} \times \underline{Y}' - \underline{\bar{\mathcal{T}}}' \cdot \frac{\dot{\omega}}{2} \times \underline{\boldsymbol{Y}}\right) dV' dV - \int_{\mathcal{P}}\int_{\mathcal{B}} \underline{\bar{\mathcal{T}}} \cdot \frac{\dot{\omega}' + \dot{\omega}}{2} \times \underline{Y}\, dV'.
\end{aligned}
\tag{26}
$$

It follows from (24), (25), and (26) that we can express the balance of energy as

$$
\dot{\mathcal{K}} + \mathcal{W}_{abs} = \mathcal{W}_{sup}, \tag{27}
$$

where $\dot{\mathcal{K}}$ is the rate of kinetic energy, $\mathcal{W}_{abs}$ is absorbed power, and $\mathcal{W}_{sup}$ is supplied power. The latter two terms are written as

$$
\mathcal{W}_{abs} = \int_{\mathcal{P}}\int_{\mathcal{B}} \left[\underline{\bar{\mathcal{T}}} \cdot \left(\dot{\underline{U}} - \frac{\dot{\omega}' + \dot{\omega}}{2} \times \underline{Y}\right) + \underline{\underline{\mathcal{M}}} \cdot \left(\dot{\omega}' - \dot{\omega}\right)\right] dV' dV. \tag{28}
$$

$$
\begin{aligned}
\mathcal{W}_{sup} = &\int_{\mathcal{P}}\int_{\mathcal{B}\setminus\mathcal{P}} \left[\underline{\bar{\mathcal{T}}} \cdot \left(\dot{u}' + \frac{\dot{\omega}'}{2} \times \underline{Y}'\right) - \underline{\bar{\mathcal{T}}}' \cdot \left(\dot{u} + \frac{\dot{\omega}}{2} \times \underline{Y}\right)\right] dV' dV \\
&+ \int_{\mathcal{P}}\int_{\mathcal{B}\setminus\mathcal{P}} \left(\underline{\boldsymbol{M}} \cdot \dot{\boldsymbol{\omega}}' - \underline{\boldsymbol{M}'} \cdot \dot{\boldsymbol{\omega}}\right) dV' dV + \int_{\mathcal{P}} \left(\rho^s \boldsymbol{g} \cdot \dot{\boldsymbol{u}} + \boldsymbol{l} \cdot \dot{\boldsymbol{\omega}}\right) dV.
\end{aligned}
\tag{29}
$$

The absorbed energy is equal to the internal energy without considering other sources of energy (e.g., thermal energy) (e.g., [3]). Thus, it follows from (28) and the assumption of small deformation that the internal energy density can be written as

$$
\begin{aligned}
\dot{\mathcal{E}} &= \int_{\mathcal{B}} \left[\underline{\bar{\mathcal{T}}} \cdot \left(\dot{\underline{U}} - \dot{\underline{\Omega}} \times \underline{Y}\right) + \underline{\mathcal{M}} \cdot \dot{\underline{\Omega}}\right] dV' \\
&= \int_{\mathcal{B}} \left[\underline{\bar{\mathcal{T}}} \cdot \left(\dot{\underline{U}} - \dot{\underline{\Omega}} \times \underline{X}\right) + \underline{\mathcal{M}} \cdot \dot{\underline{\Omega}}\right] dV' \\
&= \int_{\mathcal{B}} \left(\underline{\bar{\mathcal{T}}} \cdot \dot{\underline{\bar{U}}} + \underline{\mathcal{M}} \cdot \dot{\underline{\Omega}}\right) dV' \\
&= \underline{\bar{\mathcal{T}}} \cdot \dot{\underline{\bar{U}}} + \underline{\mathcal{M}} \cdot \dot{\underline{\Omega}}
\end{aligned}
\tag{30}
$$

Under the isothermal condition and assuming small strain, the internal energy of the skeleton of a micro-polar porous material can be written as

$$
\dot{e} = \bar{\sigma} : \dot{\bar{\varepsilon}} + m : \dot{\kappa}, \tag{31}
$$

where $\bar{\boldsymbol{\sigma}}$ the effective stress tensor, $\boldsymbol{\varepsilon}$ is the micropolar strain tensor, $\boldsymbol{m}$ is the couple stress tensor, and $\boldsymbol{\kappa}$ is the wryness tensor. Following the lines in periporomechanics [3], given $\underline{\bar{\boldsymbol{U}}}$ and $\underline{\boldsymbol{\Omega}}$ the nonlocal versions of $\bar{\boldsymbol{\varepsilon}}$ and $\boldsymbol{\kappa}$ can be written as

$$
\bar{\varepsilon} = \left[\int_{\mathcal{H}} \underline{\omega} \left(\underline{\bar{U}} \otimes \underline{\xi}\right) dV'\right] K^{-1}, \tag{32}
$$

$$\kappa = \left[\int_{\mathcal{H}} \underline{\omega}\left(\underline{\Omega} \otimes \underline{\xi}\right) dV'\right] K^{-1}, \tag{33}$$

where $\underline{\omega}$ is a weighting function and $\boldsymbol{K}$ is the shape tensor [15]. The shape function is defined as

$$K = \int_{\mathcal{H}} \underline{\omega}\, \underline{\xi} \otimes \underline{\xi}\, dV' \tag{34}$$

It follows from (31), (32), and (33) that the rate of the internal energy of a micropolar porous medium can be written as

$$\begin{aligned} \dot{e} &= \int_{\mathcal{H}} \underline{\omega}\bar{\sigma} : \left(\dot{\underline{\overline{U}}} \otimes \underline{\xi}\right) K^{-1} dV' + \int_{\mathcal{H}} \underline{\omega}\, m : \left(\dot{\underline{\Omega}} \otimes \underline{\xi}\right) K^{-1} dV' \\ &= \int_{\mathcal{H}} \left(\underline{\omega}\bar{\sigma} K^{-1} \underline{\xi}\right) \cdot \dot{\underline{\overline{U}}} dV' + \int_{\mathcal{H}} \left(\underline{\omega} m K^{-1} \underline{\xi}\right) \cdot \dot{\underline{\Omega}} dV'. \end{aligned} \tag{35}$$

It follows from (36),(30) and (35) that the effective force state and the moment state in Cosserat periporomechanics can be written in terms of the effective stress tensor and moment tensor as

$$\dot{\mathcal{E}} = \dot{e}. \tag{36}$$

From the Cosserat periporomechanics correspondence principle [3] we have

$$\underline{\overline{\boldsymbol{\mathcal{T}}}} = \underline{\omega}\bar{\sigma} K^{-1} \underline{\xi}, \tag{37}$$

$$\underline{\boldsymbol{\mathcal{M}}} = \underline{\omega} m K^{-1} \underline{\xi}. \tag{38}$$

From (37) and (38) the Cosserat effective force state and the moment state can be computed from the classical constitutive models given the composite deformation state and micro-rotational state following the lines in classical Cosserate continuum mechanics for solids. In what follows, we first demonstrate the Cosserat periporomechanics correspondence principal suffers from the zeroenergy deformation instability mode as in the non-polar periporomechanics. Then we present a stabilization scheme for the Cosserate periporomechanics through which classical micro-polar constitutive models can be incorporated into the Cosserat periporomechanics developed in this article.

*2.3. Zero energy modes and stabilization scheme*

In this part, we first demonstrate the Cosserat periporomechanics correspondence principal suffers from the zero-energy deformation mode (both translational and micro-rotational) under non-uniform deformation. Then, we present a stabilization scheme to circumvent the zero-energy deformation modes. In what follows, we show the origin of zero-energy deformation modes of the Cosserat periporomechanics correspondence principal. Let us first define the non-uniform composite (micro-polar) displacement state and non-uniform micro-rotation state as follows.

$$\underline{R_1} = \underline{\overline{U}} - \bar{\varepsilon}\underline{\xi}, \tag{39}$$

$$\underline{R_2} = \underline{\Omega} - \kappa\underline{\xi}, \tag{40}$$

where $\bar{\boldsymbol{\varepsilon}}$ and $\boldsymbol{\kappa}$ are defined in (32) and (33), respectively. Substituting (39) into (32), i.e., by replacing $\underline{\overline{\boldsymbol{U}}}$ by $\underline{\boldsymbol{R_1}}$ in (32), we have

$$\left[\int_{\mathcal{H}} \underline{\omega}\left(\underline{R_1} \otimes \underline{\xi}\right) dV'\right] K^{-1} = \left[\int_{\mathcal{H}} \underline{\omega}\left(\underline{\overline{\boldsymbol{U}}} - \overline{\boldsymbol{R}}_1\underline{\boldsymbol{\xi}}\right) \otimes \underline{\boldsymbol{\xi}}) dV'\right]\right] \boldsymbol{K}^{-1}$$

$$
\begin{aligned}
&= \left[\int_{\mathcal{H}} \underline{\omega}\left(\underline{\overline{U}} \otimes \underline{\xi}\right) dV'\right] K^{-1} - \left[\int_{\mathcal{H}} \underline{\omega}\overline{\varepsilon}\,\underline{\xi} \otimes \underline{\xi} dV'\right] K^{-1} \\
&= \overline{\varepsilon} - \overline{\varepsilon} K K^{-1} \\
&= \overline{\varepsilon} - \overline{\varepsilon} \\
&= 0
\end{aligned}
\tag{41}
$$

Similarly, substituting (40)] into (33) we have

$$
\begin{aligned}
\left[\int_{\mathcal{H}} \underline{\omega}\left(\underline{R_2} \otimes \underline{\xi}\right) dV'\right] K^{-1} &= \left[\int_{\mathcal{H}} \underline{\omega}\left(\underline{\boldsymbol{\Omega}} - \boldsymbol{R}_2\underline{\boldsymbol{\xi}}\right) \otimes \underline{\boldsymbol{\xi}} dV'\right] \boldsymbol{K}^{-\mathbf{1}} \\
&= \left[\underline{\omega}(\underline{\boldsymbol{\Omega}} \otimes \underline{\boldsymbol{\xi}}) dV'\right] \boldsymbol{K}^{-\mathbf{1}} - \left[\int_{\mathcal{H}} \underline{\omega}\boldsymbol{\kappa}\underline{\boldsymbol{\xi}} \otimes \underline{\boldsymbol{\xi}} dV'\right] \boldsymbol{K}^{-\mathbf{1}} \\
&= \kappa - \kappa K K^{-1} \\
&= \kappa - \kappa \\
&= 0
\end{aligned}
\tag{42}
$$

Through (41) and (42) we have demonstrated that the nonuniform micro-polar displacement state and micro-rotational state are smoothed out in the micro-polar periporomechanics correspondence principal. Therefore, the zero-energy deformation mode instability occurs in Cosserate periporomechanics that incorporates the correspondence material models. To resolve this issue, following the lines in [5] we develop a stabilization scheme based on an energy method. We refer to [45] for a comprehensive review of other methods for stabilization schemes of correspondence materials models in the original peridynamics for solids. In this method [5], the internal energies related to the non-uniform composite displacement state and non-uniform micro-rotational state are considered in (30). In this case, the total internal energy is written as

$$
\mathcal{E} = \mathcal{E} + \mathcal{E}_{\underline{\mathcal{R}_1}} + \mathcal{E}_{\underline{\mathcal{R}_2}}, \tag{43}
$$

where $\mathcal{E}_{\underline{\mathcal{R}_1}}$ and $\mathcal{E}_{\underline{\mathcal{R}_2}}$ are the energy terms of the non-uniform composite displacement state and non-uniform micro-rotational state, respectively. The two terms are defined as

$$
\mathcal{E}_{\underline{\mathcal{R}_1}} = \frac{1}{2}\left(\underline{\alpha}\underline{R_1}\right) \cdot \underline{R_1}, \tag{44}
$$

$$
\mathcal{E}_{\underline{\mathcal{R}_2}} = \frac{1}{2}\left(\underline{\beta}\underline{R_2}\right) \cdot \underline{R_2}, \tag{45}
$$

where $\underline{\alpha}$ and $\underline{\beta}$ are the two scalar states. Following the lines in [5], assuming a micro-polar bond-based periporomechanics the two scalar states can be written as

$$
\underline{\alpha} = \frac{\mathcal{G}_1 \mathcal{C}_1}{\omega_0} \underline{\omega}, \tag{46}
$$

$$
\underline{\beta} = \frac{\mathcal{G}_2 \mathcal{C}_2}{\omega_0} \underline{\omega}, \tag{47}
$$

where

$$
\omega_0 = \int_{\mathcal{H}} \underline{\omega} \, dV', \tag{48}
$$

and $\mathcal{G}_1$ and $\mathcal{G}_2$ are positive constants on the order of 1, $\mathcal{C}_1$ and $\mathcal{C}_2$ are two material parameters. For a three-dimensional case, the two parameters can be written as [5, 29, 46]

$$\mathcal{C}_1 = \frac{12\mathcal{D}}{|\underline{\xi}|^3}, \tag{49}$$

$$\mathcal{C}_2 = \frac{\mathcal{D}}{|\underline{\xi}|}, \tag{50}$$

where $\mathcal{D}$ is a parameter that depends on the internal length scale (i.e., the horizon and material properties). In the three-dimensional case it reads

$$\mathcal{D} = \frac{E(1-4\nu)}{4\pi\delta^2(1-\nu-2\nu^2)}, \tag{51}$$

where $E$ and $\nu$ are Young's modulus and Poisson's ration, respectively. It follows from (44), (45), (37), and (38) the effective force state and moment state with stabilization removing zero energy modes can be written as

$$\underline{\overline{\mathcal{T}}}^s = \underline{\omega}\bar{\sigma}K^{-1}\underline{\xi} + \underline{\alpha}\underline{\underline{R_1}}, \tag{52}$$

$$\underline{\underline{\mathcal{M}}}^s = \underline{\omega}mK^{-1}\underline{\xi} + \beta\underline{\underline{R_2}}. \tag{53}$$

Substituting (52) and (53) the governing equations for Cosserat periporomechanics can be written as

$$\rho^s\ddot{u} = \int_{\mathcal{H}} \left(\underline{\overline{\mathcal{T}}}'^s - \underline{\overline{\mathcal{T}}}^s\right) dV' + \rho^s g, \tag{54}$$

$$\mathcal{I}^s\ddot{\omega} = \int_{\mathcal{H}} \left(\underline{\underline{\mathcal{M}}}^s - \underline{\underline{\mathcal{M}}}'^s\right) dV' + \frac{1}{2}\int_{\mathcal{H}} \underline{Y} \times \left(\underline{\overline{T}}'^s - \underline{\overline{T}}^s\right) dV' + l. \tag{55}$$

In what follows, we introduce the classical micro-polar visco-plastic, visco-elastic, and damage models that will be implemented in the proposed Cosserat periporomechanics paradigm.

### *2.4. Micro-polar rate-dependent constitutive and damage models*

We first introduce the classical visco-plastic model and visco-elastic model cast in the framework of the Cosserat continuum theory that take into account the rate-dependency of porous geological materials [22]. We then present a micro-polar bilinear damage model and an energy damage criterion that incorporates micro-rotation energy.

#### *2.4.1. Micro-polar visco-plastic model*

The micro-polar visco-plastic model is cast using the Drucker-Prager yield surface. As an infinitesimal deformation, the strain tensor $\boldsymbol{\varepsilon}$ and the wryness tensor $\boldsymbol{\kappa}$ are additively decomposed into elastic and visco-plastic parts as

$$\varepsilon = \varepsilon^e + \varepsilon^{vp}, \tag{54}$$
$$\kappa = \kappa^e + \kappa^{vp}, \tag{55}$$

where $\boldsymbol{\varepsilon}^e$ and $\boldsymbol{\varepsilon}^{vp}$ are the elastic and visco-plastic strain tensors, respectively, and $\boldsymbol{\kappa}^e$ and $\boldsymbol{\kappa}^{vp}$are the elastic and visco-plastic wryness tensors, respectively. Given elastic strain and elastic wryness tensors, the effective stress tensor $\overline{\boldsymbol{\sigma}}$ and couple stress tensor $\boldsymbol{m}$ can be expressed as

$$\overline{\boldsymbol{\sigma}} = \lambda tr(\boldsymbol{\varepsilon}^e)1 + (\mu+\mu_c)\boldsymbol{\varepsilon}^e + (\mu-\mu_c)\boldsymbol{\varepsilon}^{e^T}, \tag{56}$$

$$\boldsymbol{m} = \alpha_1 \boldsymbol{tr}(\boldsymbol{\kappa}^e)\boldsymbol{1} + \alpha_2\boldsymbol{\kappa}^e + \alpha_3\boldsymbol{\kappa}^{e^T}, \tag{57}$$

where $\lambda$ and $\mu$ are L´ame's first elastic constant and elastic shear modulus, respectively, which can be determined from Young's modulus and Poisson's ratio, and $\mu_c$, $\boldsymbol{\alpha_1}$, $\boldsymbol{\alpha_2}$ and $\boldsymbol{\alpha_3}$ are the micro-polar parameters [28, 44, 47].

The stress tensor and couple stress tensor can be written in a vector form as form

$$\widetilde{\boldsymbol{\sigma}} = \{\overline{\sigma}_{11}, \overline{\sigma}_{22}, \overline{\sigma}_{33}, \overline{\sigma}_{12}, \overline{\sigma}_{21}, \overline{\sigma}_{13}, \overline{\sigma}_{31}, \overline{\sigma}_{23}, \overline{\sigma}_{32}\}^T \tag{58}$$

$$\widetilde{\boldsymbol{m}} = \{m_{11}/l, m_{22}/l, m_{33}/l, m_{12}/l, m_{21}/l, m_{13}/l, m_{31}/l, m_{23}/l, m_{32}/l\}^T, \tag{59}$$

where $T$ is the transpose operator and $l$ is the Cosserat length scale [28]. It follows from (60) and (61) that the mean stress $\overline{p}$ and the deviatoric stress $q$ [19, 28] can be written as

$$\overline{p} = (\widetilde{\sigma_1} + \widetilde{\sigma_2} + \widetilde{\sigma_3})/3, \tag{60}$$

$$q = \left[\frac{1}{2}\left(\widetilde{\sigma^T}\tilde{P}\widetilde{\sigma} + \widetilde{m^T}\tilde{\tilde{P}}\widetilde{m}\right)\right]^{\frac{1}{2}}, \tag{61}$$

Where

$$\tilde{P} = \begin{bmatrix} 2 & -1 & -1 & & & & & & \\ -1 & 2 & -1 & & & & & & \\ -1 & -1 & 2 & & & & & & \\ & & & 3/2 & 3/2 & & & & \\ & & & 3/2 & 3/2 & & & & \\ & & & & & 3/2 & 3/2 & & \\ & & & & & 3/2 & 3/2 & & \\ & & & & & & & 3/2 & 3/2 \\ & & & & & & & 3/2 & 3/2 \end{bmatrix} \tag{62}$$

and

$$\tilde{\tilde{P}} = 3I_{9\times 9}. \tag{63}$$

Given $\overline{p}$ and $q$, the Cosserat Drucker-Prager yield function [48] is written as

$$f = q + \sqrt{3}\mathcal{A}_{1\overline{p}} + \mathcal{A}_2. \tag{64}$$

In (66)

$$\mathcal{A}_1 = \frac{2\sin\varphi}{\sqrt{3}(3 - \sin(\varphi))}, \tag{65}$$

$$\mathcal{A}_2 = \frac{-6c\cos\varphi}{\sqrt{3}(3 - \sin(\varphi))}, \tag{66}$$

where $\varphi$ is the frictional angle. Here, we adopt a linear isotopic hardening as

$$c = c_0 + h\overline{\varepsilon}. \tag{67}$$

where $c_0$ is the initial cohesion, $h$ is the linear isotropic hardening modulus, and $\bar{\varepsilon}$ is the viscoplastic internal variable. Assuming the non-associative plasticity, the plastic flow potential [48] can be defined as

$$g = q + \sqrt{3}\mathcal{A}_{3\overline{p}} + \mathcal{A}_2, \tag{68}$$

where

$$\mathcal{A}_3 = \frac{2\sin\psi}{\sqrt{3}(3 - \sin(\psi))}, \tag{69}$$

and $\psi$ is the dilatancy angle.

Given (70), the visco-plastic strain and wryness tensors can be determined as

$$\dot{\varepsilon}^{vp} = \frac{\langle f \rangle}{\eta}\frac{\partial g}{\partial \overline{\boldsymbol{\sigma}}}, \tag{70}$$

$$\dot{\kappa}^{vp} = \frac{\langle f \rangle}{\eta}\frac{\partial g}{\partial \boldsymbol{m}}, \tag{71}$$

where $\eta$ is the viscosity of the skeleton. The internal viscoplastic variable can be defined as

$$\dot{\bar{\varepsilon}} = \left(\frac{1}{3}\dot{\varepsilon}_s^{vp}:\dot{\varepsilon}_s^{vp} + \frac{1}{3}\dot{\varepsilon}_s^{vp}:\dot{\varepsilon}_s^{vp,T} + \frac{2}{3}\dot{\kappa}^{vp}:\dot{\kappa}^{vp}\right)^{1/2}, \tag{72}$$

where $\dot{\boldsymbol{\varepsilon}}_s^{vp}$ is the deviatoric part of the visco-plastic strain tensor. By substituting (72) and (73) into (74) the internal variable can be written as

$$\dot{\bar{\varepsilon}} = \frac{\langle f \rangle}{\eta}\left[1 + \operatorname{sign}(\overline{p})\frac{\mathcal{A}_3}{\sqrt{3}}\right]. \tag{73}$$

*2.4.2. Micro-polar visco-elastic model*

In this subsection, we introduce the micro-polar visco-elastic model [49] for the simulations of cracking in Section 4. The strain and wryness tensors are additively decomposed into elastic and visco-elastic parts as

$$\varepsilon = \varepsilon^e + \varepsilon^{ve}, \tag{74}$$

$$\kappa = \kappa^e + \kappa^{ve}. \tag{75}$$

where $\boldsymbol{\varepsilon}^{ve}$ and $\boldsymbol{\kappa}^{ve}$ are visco-elastic strain and wryness tensors, respectively. Given the elastic strain and wryness tensors, (58) and (59) can be used to compute the stress and couple stress for the micropolar viscoelastic model.

Following the simple Maxwell model as shown in Figure 3, the evolution equation for the micro-polar visco-elastic model [50] can be written as

$$\dot{\varepsilon}^{ve} + \frac{\varepsilon^{ve}}{\tau_r} = \dot{\varepsilon}, \tag{76}$$

$$\dot{\kappa}^{ve} + \frac{\kappa^{ve}}{\tau_r} = \dot{\kappa}. \tag{77}$$

where $\tau_r$ is the relaxation time.

It follows from (78), (79), (58), and (59), the stress and couple stress in the integral form [50]

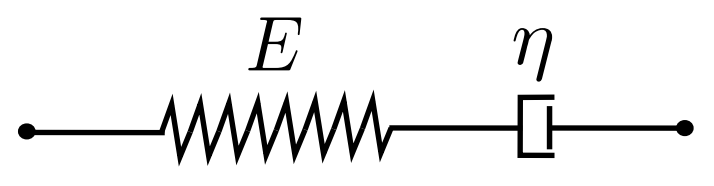

Figure 3: Schematic of Maxwell viscoelastic model.

can be written as

$$\sigma = \int_0^t \exp\left[(\tau - t)/\tau_r\right]\left[\lambda \frac{dtr(\boldsymbol{\varepsilon})1}{d\tau} + (\mu + \mu_c)\frac{d\boldsymbol{\varepsilon}}{d\tau} + (\mu - \mu_c)\frac{d\boldsymbol{\varepsilon}^T}{d\tau}\right] d\,\tau, \tag{78}$$

$$m = \int_0^t \exp\left[(\tau - t)/\tau_r\right]\left[\alpha_1 \frac{dtr(\boldsymbol{\kappa})1}{d\tau} + \alpha_2\frac{d\boldsymbol{\kappa}}{d\tau} + \alpha_3\frac{d\boldsymbol{\kappa}^T}{d\tau}\right] d\,\tau. \tag{79}$$

For the micropload visco-elastic correspondence model, the effective force and moment states [15] can be written as

$$\underline{\boldsymbol{T}} = \underline{t}\left(\frac{\underline{\widetilde{U}}}{|\underline{\widetilde{U}}|}\right), \tag{79}$$

$$\underline{\boldsymbol{M}} = \underline{\boldsymbol{m}}\left(\frac{\underline{\boldsymbol{\Omega}}}{|\underline{\boldsymbol{\Omega}}|}\right), \tag{80}$$

where $\underline{t}$ and $\underline{\boldsymbol{m}}$ are the scalar effective force and moment states, respectively, and $\underline{\widetilde{\boldsymbol{U}}}/|\underline{\widetilde{\boldsymbol{U}}}|$ and $\underline{\boldsymbol{\Omega}}/|\underline{\boldsymbol{\Omega}}|$ denote the directions of the effective force and moment vector states, respectively. For a bond-based visco-elastic micropolar material model, the effective force state can be decomposed into the part parallel to the bond $\underline{t}_1$ and part perpendicular to the bond $\underline{t}_2$ as

$$\underline{t}_1 = \frac{1}{2} f_1(u_1), \tag{84}$$

$$\underline{t}_2 = \frac{1}{2} f_2(\boldsymbol{u}_2), \tag{85}$$

where $\boldsymbol{u}_1$ and $\boldsymbol{u}_2$ are the displacements in the axial and the normal directions of the bond, respectively. Figure 4 presents a schematic of the decomposition of deformation states and force states for the micro-polar visco-elasticity model. For a bond-based micro-polar visco-elastic model,

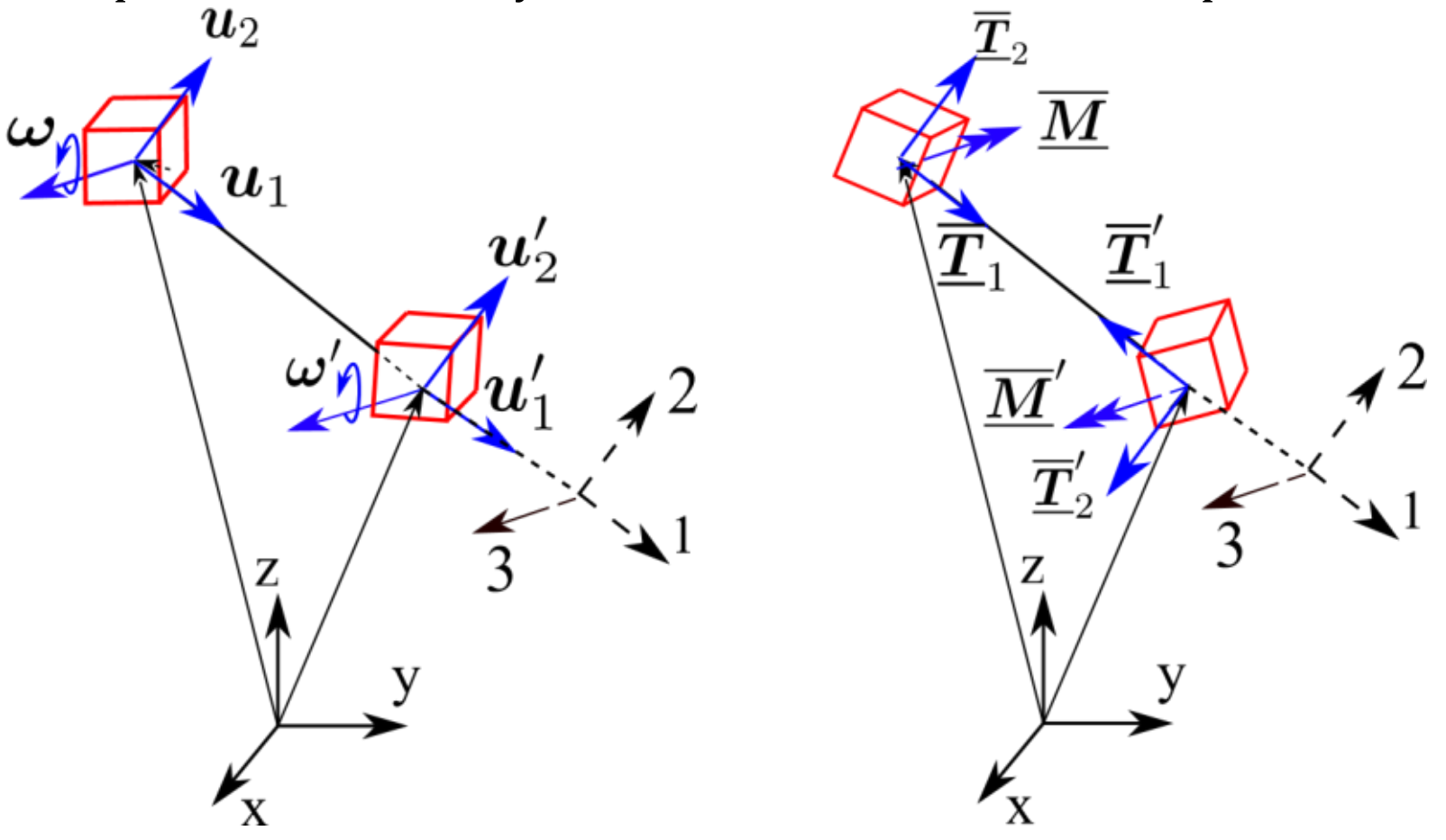

Figure 4: Schematic of the decomposition of deformation states and force states for the micro-polar visco-elasticity model.

the bond stretch rate $\dot{s}_1$ and shear deformation rate $\dot{s}_2$ [51] can be defined as

$$\dot{s}_1 = \frac{|\Delta\dot{u}_1|}{|\underline{\xi}|}, \tag{86}$$

$$\dot{s}_2 = \frac{|\Delta\dot{u}_2 - \dot{\underline{\Omega}}\underline{\xi}|}{|\underline{\xi}|}. \tag{87}$$

Figure 5 sketches the concept for a special case of the visco-elastic material model. Following

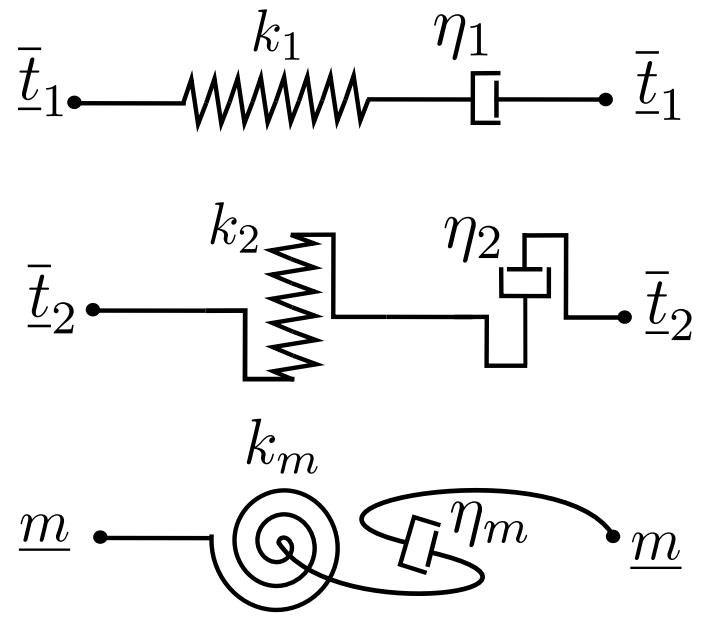

Figure 5: Schematic of a simple viscoelastic micropolar model.

the linear micro-polar elasticity, it is assumed that the total bond stretch $s_1$, shear deformation $s_2$, and the relative micro rotation $\Delta\omega$ (magnitude) can be decomposed into elastic and visco-elastic parts as

$$s_1 = s_1^e + s_1^{ve} \tag{88}$$

$$s_2 = s_2^e + s_2^{ve}, \tag{89}$$

$$\Delta\omega = \Delta\omega^e + \Delta\omega^{ve}. \tag{90}$$

where $s_1^e$ and $s_1^{ve}$ are elastic and visco-elastic stretches, respectively, $s_2^e$ and $s_2^{ve}$ are elastic and visco-elastic shear deformations, respectively, and $\Delta\omega^e$ and $\Delta\omega^{ve}$ are the elastic and visco-elastic relative rotations (magnitude), respectively. The scalar axial and normal force and scalar moment states can be defined as

$$\underline{\bar{t}}_1 = k_1 s_1^e, \tag{91}$$
$$\underline{\bar{t}}_2 = k_2 s_2^e, \tag{92}$$
$$\underline{\mathrm{m}} = \mathrm{k_m} \Delta\omega^{\mathrm{e}}, \tag{93}$$

where $k_1$, $k_2$, and $k_m$ are the material constants [52, 53]. The evolution equations for the bondbased micro-polar visco-elastic model can be written as

$$\dot{s}^{ve} + \frac{s^{ve}}{\tau_1} = \dot{s}, \tag{94}$$
$$\dot{s}_2^{ve} + \frac{s_2^{ve}}{\tau_2} = \dot{s}_2 \tag{95}$$

$$\Delta\dot{\omega}^{ve} + \frac{\Delta\omega^{ve}}{\tau_m} = \Delta\dot{\omega}. \tag{96}$$

where $\tau_1$, $\tau_2$ and $\tau_m$ are the relaxation time for the axial force, normal force, and moment, respectively. For simplicity, in this study it is assumed that $\tau_1 = \tau_2 = \tau_m = \tau_r$. This constitutive model in the integration form can be written as

$$\underline{\bar{t}}_1 = \int_0^t k_1 \exp[(\tau - t)/\tau_r] \frac{ds_1}{d\tau} d\tau, \tag{97}$$

$$\underline{\bar{t}}_2 = \int_0^t k_2 \exp[(\tau - t)/\tau_r] \frac{ds_2}{d\tau} d\tau, \tag{98}$$

$$\underline{m} = \int_0^t k_m \exp[(\tau - t)/\tau_r] \frac{d\Delta\boldsymbol{\omega}}{d\tau} d\tau. \tag{99}$$

In what follows, we introduce a bilinear micro-polar damage model.

*2.4.3. Ordinary micro-polar damage model*

In this study, we present an ordinary micro-polar bilinear damage model following [54] for modeling the softening behavior of quasi brittle porous media. Figure 6 presents a schematic of the bilinear micro-polar damage model.

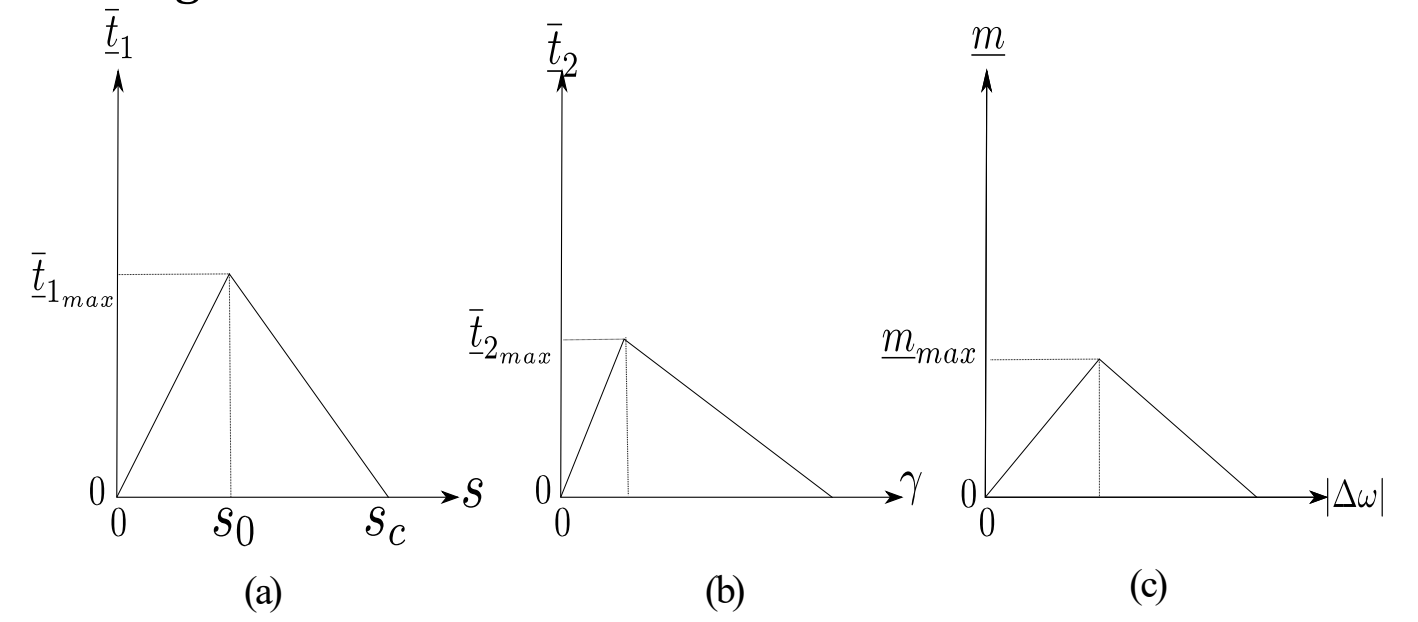

Figure 6: Schematic of the bilinear micro-polar damage model: (a) axial force, (b) normal force, and (c) moment.

It follows from Figure 6 that the bilinear damage model [52, 54] can be written as

$$\underline{\bar{t}}_1 = \begin{cases} k_1 s_1^e & if\ s_1 < s_0, \\ 0 & if\ s_1 \geq s_c, \end{cases} \tag{100}$$

$$\underline{\bar{t}}_2 = \begin{cases} k_2 s_2^e & if\ s_1 < s_0, \\ 0 & if\ s_1 \geq s_c, \end{cases} \tag{101}$$

$$\underline{m} = \begin{cases} k_{\mathrm{m}} \Delta\omega^e & if\ s_1 < s_0, \\ 0 & if\ s_1 \geq s_c, \end{cases} \tag{102}$$

The influence function $\underline{\omega}$ can be calculated as

$$\underline{\omega} = \begin{cases} k_1 s_1^e & if\ s_1 < s_0, \\ 0 & if\ s_1 \geq s_c, \end{cases} \tag{103}$$

where $s_0$ and $s_c$ are the bond stretch at the peak load and the critical bond stretch, respectively. The local damage parameter for a material point [39] can be defined as

$$D = 1 - \frac{\int_{\mathscr{H}} \underline{\omega} dV'}{\int_{\mathscr{H}} dV'}. \tag{104}$$

*2.4.4. Energy-based bond breakage criterion*

In this study, we also formulate an energy-based criteria to detect the bond breakage. In this case, the bond-breakage criterion depends on the deformation energy in a bond. The effective force state and the moment state (the energy conjugates of the composite state and relative rotation state, respectively) are used to determine the deformation energy [55]. Therefore, the energy density in bond $\xi$ can be written as

$$\mathcal{W} = \int_0^{t_f} (\underline{\mathscr{T}} - \underline{\mathscr{T}}')\dot{\underline{\tilde{U}}} dt + \int_0^{t_f} (\underline{\mathscr{M}} - \underline{\mathscr{M}}')\dot{\underline{\Omega}} dt. \tag{105}$$

where $t_f$ is the final loading time. In this case, the influence function for a bond $\xi$ can be defined as

$$\underline{\omega} = \begin{cases} 1 & \text{if } \mathcal{W} < \mathcal{W}_{cr} \\ 0 & \text{if } \mathcal{W} \geq \mathcal{W}_{cr}. \end{cases} \tag{106}$$

It is noted that the critical energy density for bond breakage can be calculated from the critical energy release rate as

$$\mathcal{W}_{cr} = \frac{4\mathcal{G}_{cr}}{\pi\delta^4}. \tag{107}$$

In the linear elastic fracture mechanics for mode I fracture the critical energy release rate reads

$$\mathcal{G}_{cr} = \mathcal{K}_I^2 \frac{1-\nu^2}{E}, \tag{108}$$

where $K_I$ is the fracture toughness of mode I crack. It follows from (106) that the local damage parameter at a material point can be determined from (104).

In the following section, we present the numerical implementation of the proposed Cosserat periporomechanics paradigm for porous media under dry conditions.

## 3. Numerical implementation

### *3.1. Discretization in space*

The motion equation (54) and the moment equation (55) are discretized in space by the total Lagrangian meshfree scheme [5–7]. In this method, a porous continuum material is discretized into a finite number of mixed material points (i.e., solid skeleton and pore water). Under the dry condition, each material point has two kinds of degrees of freedom (i.e., displacement and micro-rotation). The uniform grid is used to spatially discretize the problem domain in which all material points have the same dimensions. The spatial discretization forms of (54) and (55) are written as

$$\mathbf{0} = \mathscr{A}_{i=1}^{\mathcal{P}}(\overline{\mathscr{M}}_i \ddot{\mathbf{u}}_i - \overline{\mathscr{T}}_i + \mathcal{F}_i), \tag{109}$$

$$\mathbf{0} = \mathscr{A}_{i=1}^{\mathcal{P}}(\mathscr{I}_i \ddot{\boldsymbol{\omega}}_i - \mathscr{M}_i + \widetilde{\mathscr{M}}_i + \mathcal{L}_i), \tag{110}$$

where $\mathcal{A}$ is a global linear assembly operator [4, 40], $\overline{\mathcal{M}}_i$ is the mass matrix at material point $i$, $\overline{\mathcal{T}}_i$ is the vector of effective force, $\mathcal{F}_i$ is the vector gravity force [7], $\mathcal{M}_i$ is the vector of rotational moment, $\widetilde{\mathcal{M}}_i$ is the moment by the effective force state, and $\mathcal{L}_i$ is the body couple vector. At material point $i$, these five vectors can be written as

$$\overline{\boldsymbol{\mathscr{M}}}_i = \rho_s \phi_i \mathscr{V}_i \mathbf{1}, \tag{111}$$

$$\overline{\boldsymbol{\mathscr{T}}}_i = \sum_{j=1}^{\mathcal{N}_i} \left( \overline{\underline{\mathscr{T}}}^s_{(ij)} - \overline{\underline{\mathscr{T}}}^s_{(ji)} \right) \mathscr{V}_j \mathscr{V}_i, \tag{112}$$

$$\boldsymbol{\mathscr{I}}_i = \mathscr{I}_i \mathscr{V}_i \mathbf{1}, \tag{113}$$

$$\boldsymbol{\mathscr{M}}_i = \sum_{j=1}^{\mathcal{N}_i} \left( \underline{\mathscr{M}}^s_{(ij)} - \underline{\mathscr{M}}^s_{(ji)} \right) \mathscr{V}_j \mathscr{V}_i, \tag{114}$$

$$\widetilde{\boldsymbol{\mathscr{M}}}_i = \sum_{j=1}^{\mathcal{N}_i} \left[ \frac{1}{2} \underline{\mathbf{Y}}_{(ij)} \left( \underline{\mathscr{T}}^s_{(ij)} - \underline{\mathscr{T}}^s_{(ji)} \right) \right] \mathscr{V}_j \mathscr{V}_i, \tag{115}$$

where $\mathcal{V}_i$ and $\mathcal{V}_j$ are the volumes of material points *i* and *j*, respectively, in the reference configuration. In (112), (114) and (115), the effective force state and the moment state are written as

$$\underline{\mathscr{T}}^{s}_{(ij)} = \underline{\omega}_{(ij)}\overline{\boldsymbol{\sigma}}_{(i)}\mathscr{K}^{-1}_{(i)}\underline{\boldsymbol{\xi}}_{(ij)} + \underline{\alpha\mathscr{R}}_{1(ij)}, \tag{116}$$

$$\underline{\mathscr{T}}^{s}_{(ji)} = \underline{\omega}_{(ji)}\overline{\boldsymbol{\sigma}}_{(j)}\mathscr{K}^{-1}_{(j)}\underline{\boldsymbol{\xi}}_{(ji)} + \underline{\alpha\mathscr{R}}_{1(ji)}, \tag{117}$$

$$\underline{\mathscr{M}}^{s}_{(ij)} = \underline{\omega}_{(ij)}\boldsymbol{m}_{(i)}\mathscr{K}^{-1}_{(i)}\underline{\boldsymbol{\xi}}_{(ij)} + \underline{\beta\mathscr{R}}_{2(ij)}, \tag{118}$$

$$\underline{\mathscr{M}}^{s}_{(ji)} = \underline{\omega}_{(ji)}\boldsymbol{m}_{(j)}\mathscr{K}^{-1}_{(j)}\underline{\boldsymbol{\xi}}_{(ji)} + \underline{\beta\mathscr{R}}_{2(ji)}. \tag{119}$$

The micro-polar strain tensor and the wryness tensor can be written as

$$\overline{\varepsilon}_{(i)} = \left[\sum_{j=1}^{\mathcal{N}_i} \underline{\omega}_{(ij)}(\underline{\overline{U}}_{(ij)} \otimes \underline{\boldsymbol{\xi}}_{(ij)})\mathscr{V}_j\right]\mathscr{K}^{-1}_{(i)}, \tag{120}$$

$$\boldsymbol{\kappa}_{(i)} = \left[\sum_{j=1}^{\mathcal{N}_i} \underline{\omega}_{(ij)}(\underline{\boldsymbol{\Omega}}_{(ij)} \otimes \underline{\boldsymbol{\xi}}_{(ij)})\mathscr{V}_j\right]\mathscr{K}^{-1}_{(i)}. \tag{121}$$

Given $\overline{\boldsymbol{\varepsilon}}_{(i)}$ and $\boldsymbol{\kappa}_{(i)}$, classical constitutive models can be used to compute $\overline{\boldsymbol{\sigma}}_{(i)}$ and $\boldsymbol{m}_{(i)}$ as described in Section 2.4. In what follows we present the discretization in time through an explicit Newmark scheme [11, 40].

### *3.2. Discretization in time*

The Newmark scheme [40] is applied to integrate the equations of motion and moment in time. Let $\boldsymbol{u}_n, \dot{\boldsymbol{u}}_n$, and $\ddot{\boldsymbol{u}}_n$ be the displacement, velocity, and acceleration vectors at time step *n*. The predictors of displacement and velocity in a general Newmark scheme read

$$\dot{\tilde{u}}_{n+1} = \dot{u}_n + (1-\beta_1)\Delta\ddot{u}_n, \tag{122}$$

$$\tilde{u}_{n+1} = u_n + \Delta t\dot{u}_n + \frac{\Delta t^2}{2}(1-2\beta_2)\,\ddot{u}_n, \tag{123}$$

$$\dot{\tilde{\omega}}_{n+1} = \dot{\omega}_n + (1-\beta_1)\Delta\ddot{\omega}_n, \tag{124}$$

$$\tilde{\omega}_{n+1} = \omega_n + \Delta t\dot{\omega}_n + \frac{\Delta t^2}{2}(1-2\beta_2)\,\ddot{\omega}_n, \tag{125}$$

where $\beta_2$ and $\beta_1$ are numerical integration parameters. Given (123) and (122), the accelerations $\ddot{u}_{n+1}$ and $\omega$ are determined by the recursion relation

$$\ddot{u}_{n+1} = \boldsymbol{\mathscr{M}}^{-1}_{n+1}(\boldsymbol{\mathcal{F}}_{n+1} - \widetilde{\widetilde{\boldsymbol{\mathscr{T}}}}_{n+1}), \tag{126}$$

$$\ddot{\omega}_{n+1} = \boldsymbol{\mathscr{I}}^{-1}_{n+1}(\boldsymbol{\mathcal{L}}_{n+1} - \widetilde{\widetilde{\boldsymbol{\mathscr{M}}}}_{n+1} + \widetilde{\widetilde{\boldsymbol{\mathscr{M}}}}_{n+1}), \tag{127}$$

Where $\widetilde{\widetilde{\boldsymbol{\mathscr{T}}}}_{n+1}$ $\widetilde{\mathscr{M}}_{n+1}$, and $\widetilde{\widetilde{\mathscr{T}}}_{n+1}$ are determined from (123) and (125) and the local constitutive models. From (126) and (127), the displacement, velocity, rotation, and rotation rate at time step *n* + 1 can be updated as

$$\dot{u}_{n+1} = \dot{\tilde{u}}_{n+1} + \beta_1\Delta t\ddot{u}_{n+1}, \tag{128}$$

$$u_{n+1} = \tilde{u}_{n+1} + \beta_2\Delta t^2\ddot{u}_{n+1}. \tag{129}$$

$$\dot{\omega}_{n+1} = \dot{\tilde{\omega}}_{n+1} + \beta_1\Delta t\ddot{\omega}_{n+1}, \tag{130}$$

$$\omega_{n+1} = \tilde{\omega}_{n+1} + \beta_2\Delta t^2\ddot{\omega}_{n+1}. \tag{131}$$

In this study, we adopt the explicit central difference solution scheme [11, 40] in which $\beta_1 = 1/2$ and $\beta_2 = 0$. We note that the explicit method is efficient and robust to model dynamic problems

[39].

The energy balance check is used to ensure numerical stability of the algorithm in time [41]. We define the internal energy, external energy, and kinetic energy of the system at time step $n$+1 as

$$\begin{aligned}\mathscr{W}_{\mathrm{int},n+1} =& \mathscr{W}_{\mathrm{int},n} + \frac{\Delta t}{2}\left(\dot{u}_n + \frac{\Delta t}{2}\ddot{u}_n\right)\left(\overline{\mathscr{T}}_n + \overline{\mathscr{T}}_{n+1}\right) + \frac{\Delta t}{2}\left(\dot{\omega}_n + \frac{\Delta t}{2}\ddot{\omega}_n\right)\left[\left(\boldsymbol{\mathscr{M}}_n - \widetilde{\boldsymbol{\mathscr{M}}}_n\right)\right.\\ &\left.+\left(\boldsymbol{\mathscr{M}}_{n+1} - \widetilde{\boldsymbol{\mathscr{M}}}_{n+1}\right)\right], \qquad (132)\end{aligned}$$

$$\mathscr{W}_{\mathrm{ext},n+1} = \mathscr{W}_{\mathrm{ext},n} + \frac{\Delta t}{2}\left(\dot{u}_n + \frac{\Delta t}{2}\ddot{u}_n\right)\left(\boldsymbol{\mathcal{F}}_n + \boldsymbol{\mathcal{F}}_{n+1}\right) + \frac{\Delta t}{2}\left(\dot{\omega}_n + \frac{\Delta t}{2}\ddot{\omega}_n\right)\left(\boldsymbol{\mathcal{L}}_n + \boldsymbol{\mathcal{L}}_{n+1}\right), \qquad (133)$$

$$\mathscr{W}_{\mathrm{kin},n+1} = \frac{1}{2}\dot{u}_{n+1}\overline{\boldsymbol{\mathscr{M}}}_{n+1}\dot{u}_{n+1} + \frac{1}{2}\dot{\omega}_{n+1}\boldsymbol{\mathscr{I}}_{n+1}\dot{\omega}_{n+1}. \qquad (134)$$

where $\hat{\varepsilon}$ is a small tolerance on the order of $10^{-2}$ [41].

$$\left|\mathscr{W}_{\mathrm{int},n+1} + \mathscr{W}_{\mathrm{kin},n+1} - \mathscr{W}_{\mathrm{ext},n+1}\right| \leqslant \hat{\varepsilon}\,\max\left(\mathscr{W}_{\mathrm{int},n+1}, \mathscr{W}_{\mathrm{kin},n+1}, \mathscr{W}_{\mathrm{ext},n+1}\right), \qquad (135)$$

Then it follows from the energy conservation criterion that

For the numerical implementation algorithms for the micro-polar visco-plastic and visco-elastic models, we refer to the celebrated literature on the subject (e.g., [23, 28, 56]).

## 4. Numerical examples

This section presents four numerical examples to show the efficacy and robustness of the proposed Cosserat periporomechanics model in simulating the dynamic shear banding and fracturing in porous media. Example 1 deals with the single shear banding in a sand specimen. Example 2 concerns the cracking of a quasi-brittle porous medium under a three-point bending test. Example 3 deals with the conjugate shear banding under the dynamic loading condition. Example 4 concerns crack branching in porous media under high loading rates.

### *4.1. Example 1: Single shear banding*

This example simulates a single shear banding in granular materials under non-symmetrical loading conditions. The numerical test is inspired by the experimental testing of shear banding of granular materials in [18]. In this example, we study the influence of spatial discretizations, internal length scales, and loading rates on the shear band formation. Figure 7 shows the boundary conditions and the loads. The specimen's length, width, and thickness are 140 mm, 40 mm, and 80 mm, respectively. The constant confining pressure is 0.2 MPa. The total vertical displacement imposed on the top and bottom boundaries is $u_y$ = 3.5 mm. The lateral displacement $u_x$ = 1 mm is applied on the top boundary to induce the single shear band in the specimen. The rates for the two displacement loads are $\dot{u}_y$ = 0.1 m/s and $\dot{u}_x$ = 0.029 m/s, respectively. Following [18], the input material parameters are: the density $\rho^s$ = 1650 kg/m$^3$, Young's modulus $E$ = 50.4 MPa, Poisson ratio $\nu$ = 0.4, Cosserat shear modulus $\mu_c = 2\mu$, initial volume fraction $\varphi_0$ = 0.65, and Cosserat length scale $l$ = 1 mm. For the visco-plastic model, initial cohesion $c_0$ = 0.13 MPa, frictional angle $\phi = 42^\circ$, dilation angle $\psi = 33^\circ$, softening modulus $h$ = −1.5 MPa, and viscosity $\eta$ = 0.003 MPa.s. Two stabilization parameters $G_1$ = 0.0175 and $G_2$ = 0.0017 [5] are used for stabilization on the localization and post-localization stage, respectively. The specimen is discretized into 28×98 material points with the uniform grid size $\Delta x$ = 1.43 mm and the horizon size $\delta = 2.05\Delta x$. The simulation time is $t = 3 \times 10^4\,\mu$s, and the time increment $\Delta t = 3\,\mu$s.

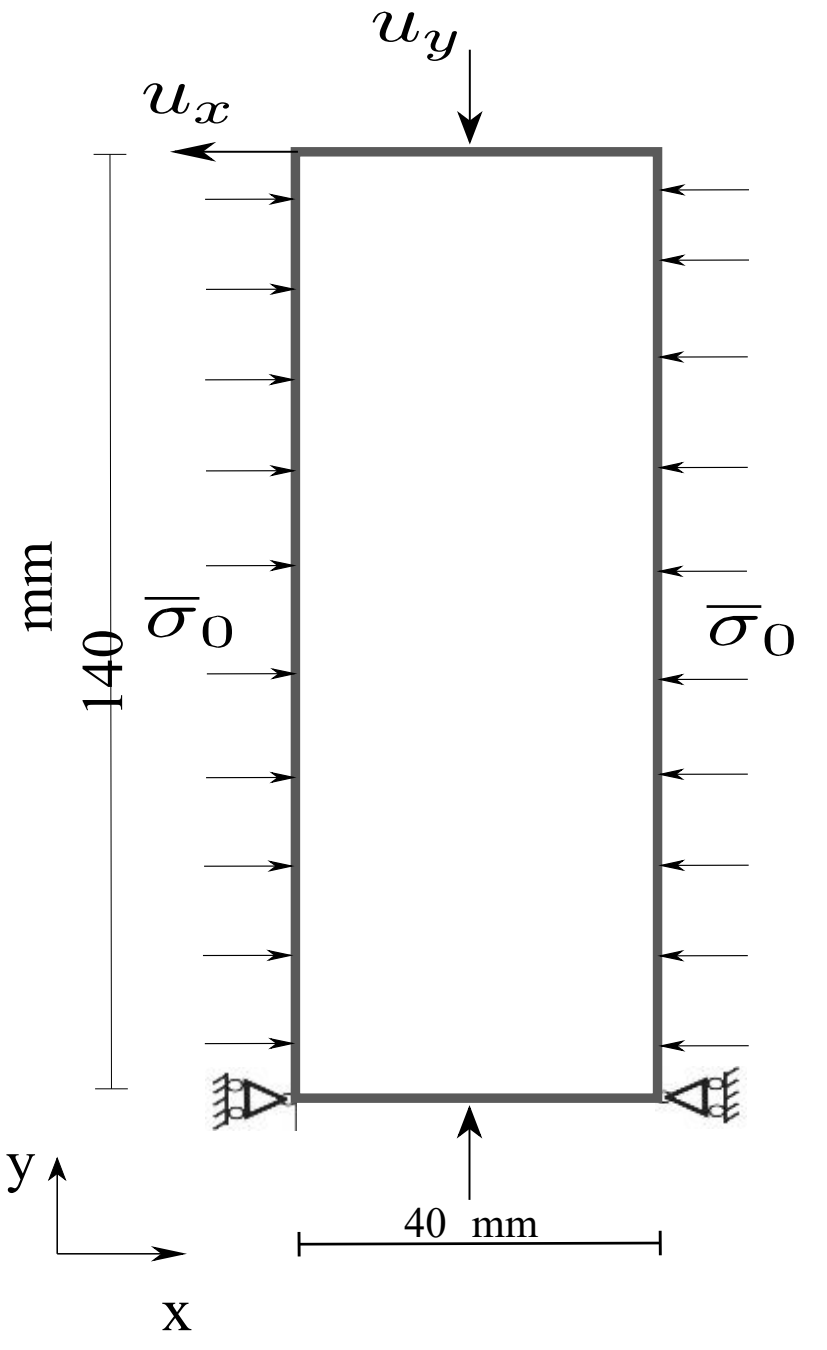

Figure 7: Problem setup for example 1.

Figure 8 plots the reaction force versus the vertical displacement on the top boundary (i.e., loading curve), which shows a softening stage following the peak load. Three points on the loading curve are selected in the softening stage to demonstrate the development of the single shear band in what follows. Figures 9, 10, 11 present the snapshots of equivalent plastic shear strain, plastic volume strain, and micro-rotation at the three loading stages, respectively. From the results in Figure 9 and 10, it follows that the shear band develops gradually with an inclination angle of $54^\circ$. The inclination angle of the shear band is comparable with the biaxial test reported in [18]. Figure 10 shows that the plastic volume strain is positive within the shear band, which is typical for a medium-dense granular material. Figure 11 plots the contour of the micro-rotation at the three loading stages. The results show the micro-rotation mainly occurs with the material points within the shear band, which is consistent with the experimental results of shear banding in granular materials [18]. In what follows, we study the influence of spatial discretizations, Cosserat length scales, and loading rates.

*4.1.1. Influence of spatial discretization*

We study the influence of spatial discretization on the results under the other same conditions. For this purpose, we consider two spatial discretization schemes, i.e., 28×98 points with $\Delta x$ = 1.43 mm and 42 × 147 points with $\Delta x$ = 0.95 mm. The same horizon $\delta$ = 2.86 mm is chosen for both cases. All other parameters remain unchanged. Figure 12 shows the vertical loading curves of the simulations with the two spatial discretization schemes. The two loading curves are the same until the peak value. In the softening stage, spatial discretization slightly influences the loading curve. Figures 13, 14 and 15 compare the contours of equivalent plastic shear strain, plastic volumetric strain, and micro rotation at the final loading stages of the two simulations, respectively. The results in Figures 13, 15 and 14 demonstrate that the spatial discretization scheme given the same horizon has a mild influence on the shear banding formation.

*4.1.2. Influence of Cosserat length scales*

We study the influence of Cosserat length scales on the shear band formation. Three different Cosserat length scales considered are $l_1$ = 1 mm $l_2$ = 0.5 mm, and $l_3$ = 0.25 mm. The same

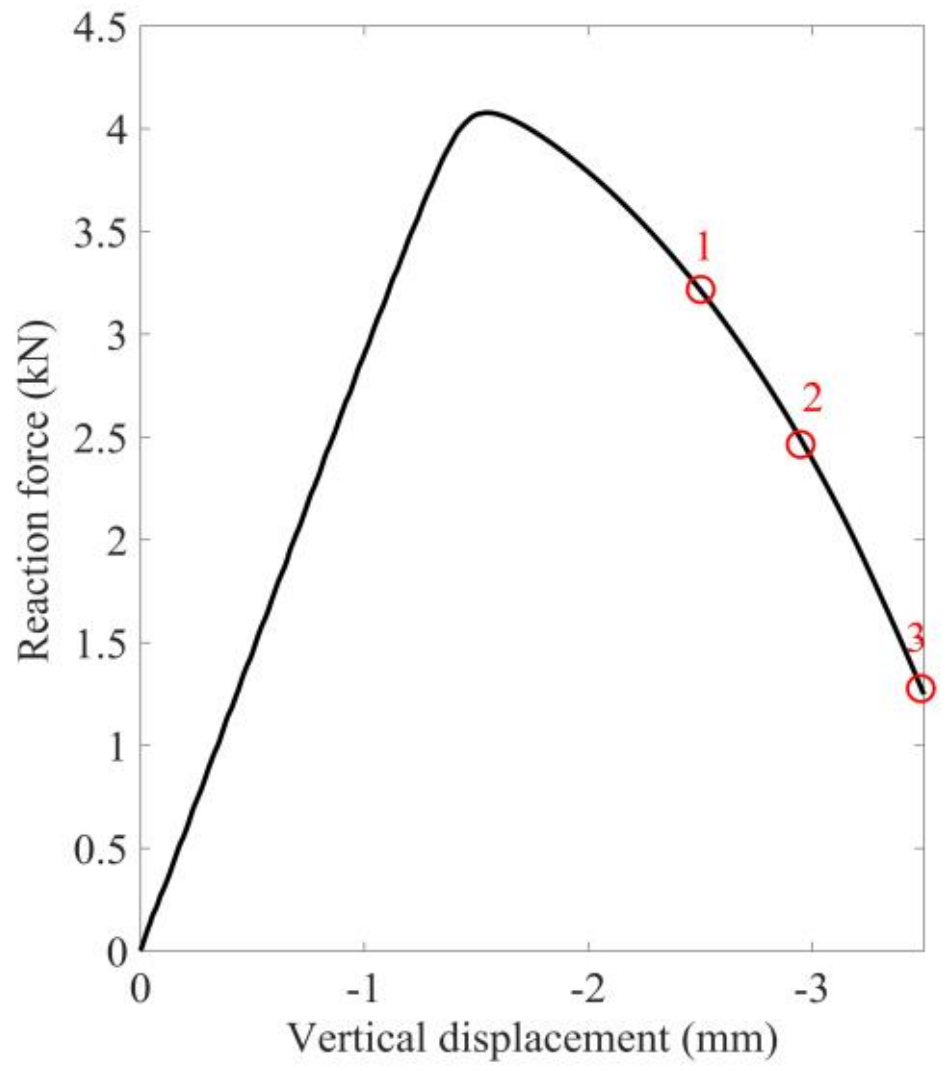

Figure 8: Reaction force versus vertical displacement on the top boundary. Note: The displacements at points 1, 2, 3 are $u_{y,1}$ = 2.5 mm, $u_{y,2}$ = 3 mm and $u_{y,3}$ = 3.5 mm, respectively.

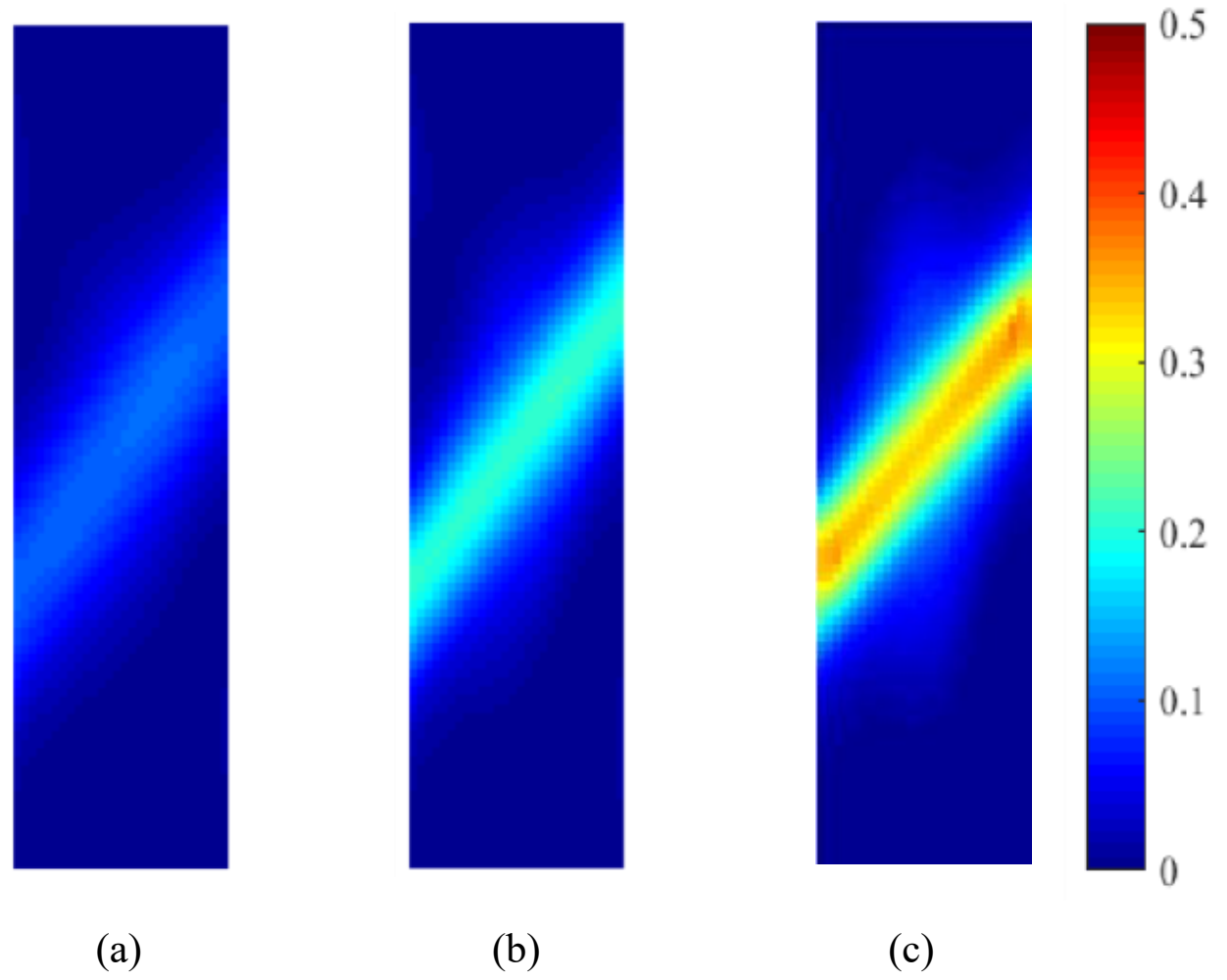

Figure 9: Contours of the equivalent plastic shear strain at three loading stages: (a) $u_{y,1}$ = 2.5 mm, (b) $u_{y,2}$ = 3 mm, and (c) $u_{y,3}$ = 3.5 mm.

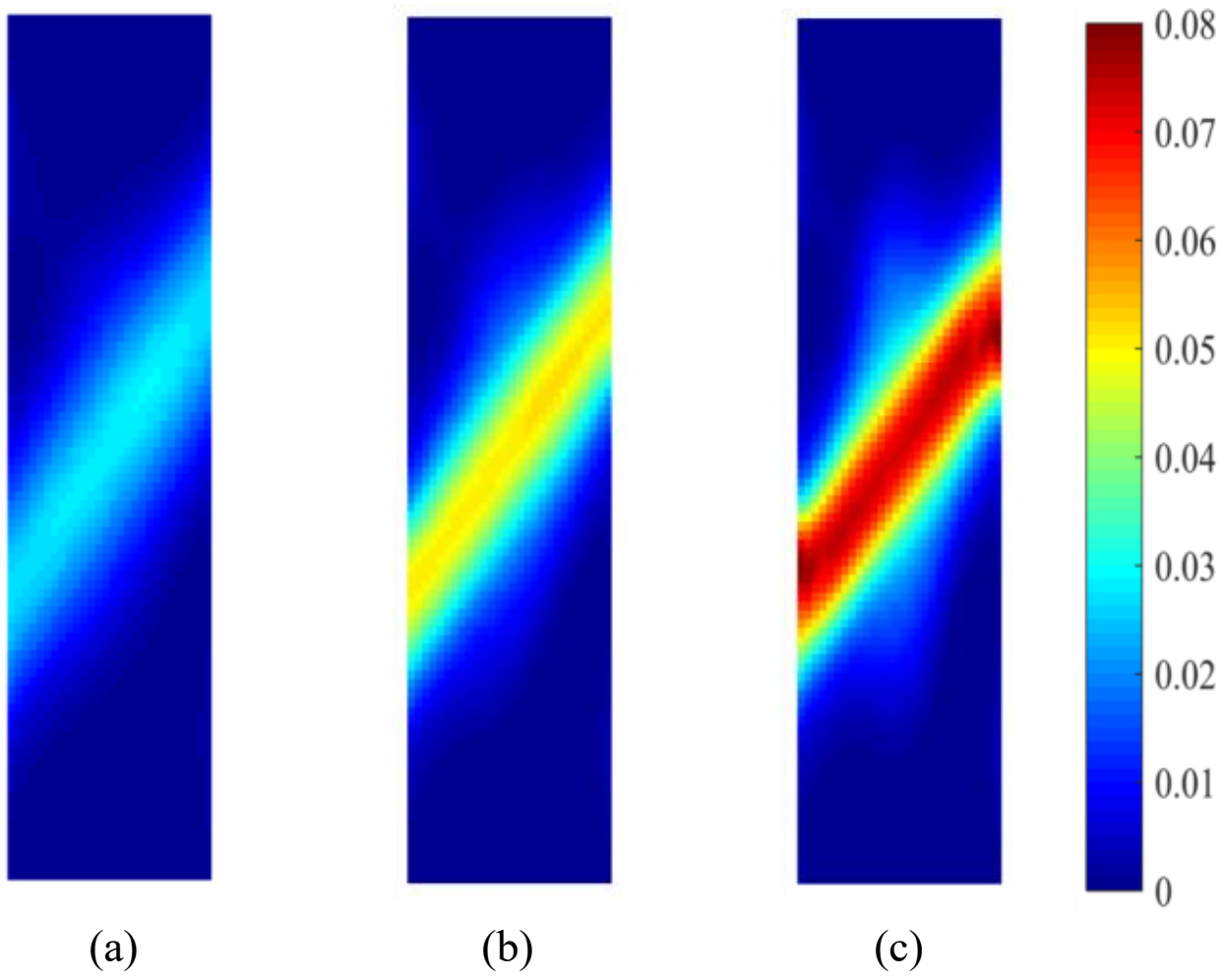

Figure 10: Contours of the plastic volume strain at at three loading stages: (a) $u_{y,1}$ = 2.5 mm, (b) $u_{y,2}$ = 3 mm, and (c) $u_{y,3}$ = 3.5 mm.

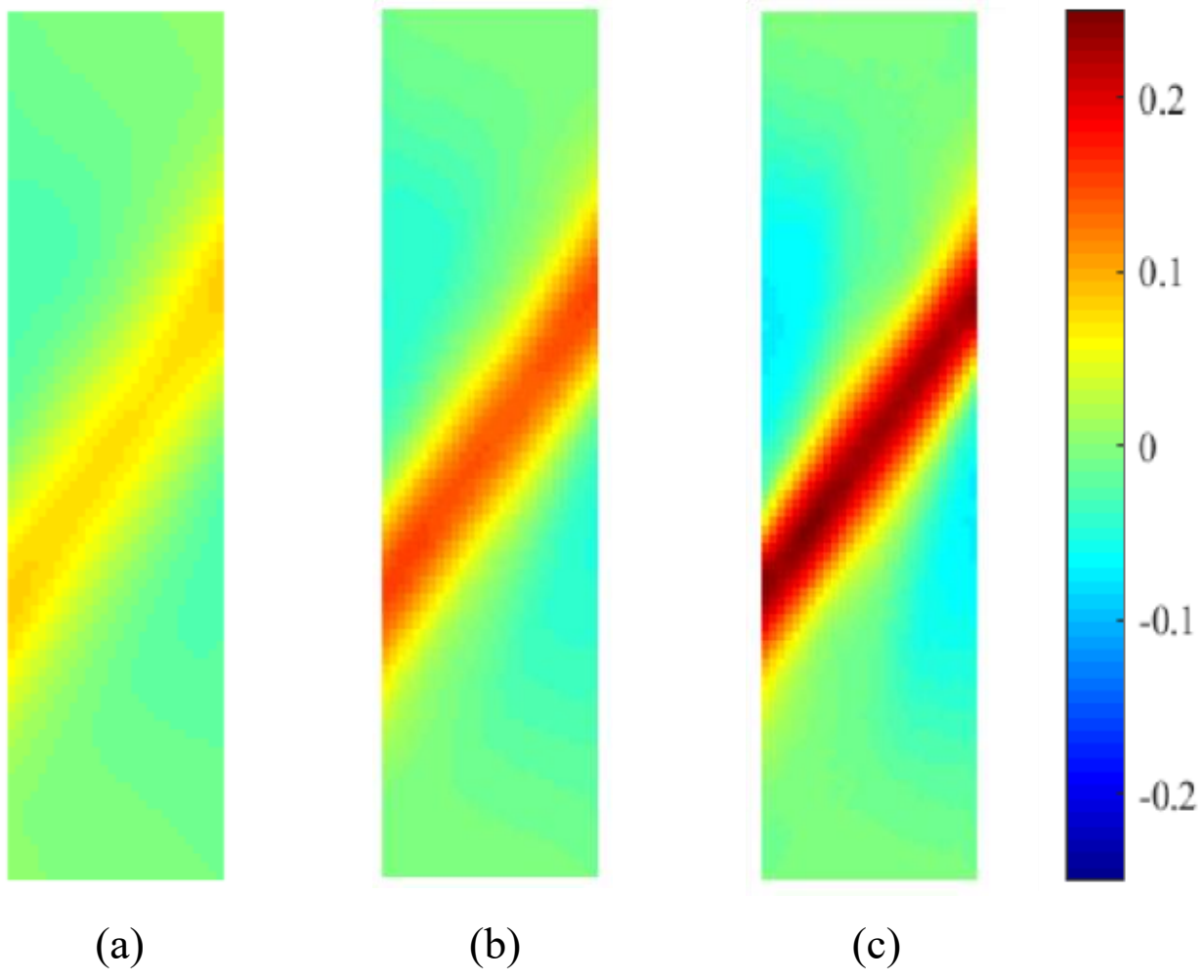

Figure 11: Contours of the micro rotation (rad) at three loading stages: (a) $u_{y,1}$ = 2.5 mm, (b) $u_{y,2}$ = 3 mm, and (c) $u_{y,3}$ = 3.5 mm.

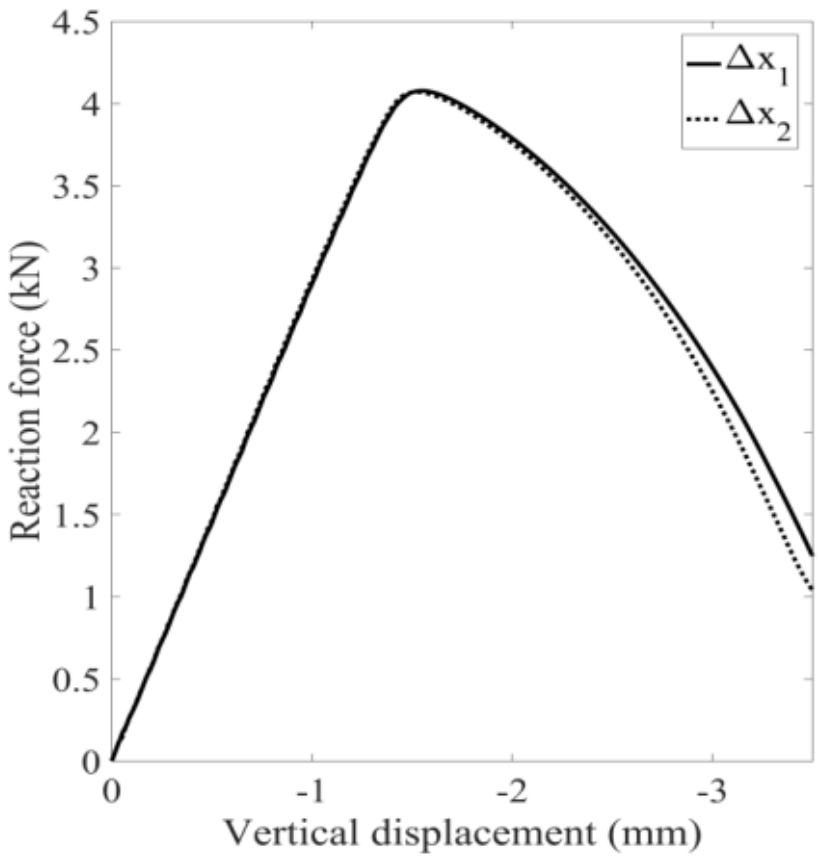

Figure 12: Comparison of vertical loading curves from the simulations with the two spatial discretization schemes.

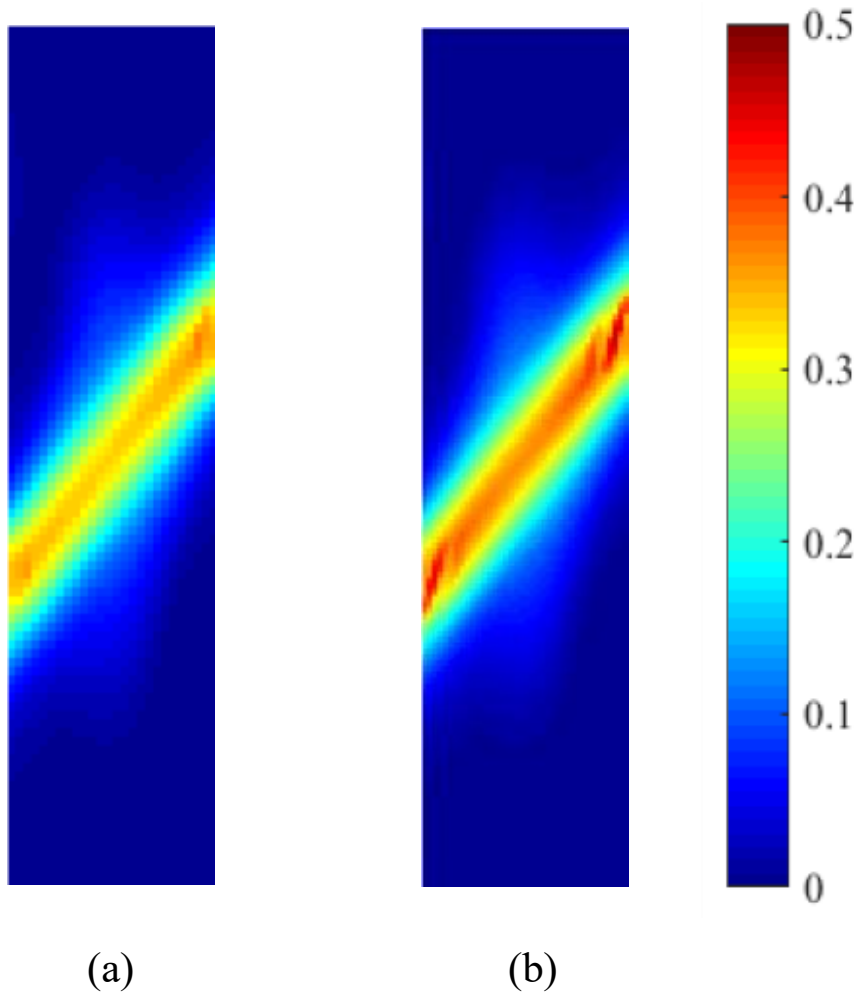

Figure 13: Contours of the equivalent plastic shear strain from the simulations with the two spatial discretization schemes.

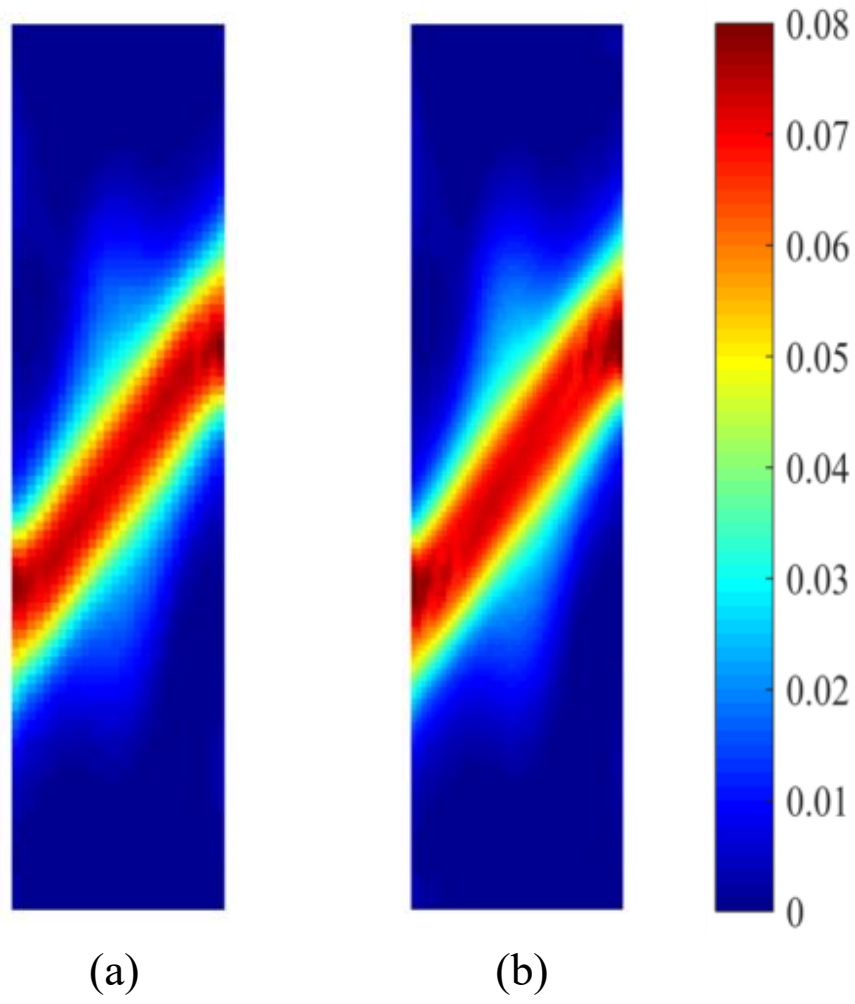

Figure 14: Contours of the plastic volume strain from the simulations with the two spatial discretization schemes.

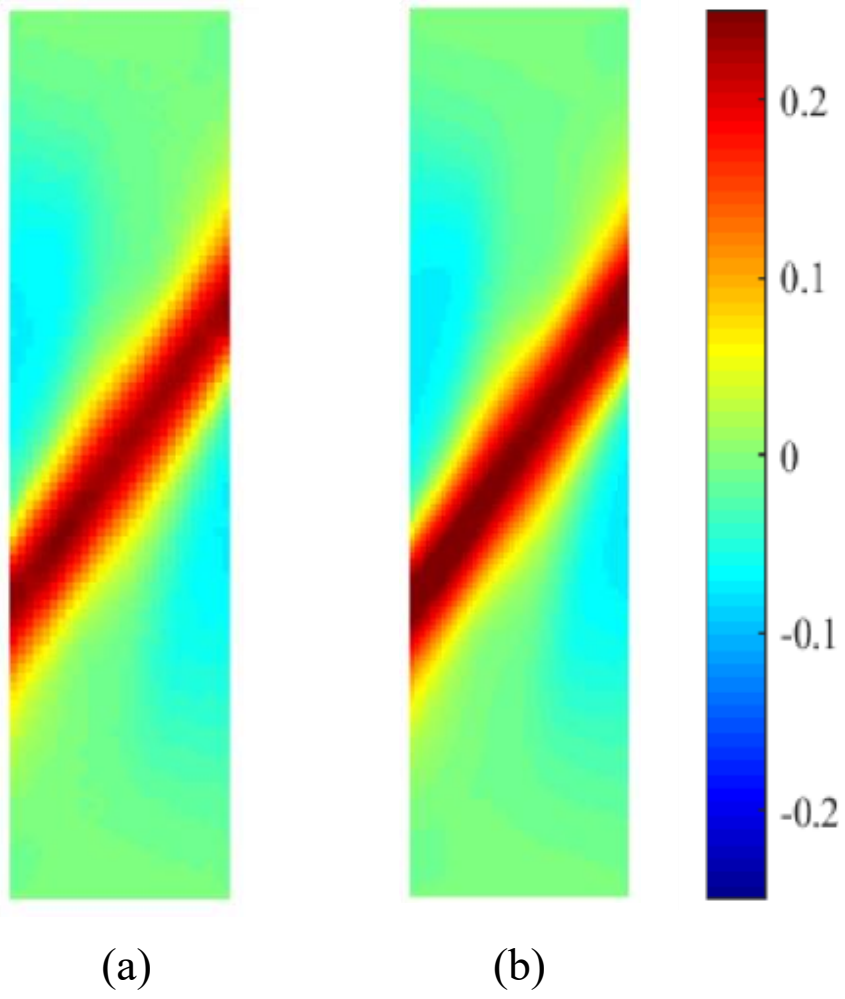

Figure 15: Contours of the micro rotation (rad) from the simulations with the two spatial discretization schemes.

horizon $\delta$ = 2.86 mm is adopted for the three cases. All other material parameters remain the same. Figure 16 plots the loading curve of mentioned models. Figure 16 shows a steeper post-peak response for the model with the smaller internal length scale. Furthermore, the contours of plastic shear strain, plastic volumetric strain, and micro rotation for the mentioned cases at the final stage are shown in Figures17, 18, and 19, respectively. As shown in Figure 17, the width of the shear band depends on the Cosserat length scale, which agrees with the experiment.

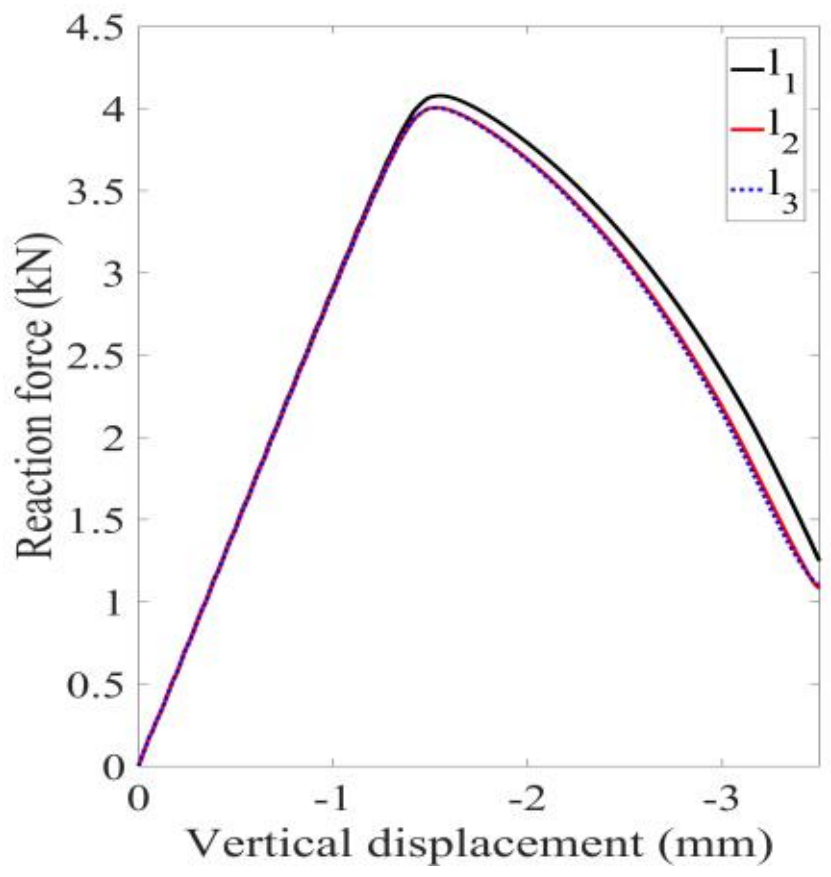

Figure 16: Comparison of the vertical loading curves from the simulations with $l_1 = 1$ mm, $l_2 = 0.5$ mm, and $l_3 = 0.25$ mm.

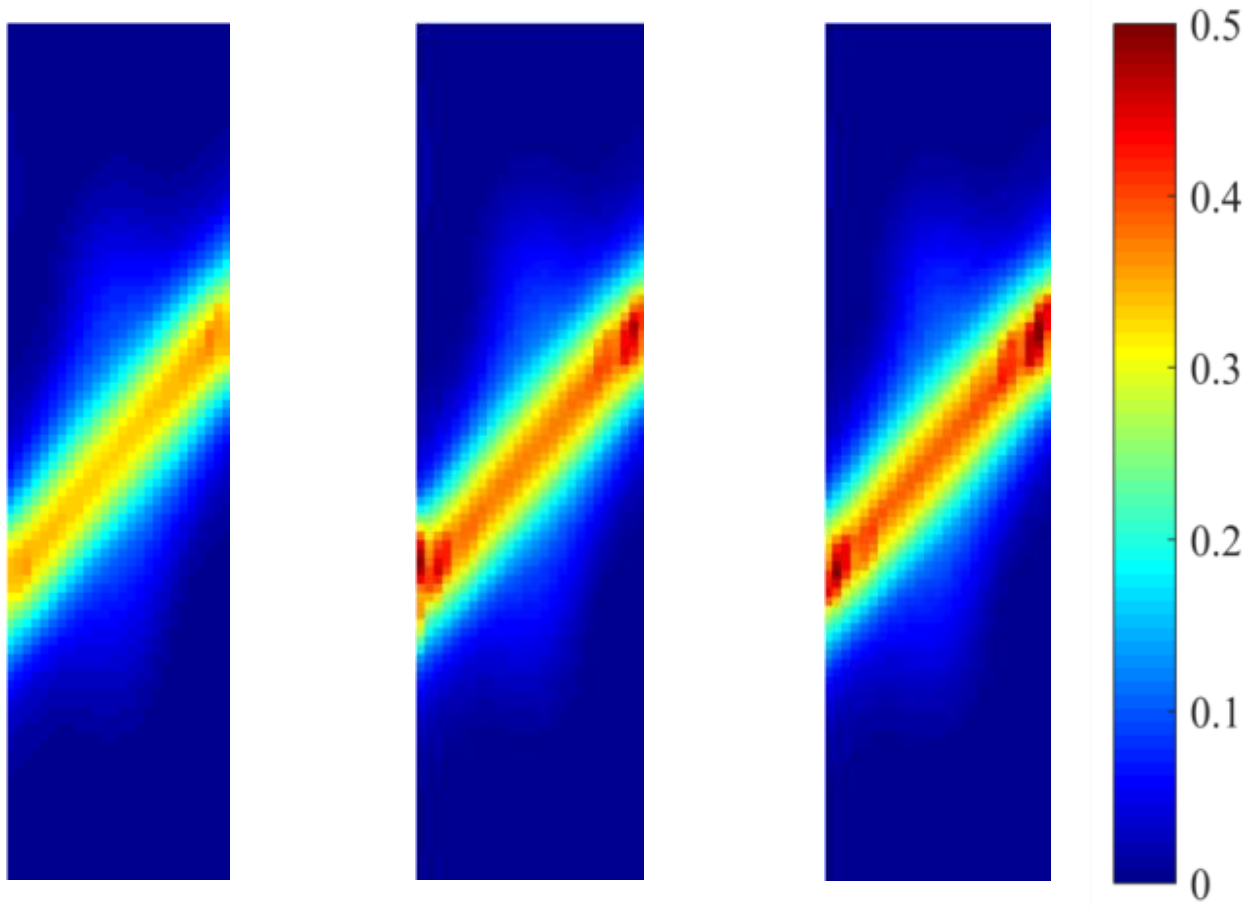

Figure 17: Contours of the equivalent plastic shear strain from the simulations with $l_1 = 1$ mm, $l_2 = 0.5$ mm, and $l_3 = 0.25$ mm.

### *4.1.3. Influence of loading rates*

In this part, we study the influence of the loading rate on the shear band formation. Figure 20 presents the loading curve for the models with the same material parameters under loading rate $\dot{u}_{y,1} = 0.175$ m/s and $\dot{u}_{y,2} = 0.12$ m/s and $\dot{u}_{y,3} = 0.0875$ m/s. As shown in Figure 20, by decreasing the loading rate, the post-peak response becomes steeper, and the plastic deformation increases. For the mentioned loading rates, the plastic shear strain, plastic volume strain, and micro rotation at the final loading stage are shown in Figures 21, 22, and 23. As demonstrated in Figure 21, the width of the shear band decreases by reducing the loading rate.

### *4.2. Example 2: Cracking in a three-point bending test*

This example simulates crack propagation and material softening in a quasi-brittle porous material. This example is inspired by the three-point bending test for cracking in a rectangle specimen of lean clay, sand, and straw fibers in [57]. Figure 24 (a) presents the experimental 3-point bending test with an initial notch in [57]. Figure 24 (b) shows the model setup for the numerical simulations in this study. The initial crack length is 44 mm. A vertical displacement load is applied at the top boundary as shown in Figure 24 (b). The vertical displacement load [30] is defined as

$$u_y = \frac{u_1}{2}\left[1 - \cos\left(\frac{\pi t}{t_1}\right)\right], \tag{136}$$

where $u_1$ = 8 mm, $t_1 = 2.28 \times 10^4\,\mu s$. The time increment is $\Delta t = 2\ \mu s$. The specimen is discretized into 4584 uniform material points with $\Delta x$ = 0.3 mm. The horizon size is chosen as $\delta = 2.05\ \Delta x$. In this example, the bilinear visco-elastic damage model is adopted for modeling the material softening. The damage parameters are assumed as $S_0$ = 0.0043 and $S_{cr}$ = 0.056. The others material parameters chosen are solid density $\rho^s$ = 2750 kg/m$^3$, initial volume fraction $\varphi_0$ = 0.89, Young's modulus $E$ = 450 MPa, Poisson's ratio $\nu$ = 0.2, Cosserat shear modulus $\mu_c = \mu/3$, relaxation time $\tau_r = 8 \times 10^3\ \mu$s, and Cosserat length scale $l$ = 2 mm.

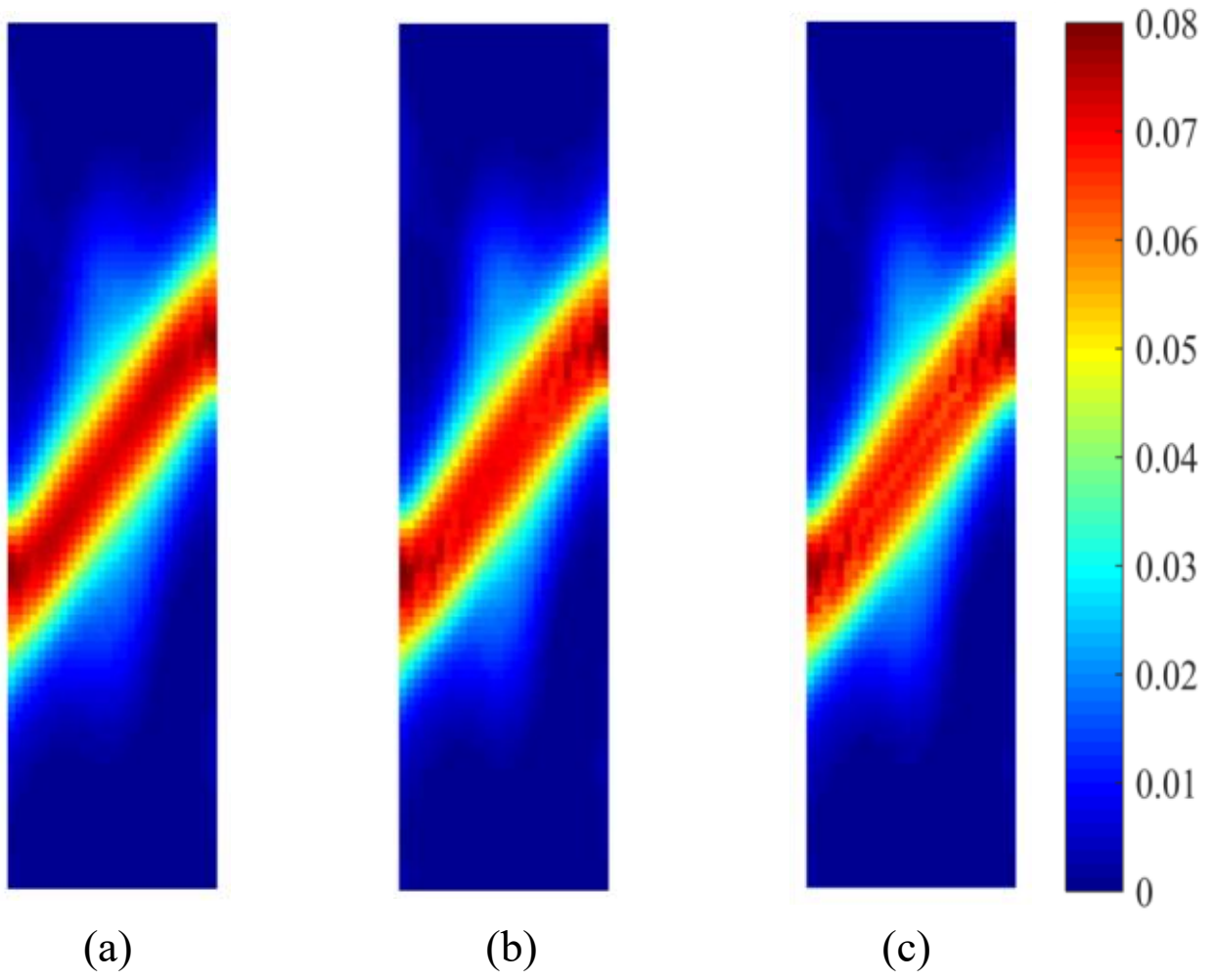

Figure 18: Contours of the plastic volume strain from the simulations with $l_1$ = 1 mm, $l_2$ = 0.5 mm, and $l_3$ = 0.25 mm.

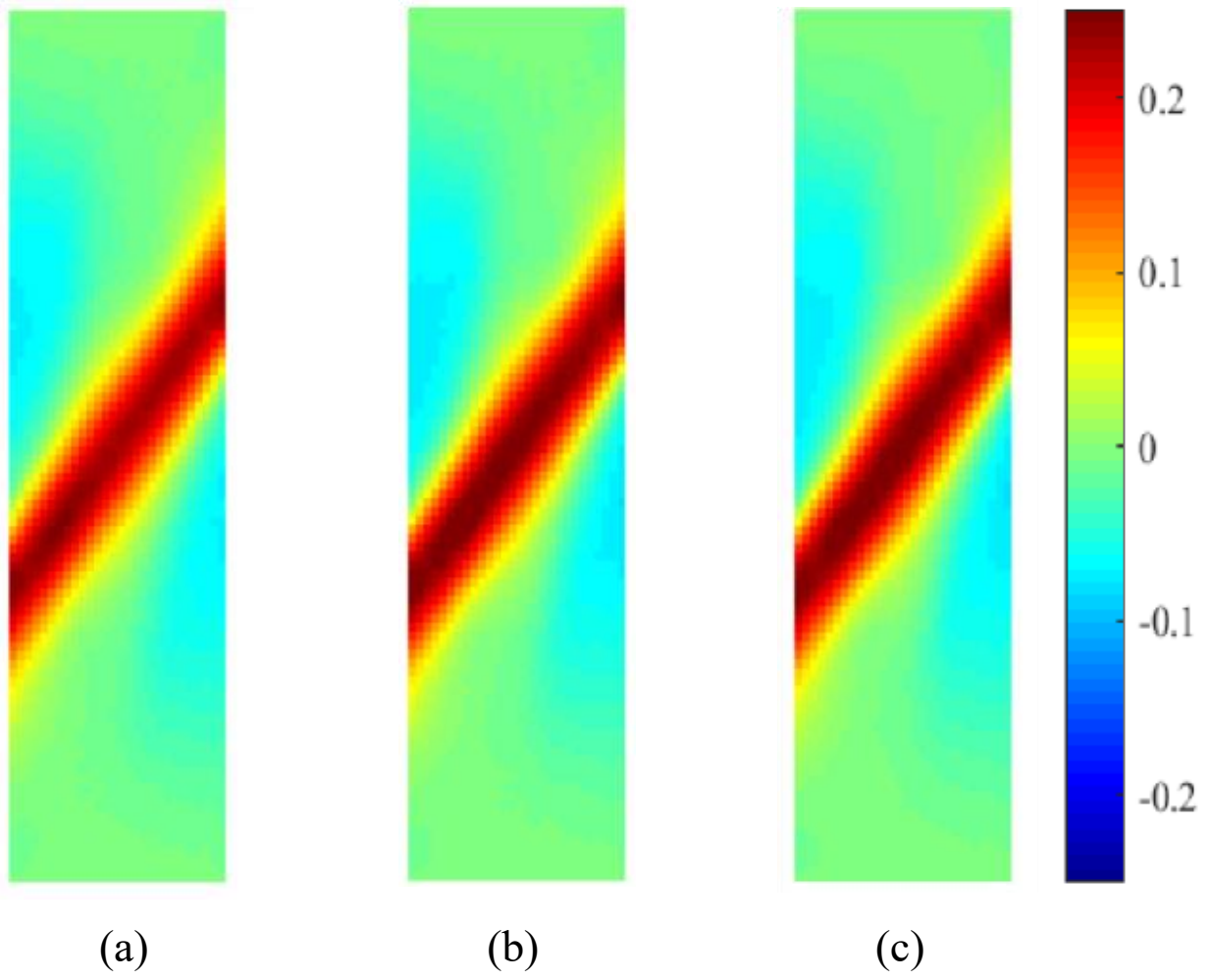

Figure 19: Contours of the micro rotation (rad) from the simulations with $l_1$ = 1 mm, $l_2$ = 0.5 mm, and $l_3$ = 0.25 mm.

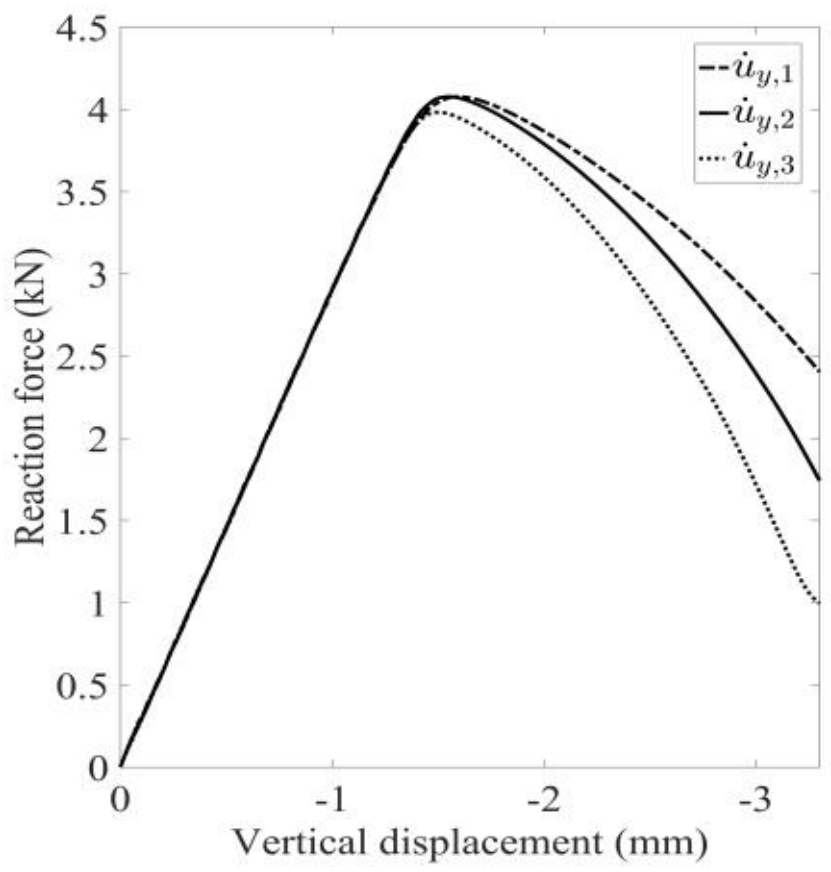

Figure 20: Comparison of the vertical loading curves from the simulations with loading rates $\dot{u}_{y,1}$ = 0.175 m/s and $\dot{u}_{y,2}$ = 0.12 m/s and $\dot{u}_{y,3}$ = 0.0875 m/s.

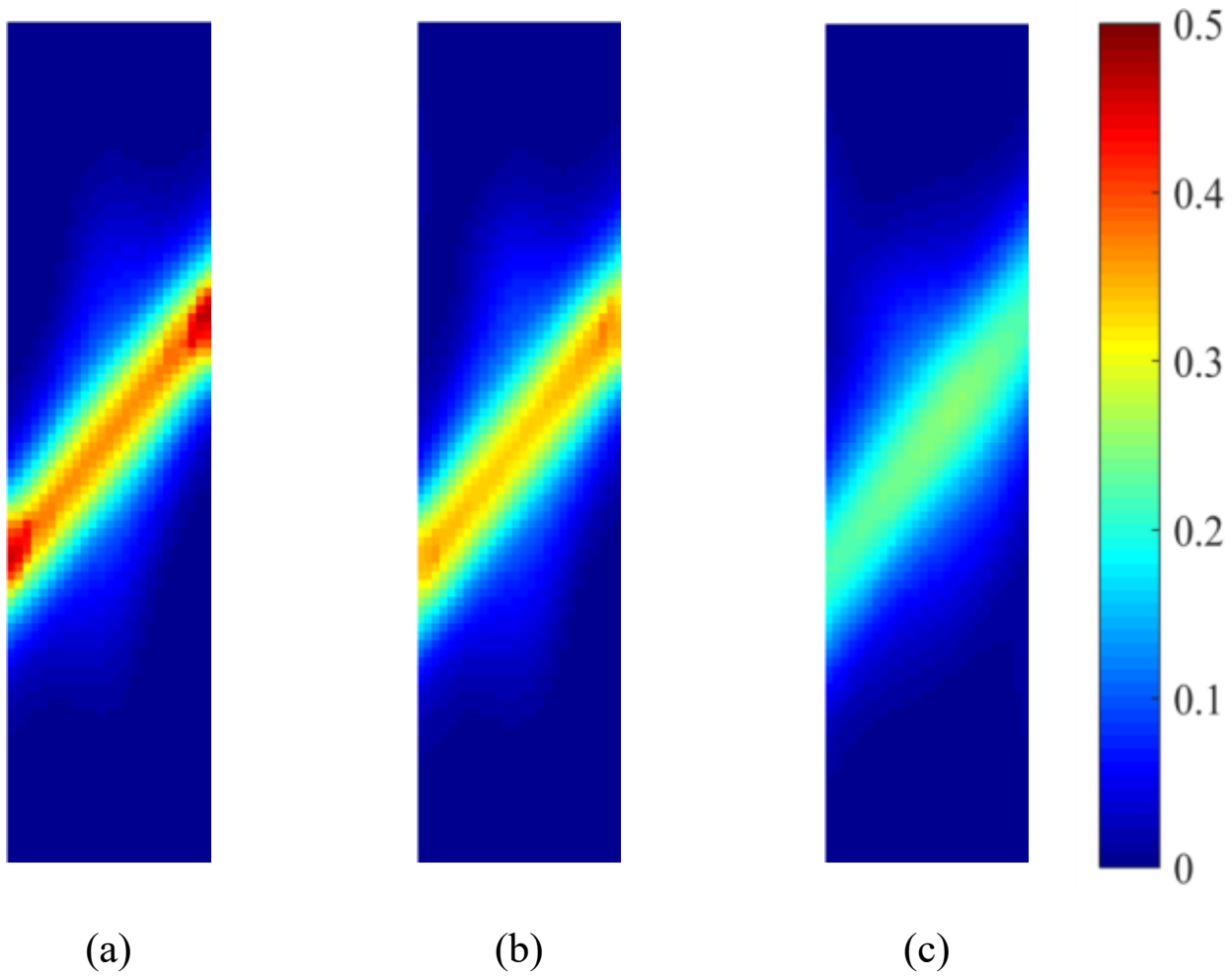

Figure 21: Contours of the equivalent plastic shear strain from the simulations with loading rates: (a) $\dot{u}_{y,1}$ = 0.0875 m/s, (b) $\dot{u}_{y,2}$ = 0.12 m/s, and (c) $\dot{u}_{y,3}$ = 0.175 m/s.

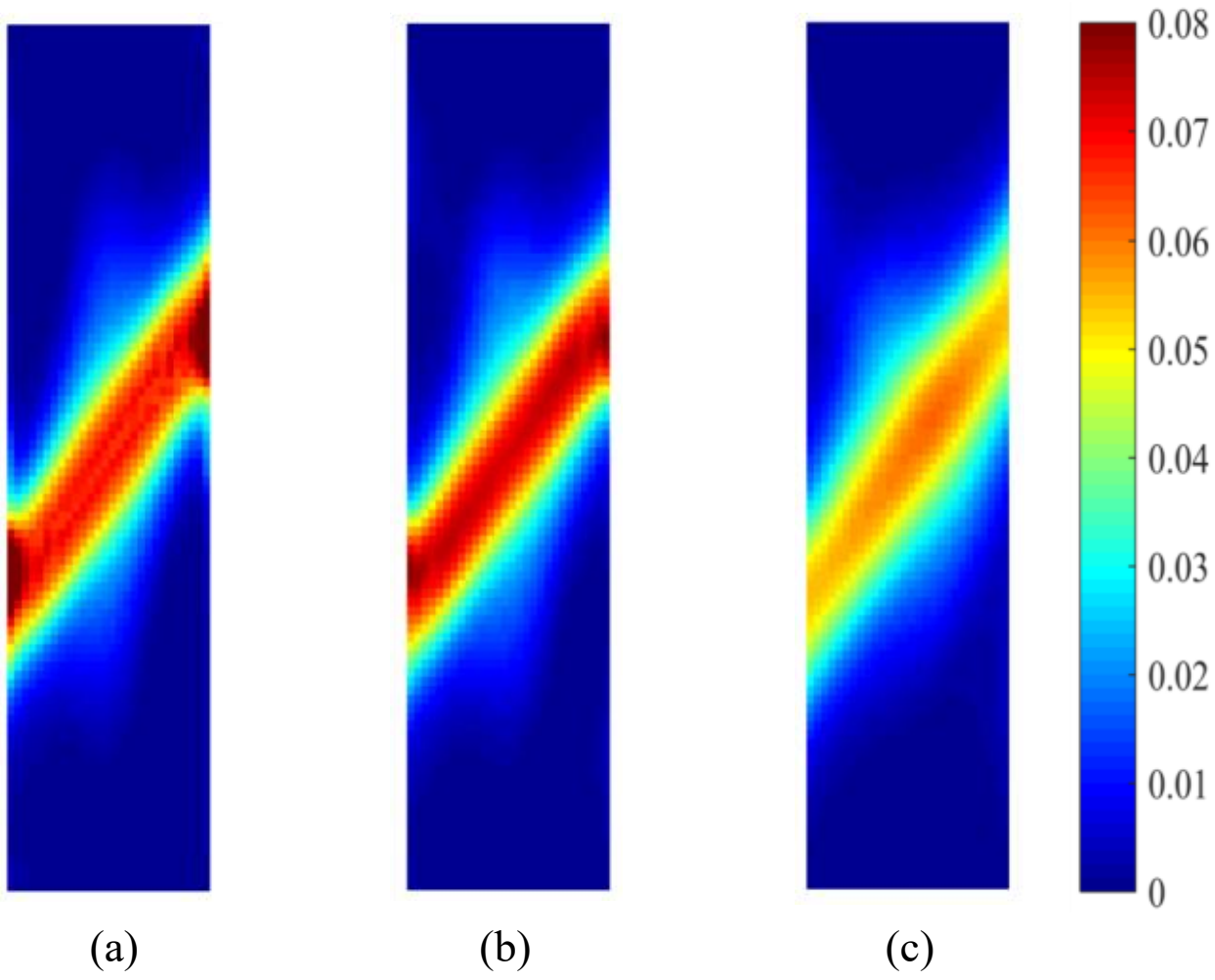

Figure 22: Contours of the plastic volume strain from the simulations with loading rates: (a) $\dot{u}_{y,1}$ = 0.0875 m/s, (b) $\dot{u}_{y,2}$ = 0.12 m/s, and (c) $\dot{u}_{y,3}$ = 0.175 m/s.

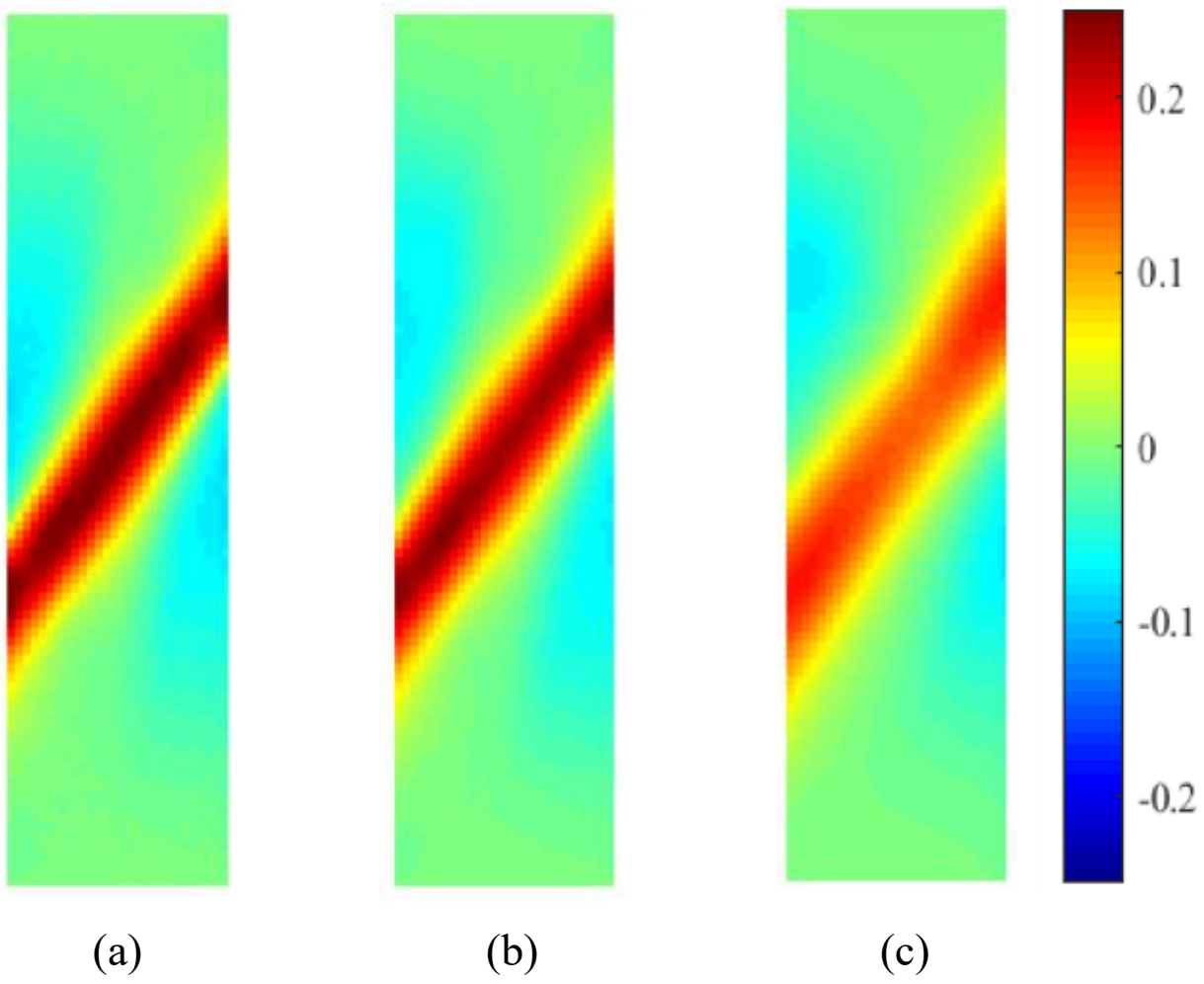

Figure 23: Contours of the micro rotation (rad) from the simulations with loading rates: (a) $\dot{u}_{y,1}$ = 0.0875 m/s, (b) $\dot{u}_{y,2}$ = 0.12 m/s, and (c) $\dot{u}_{y,3}$ = 0.175 m/s.

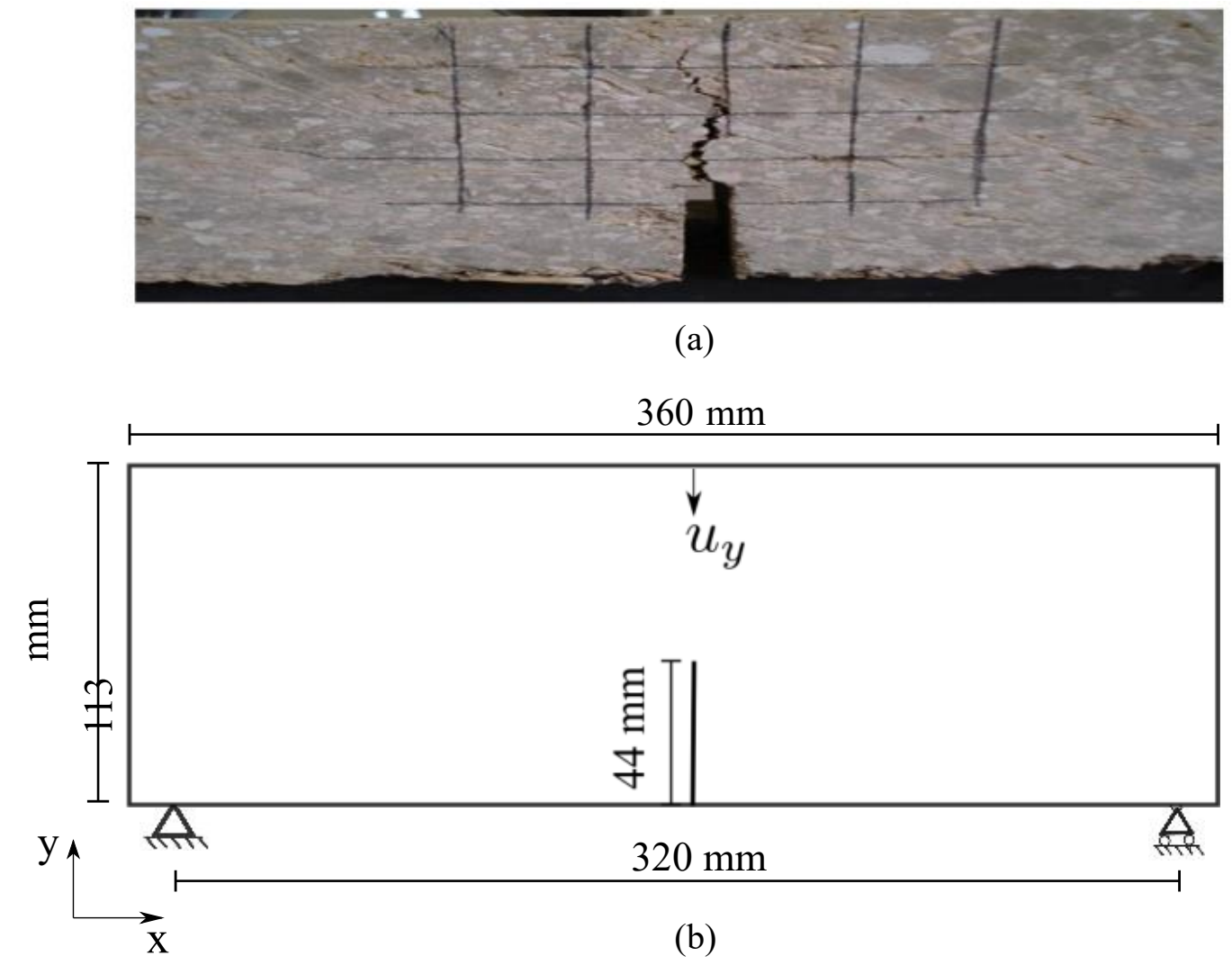

Figure 24: (a) The three-point bending test of cracking propagation of a rectangle specimen with an initial notch in [57] and (b) Model setup of the numerical simulation for the 3-point bending test.

Figure 25 plots the loading curve of the experimental testing and the numerical simulations in this study. Figure 25 shows that the numerical result is consistent with the experimental results, and the biaxial damage model can reasonably predict the softening regime of the quasibrittle material after the peak load. Nonetheless, the loading curve of the numerical simulations is not exact with the experimental data. This could be due to the inherent and induced material heterogeneity of the rectangle specimen. The four points marked on the loading curve in Figure 25 correspond to the vertical displacements $u_{y,1}$ = 0.46 mm , $u_{y,2}$ = 0.72 mm, $u_{y,3}$ = 1 mm, and $u_{y,4}$ = 1.48 mm. Figures 26 and 27 plot the contours of micro-rotation and damage variable in the deformed configurations at the four displacements, respectively. Figure 27 show that the crack propagates vertically toward the top boundary. The crack path agrees with the experimental result (see Figure 24 b). Figure 26 shows that the material points with micro-rotation are concentrated on the crack tip. The magnitude of micro-rotation increases as the crack propagates upwards.

Next, we present the results of the simulations with two spatial discretization schemes with the same horizon $\delta$ = 6 mm. For the two cases, the specimen is discretized into 4584 uniform points ($\Delta x$ = 3 mm) and 10314 uniform points ($\Delta x$ = 2 mm), respectively. Figure 28 shows the loading curves of the two simulations compared with the experimental result. The point marked in Figure 28 is at the vertical displacement $u_{y,1}$ = 1.48 mm. Figures 29 and 30 plot the contours of the micro-rotation and damage variable in the deformed configurations $u^{y,1}$ = 1.48 mm, respectively. The results in Figures 29 and 30 imply that with the same horizon the numerical results are less influenced by spatial discretization. In what follows, we investigate the factors such as loading rates and relaxation time, bilinear damage parameters $S_0$ and $S_{cr}$, and the initial volume fraction that can impact the numerical results under the same other conditions. Specifically, we present how the loading curve is sensitive to these parameters.

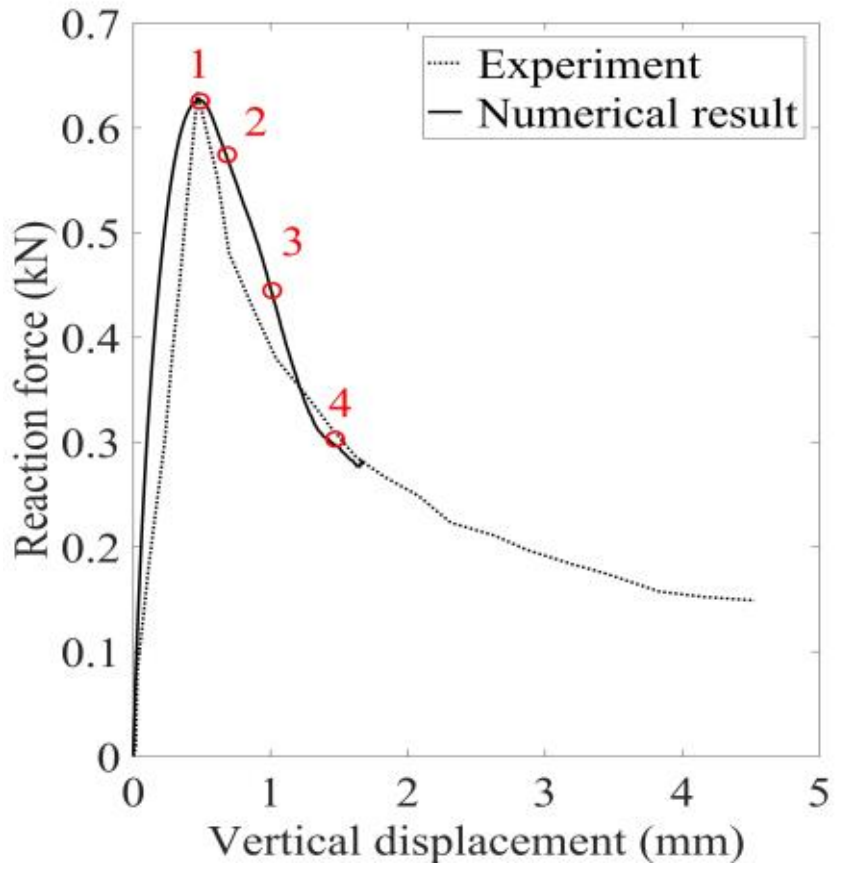

Figure 25: Comparison of loading curves of the numerical result and the experimental data.

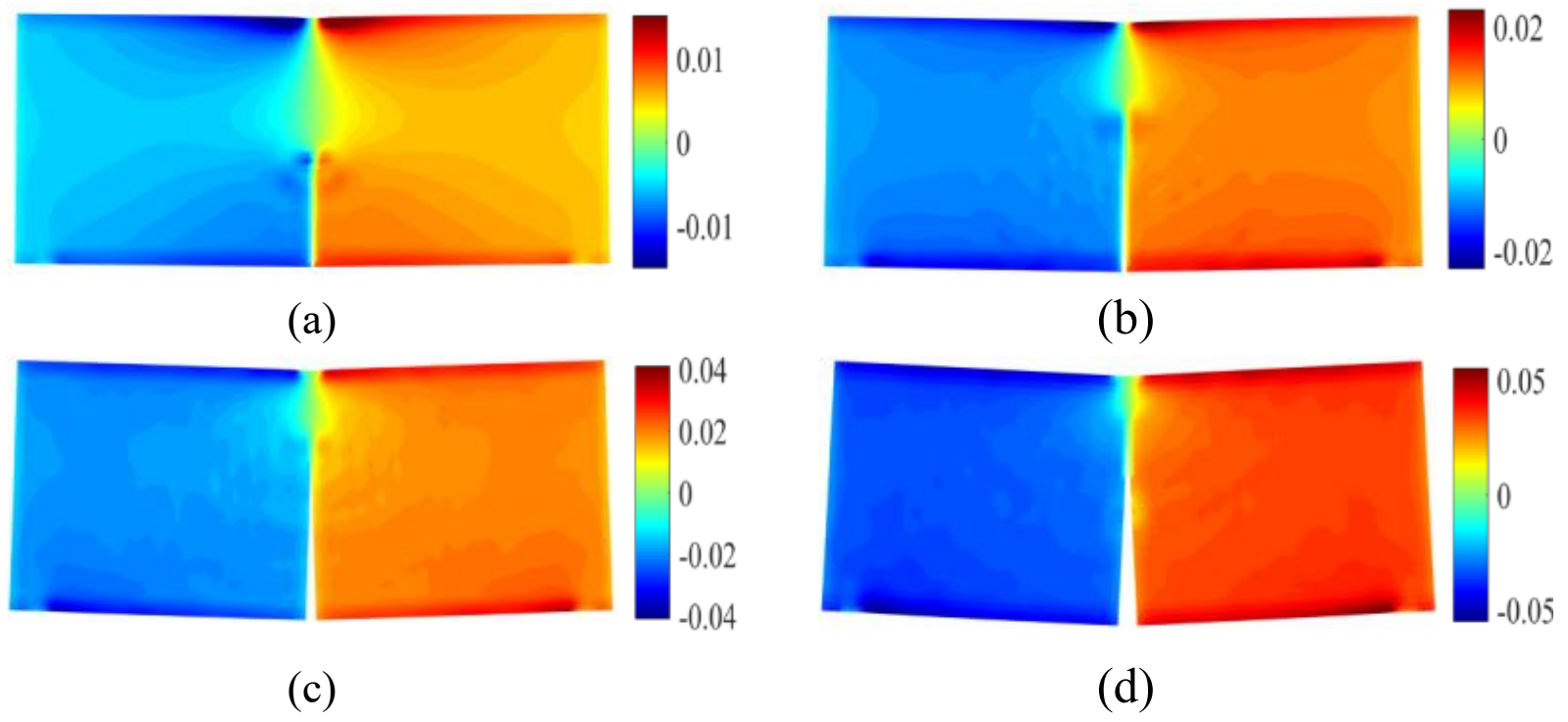

Figure 26: Contours of micro-rotation (rad) in the deformed configurations (magnification factor = 5 ) at (a) $u_{y,1}$ = 0.46 mm, (b) $u_{y,2}$ = 0.72 mm, (c) $u_{y,3}$ = 1 mm, and (d) $u_{y,4}$ = 1.48 mm.

*4.2.1. Impact of the loading rate, damage parameter, and initial volume fraction*

First, we study the influence of loading rates and relaxation time on the numerical results. The three loading rates on the top boundary are $\dot{\bar{u}}_{y,1}$ = 3.48 m/s, $\dot{\bar{u}}_{y,2}$ = 3.56 m/s and $\dot{\bar{u}}_{y,3}$ = 3.63 m/s. The other material parameters remain the same. Figure 31 compares the loading curves from the three loading rates and the experimental data. The results show that the loading rate slightly affects the peak load, e.g., a larger loading rate generates a larger peak load. The relaxation time affects the mechanical behavior of the viscoelastic material in the three-point bending test. We test the three relaxation times, $\tau_{r,1} = 8 \times 10^3\ \mu$s, $\tau_{r,2} = 6 \times 10^3\ \mu$s, and $\tau_{r,3} = 4 \times 10^3\ \mu$s. Figure 32 shows the loading curve of models with the three relaxation times and the experimental data. Figure 32 shows that the shorter relaxation time generates a smaller peak load. For the following analysis, we adopt the relaxation time $\tau_r = 8 \times 10^3\ \mu$s and the average load rate $\dot{\bar{u}}_y$ = 3.56 m/s.

Second, we study the effect of bilinear damage parameters on the numerical results. We considering different bilinear damage parameters $S_0$ and $S_{cr}$. The other material parameters are the same. The three values, $S_{0,1}$ = 0.0043, $S_{0,2}$ = 0.0086 and $S_{0,3}$ = 0.0172 are considered. Figure 33 plots the loading curves from the three $S_0$ and the experimental data. The results show the effect of $S_0$ in the softening stage. The simulation result with $S_{0,1}$ = 0.0043 is more consistent with the experimental data than the simulation results with other values of $S_0$. We also compare the results from the simulations with three values of $S_{cr}$, $S_{cr,1} = 0.056$, $S_{cr,2} = 0.112$, and $S_{cr,3} = 0.224$, as shown in Figure 34. The results show that the loading curve is less sensitive to $S_{cr}$ in the softening stage.

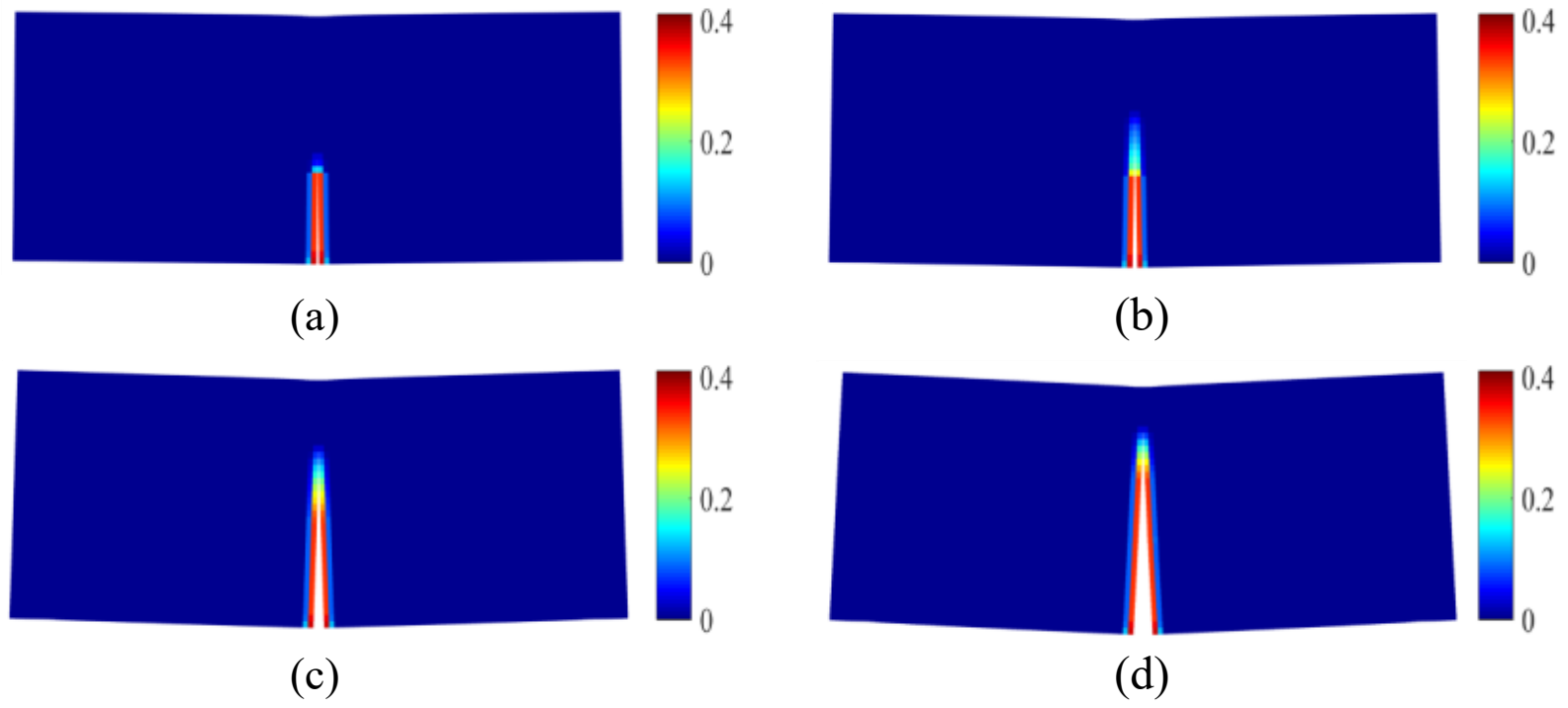

Figure 27: Contours of the damage variable in the deformed configurations (magnification factor = 5) at (a) $u_{y,1} = 0.46$ mm, (b) $u_{y,2} = 0.72$ mm, (c) $u_{y,3} = 1$ mm, and (d) $u_{y,4} = 1.48$ mm.

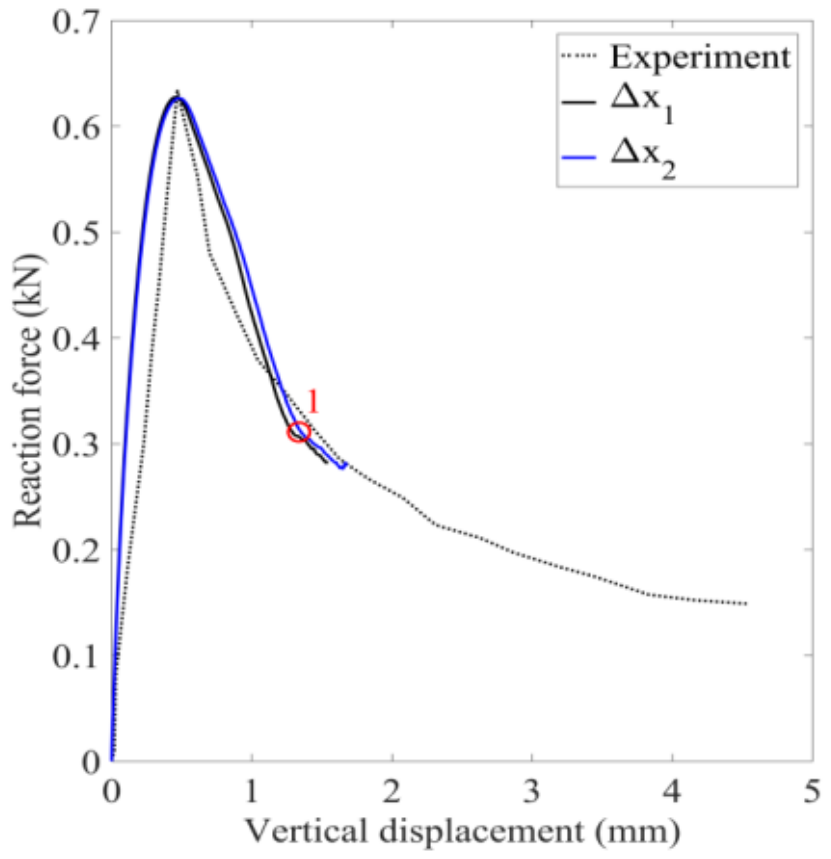

Figure 28: Comparison of the loading curves for the numerical simulations and the experimental data.

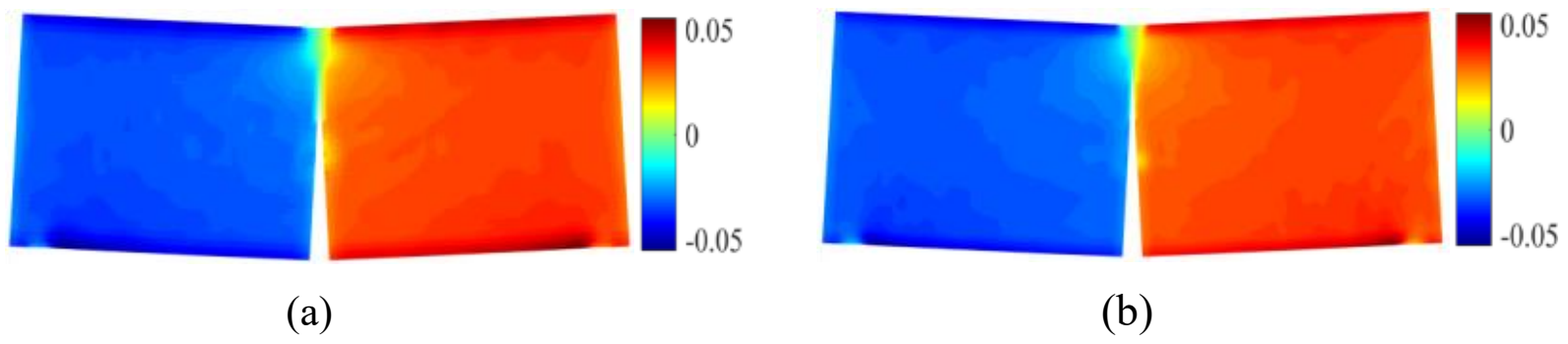

Figure 29: Contours of the micro rotation (rad) in the deformed configurations at $u_{y,1} = 1.48$ mm (magnification factor = 5) for the cases ($\delta = 6$ mm): (a) $\Delta x = 3$ mm and (b) $\Delta x = 2$ mm.

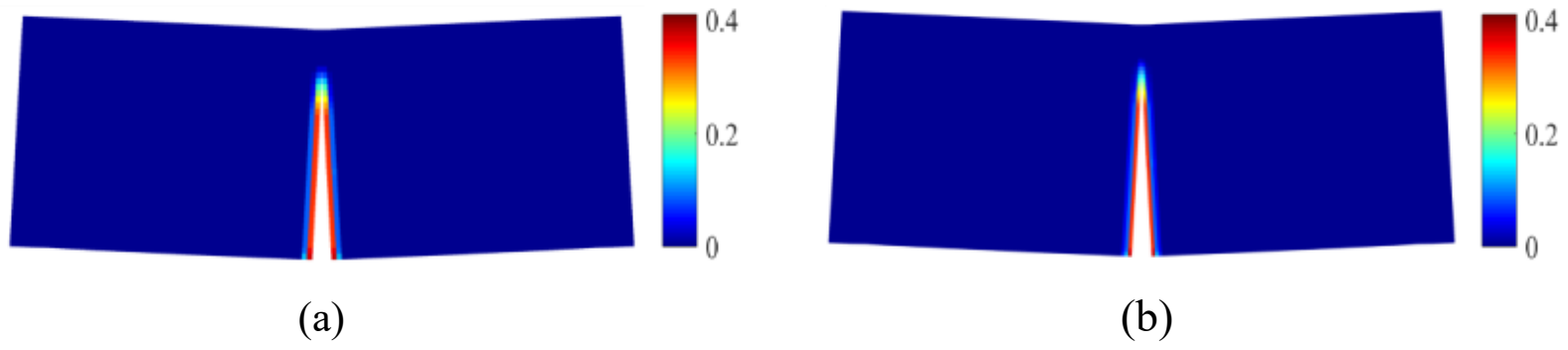

Figure 30: Contours of the damage variable in the deformed configurations at $u_{y,1} = 1.48$ mm (magnification factor = 5) for the cases ($\delta = 6$ mm): (a) $\Delta x = 3$ mm (b) $\Delta x = 2$ mm

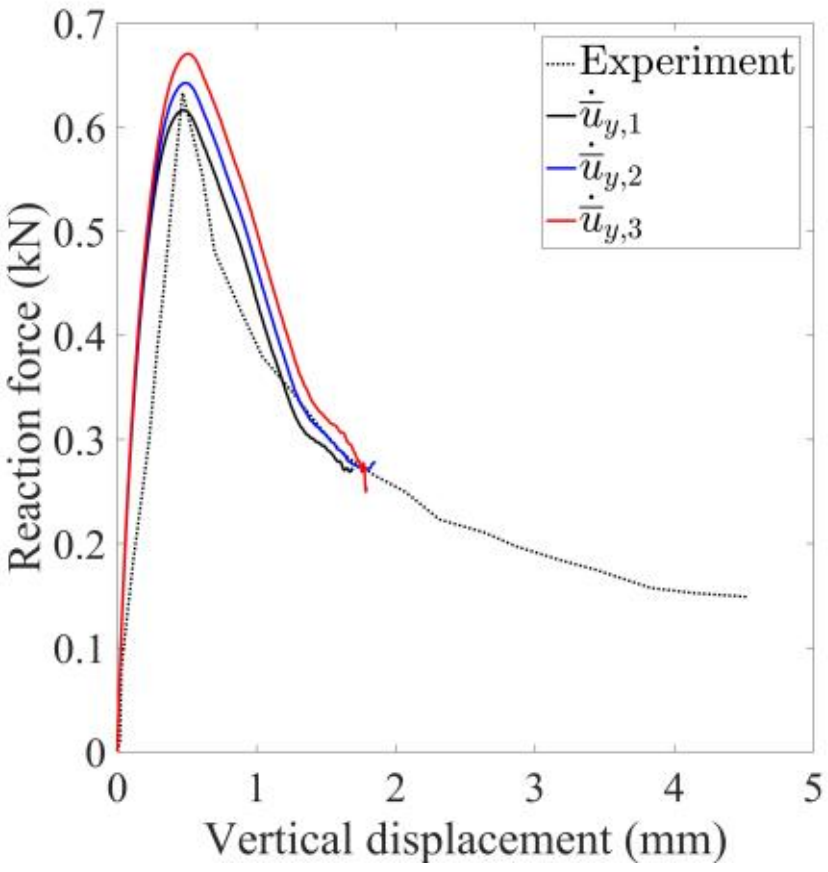

Figure 31: Comparison of the loading curves from the simulations with the average loading rates $\bar{\dot{u}}_{y,1} = 3.48$ m/s, $\bar{\dot{u}}_{y,2} = 3.56$ m/s and $\bar{\dot{u}}_{y,3} = 3.63$ m/s and the experimental data.

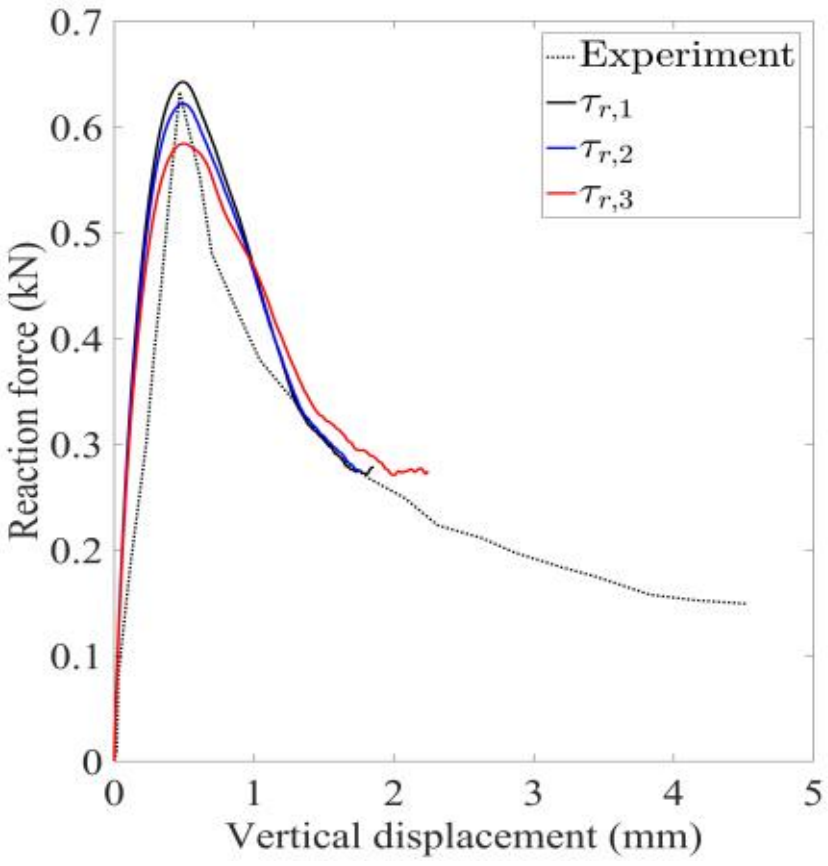

Figure 32: Comparison of the loading curves from the simulations with $\tau_{r,1} = 8 \times 10^3\,\mu$s, $\tau_{r,2} = 6 \times 10^3\,\mu$s, and $\tau_{r,3} = 4 \times 10^3\,\mu$s and the experimental data.

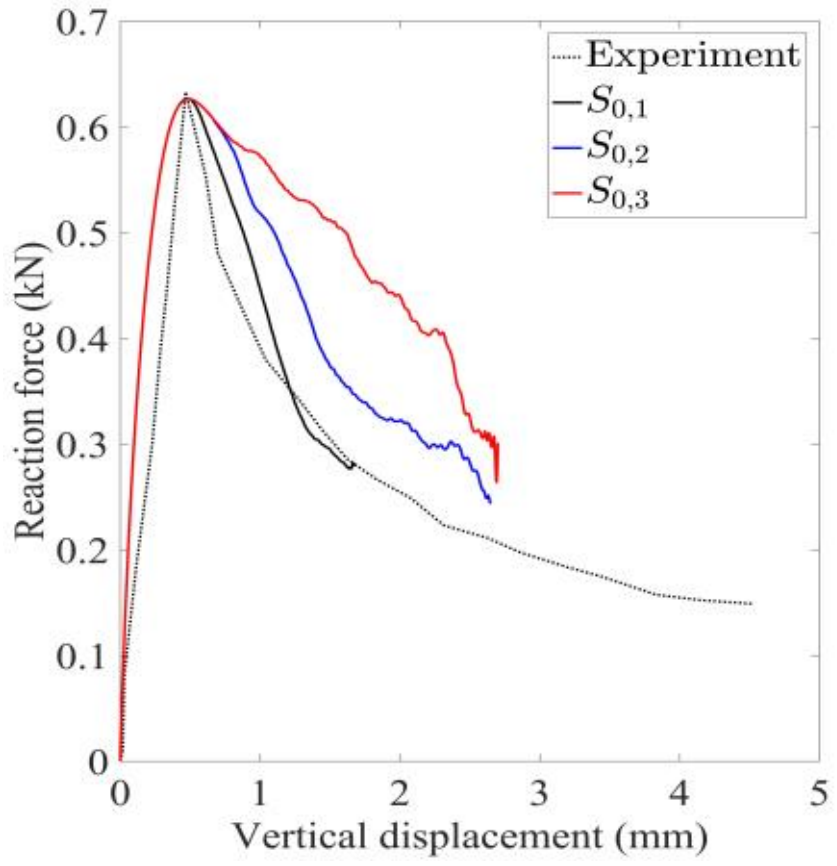

Figure 33: Comparison of the loading curves with $S_{0,1}$ = 0.0043, $S_{0,2}$ = 0.0086, and $S_{0,3}$ = 0.0172 and the experimental data.

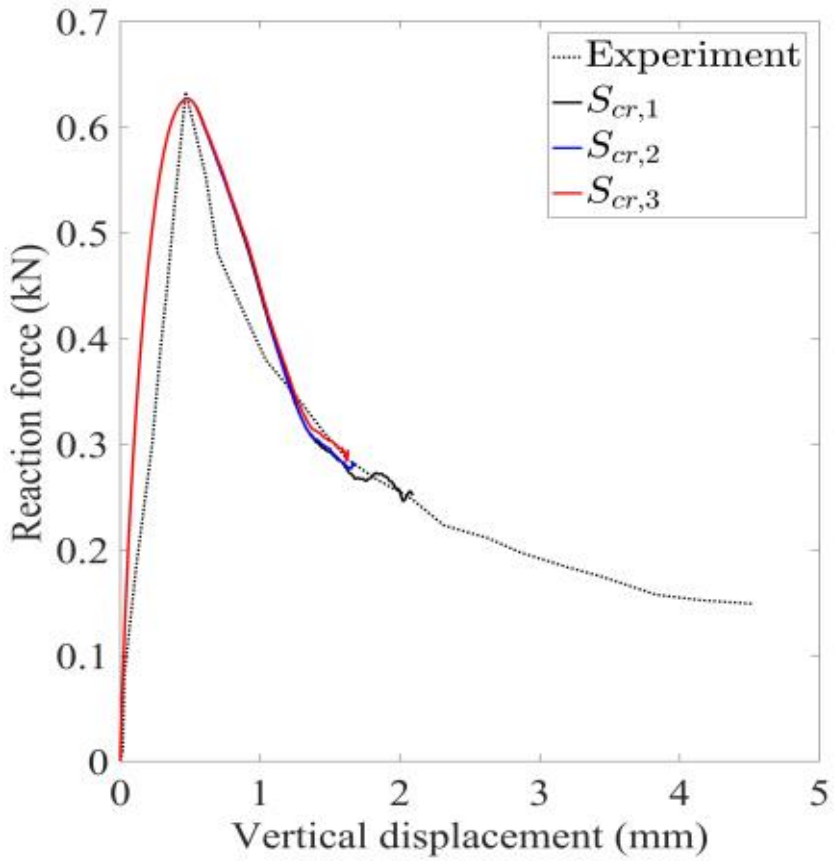

Figure 34: Comparison of the loading curves with $S_{cr,1}$ = 0.056, $S_{cr,2}$ = 0.112, and $S_{cr,3}$ = 0.224 and the experimental data.

Finally, we present the results of the loading curves assuming three initial volume fractions of solids, $\varphi_0$=0.82, 0.89, and 0.96. Figure 35 plots the loading curves from the simulations with the three initial volume fractions. The results show that the initial volume fraction does not influence the elastic regime. However, the larger volume fraction of solids generates a larger loading capacity of the specimen under other same conditions.

*4.3. Example 3: Conjugate shear banding*

This example deals with the conjugate shear banding in viscoplastic porous media under dry conditions. We investigate the impact of the dilation angle on the directions of the two conjugate shear bands in the specimen under symmetrical loading conditions. Figure 36 plot the geometry, the dimensions, the boundary conditions, and the load for this example.

The specimen is discretized into 40×80 material points with a uniform grid size $\Delta x$ = 2.5 mm and the horizon size $\delta$ = 2.05 $\Delta x$. The lateral confining pressure of 0.1 MPa is applied on the left and right boundaries. A vertical displacement is applied on the top and bottom boundaries $u_y$ = 4.5 mm with the rate $\dot{u}_y$ = 0.05 m/s. The simulation time $t = 1 \times 10^4$ $\mu$s with a stable time increment $\Delta t$ = 7 $\mu$s. The stabilization parameters $G_1$ = 0.01 and $G_2$ = 0.001 are used for stabilization in localization and post-localization, respectively [5].

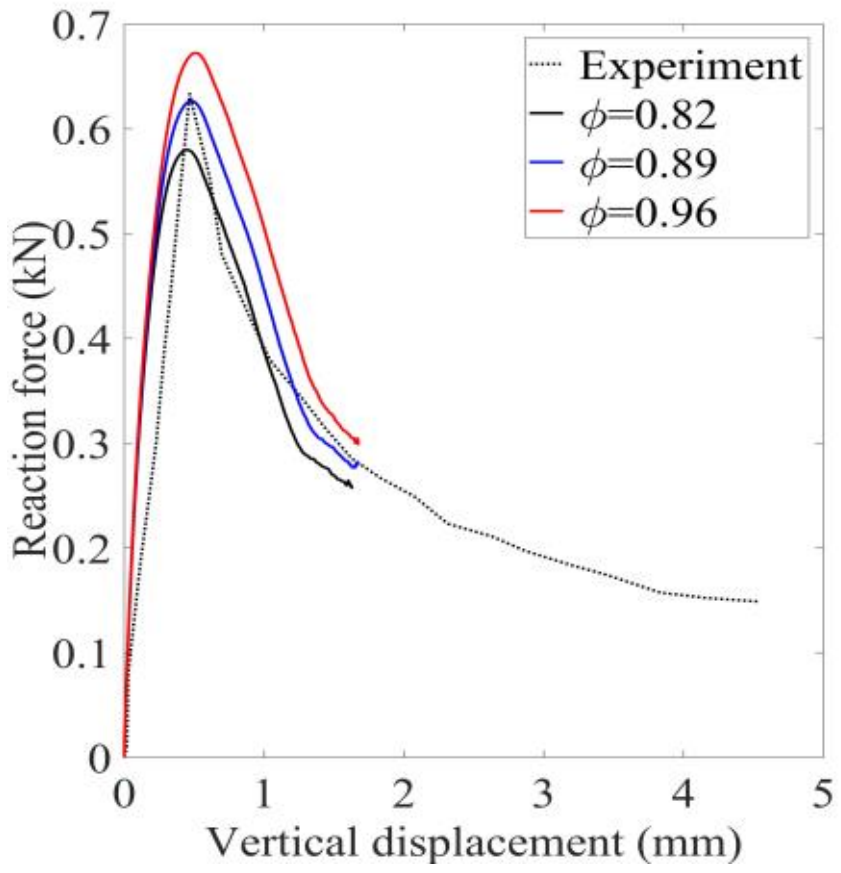

Figure 35: Comparison of the loading curves with different initial volume fractions and the experimental data.

The micro-polar visco-plastic material model is adopted for this example. The input material parameters are as follows. The density $\rho^s$ = 2000 kg/m$^3$, Young's modulus $E$ = 50 MPa, Poisson's ratio $\nu$ = 0.2, Cosserat shear modulus $\mu_c = 2\mu$, initial volume fraction $\varphi_0$ = 0.65, and Cosserat length scale $l$ = 2 mm. The viscoplastic parameters are initial cohesion $c_0$ = 0.5 MPa, softening modulus $h$ = −1 MPa, and viscosity $\eta$ = 0.01 MPa.s.

First, we investigate the influence of the dilatation angle on the inclination angle of the shear band. We consider three dilatation angles $\psi = 0^\circ, 10^\circ$, and $20^\circ$ assuming the same frictional angle $\phi = 35^\circ$. Figure 37 plots the loading curves from the simulations with the three dilation angles. Figure 37 shows that the loading curves are the same until the peak load. In the postlocalization regime, the dilatation angle has little influence on the loading curve. At the same last load step, the simulation with the null dilation angle generates the smallest reaction force. Figures 38 and 39 plot the contours of equivalent plastic shear strain and plastic volumetric strain from the simulations with three dilation angles at the same last load step, respectively. Figure 40 plots the micro rotation contours on the deformed configurations from the three simulations at the same end load step. The results in Figures38, 39, and 40 demonstrate that the dilation angle affects the inclination of the two conjugate shear bands. Figure 40 shows that the micro rotation of material points is localized within the shear band.

Table 1 compares the inclination angle of the shear band in this example with the classical Roscoe solution. Our numerical solution is consistent with the Roscoe solution [18]. The secondorder work is useful to detect shear bands in porous media [58–60]. Therefore, the second-order work criterion is used to validate our numerical results. The second-order work $dW$ for a micropolar PPM material can be written as

$$dW = \overline{d\sigma} : d\epsilon + dm : d\kappa. \tag{137}$$

Figure 41 plots the contours of the second-order work from the three dilation angles. The results show that the second-order work within the shear band is negative for all three cases.

Table 1: Comparison of the inclination angles ($^\circ$) from the close-form solution and the numerical results.

| $\phi(^\circ)$ | $\Psi(^\circ)$ | Roscoe solution | Numerical Solution |
|---|---|---|---|
| 35 | 0 | 45 | 41.3 |
| 35 | 10 | 50 | 46.6 |
| 35 | 20 | 55 | 50.9 |

Second, we study the influence of loading rates on the formation of shear bands with three loading rates, $\dot{u}_{y,1} = 0.083$ m/s, $\dot{u}_{y,2} = 0.067$ m/s and $\dot{u}_{y,3} = 0.05$ m/s. The frictional angle $\phi = 35^\circ$ and the dilatation angle $\psi = 15^\circ$. The other material parameters remain the same. Figure 42 presents the loading curves from the three loading rates. It shows that the loading rate mainly affects the peak and post-localization regimes. Figures43, 44, and 45 present the

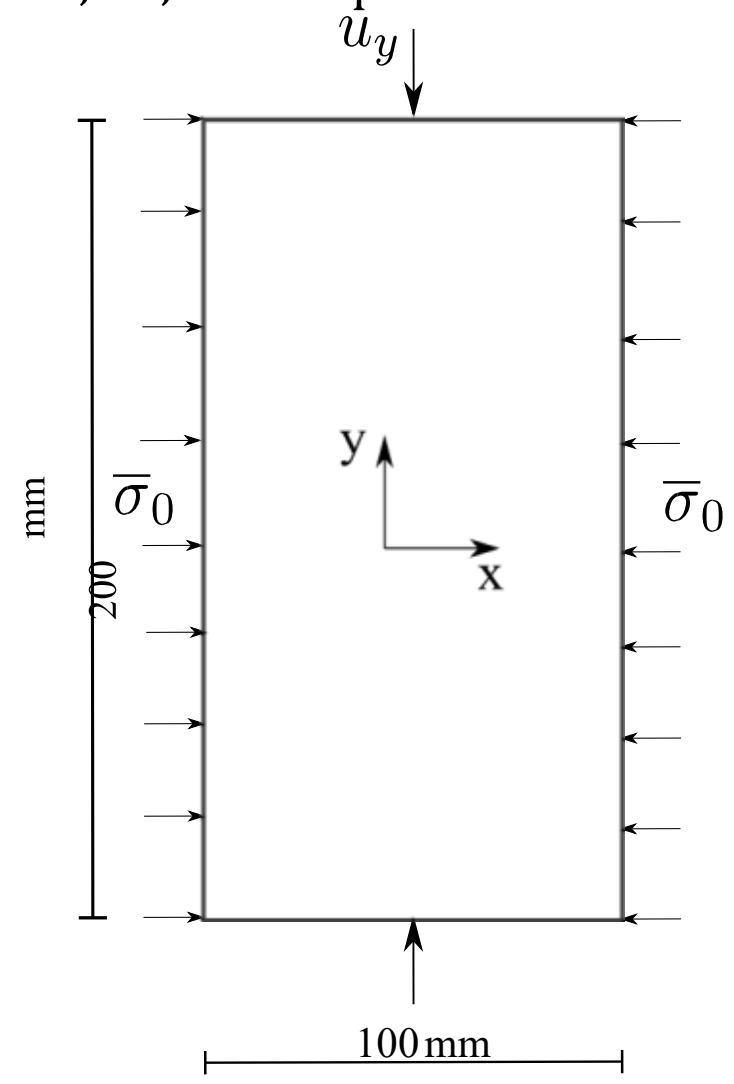

Figure 36: Geometry, boundary conditions, and the loading of the conjugate shear banding.

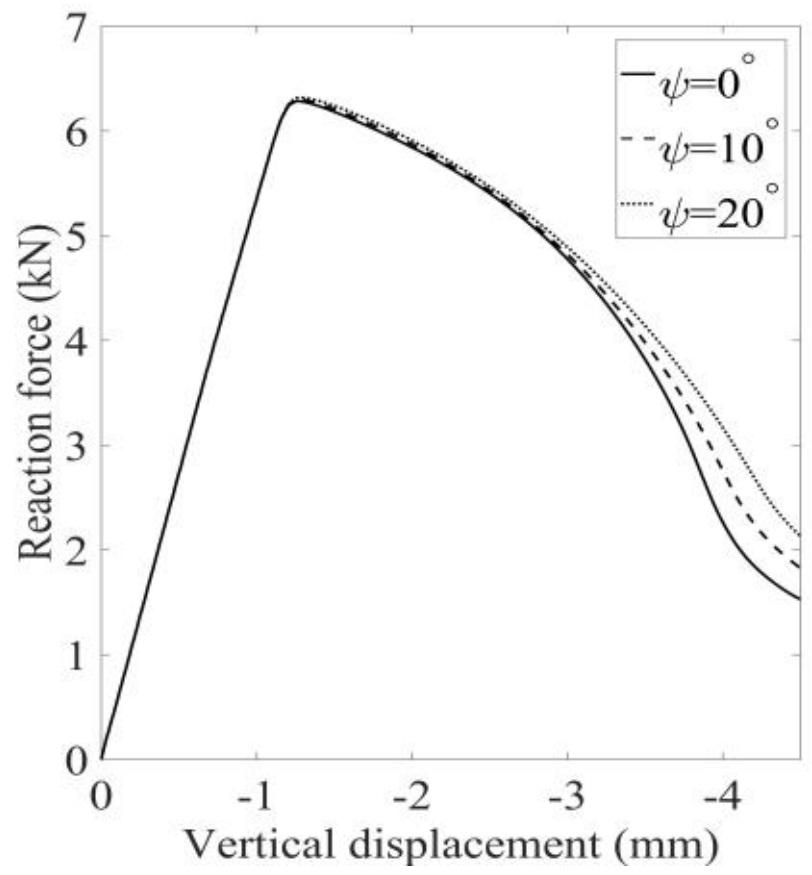

Figure 37: Loading curves assuming three dilatation angles $\psi = 0^\circ$, $\psi = 10^\circ$, and $\psi = 20^\circ$ (the same frictional angle $\phi = 35^\circ$).

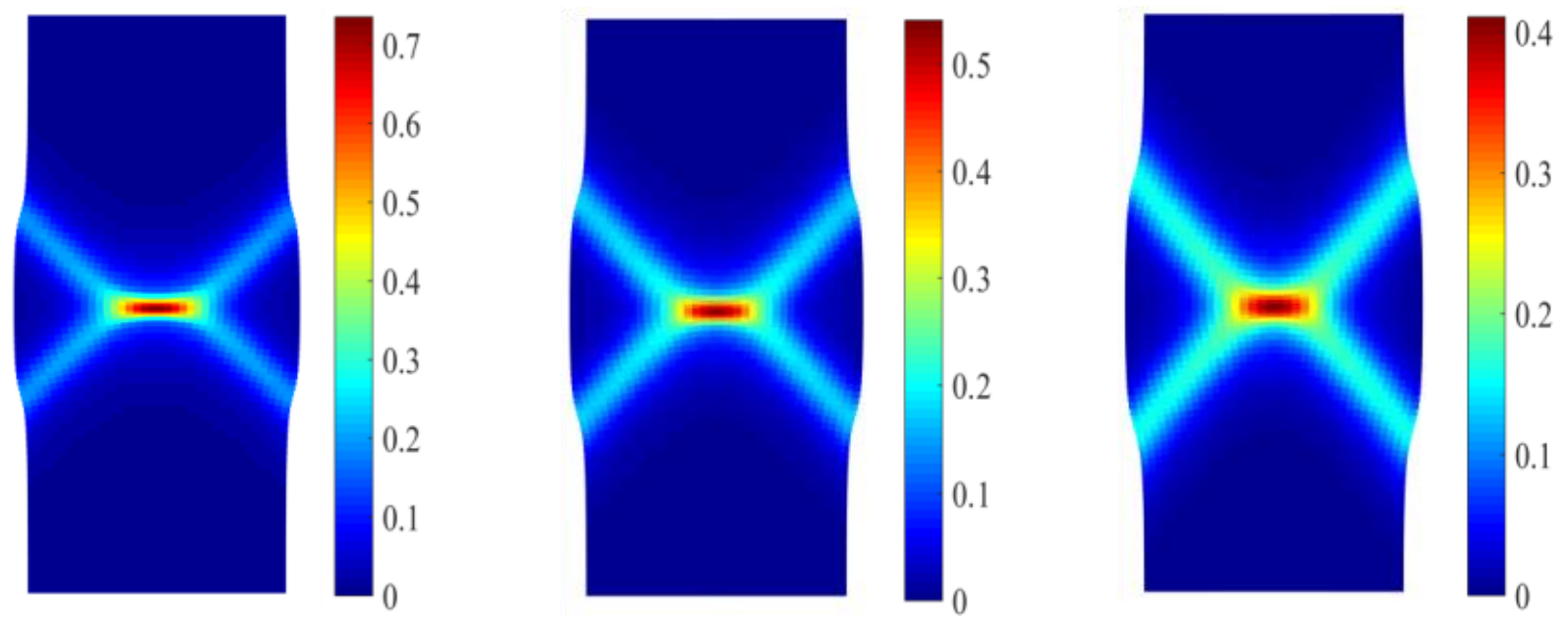

(a) (b) (c)

Figure 38: Contours of the plastic shear strain assuming three dilatation angles: (a) $\psi = 0^\circ$ (inclination angle $\theta = 41.3^\circ$), (b) $\psi = 10^\circ$ (inclination angle $\theta = 46.3^\circ$), and (c) $\psi = 20^\circ$ (inclination angle $\theta = 50.9^\circ$) (the same frictional angle $\phi = 35^\circ$).

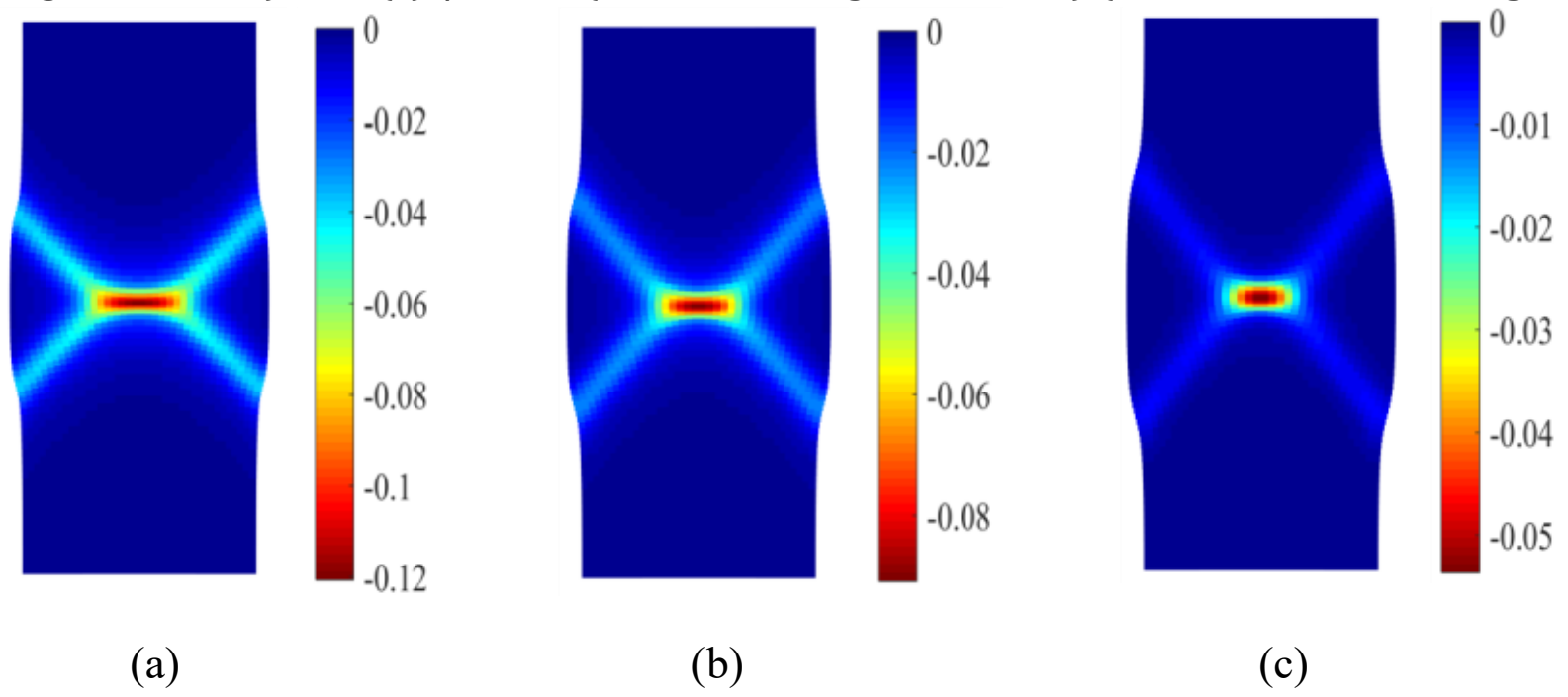

Figure 39: Contours of the plastic volume strain assuming three dilatation angles: (a) $\psi = 0^\circ$, (b) $\psi = 10^\circ$, and (c) $\psi = 20^\circ$ (the same frictional angle $\phi = 35^\circ$).

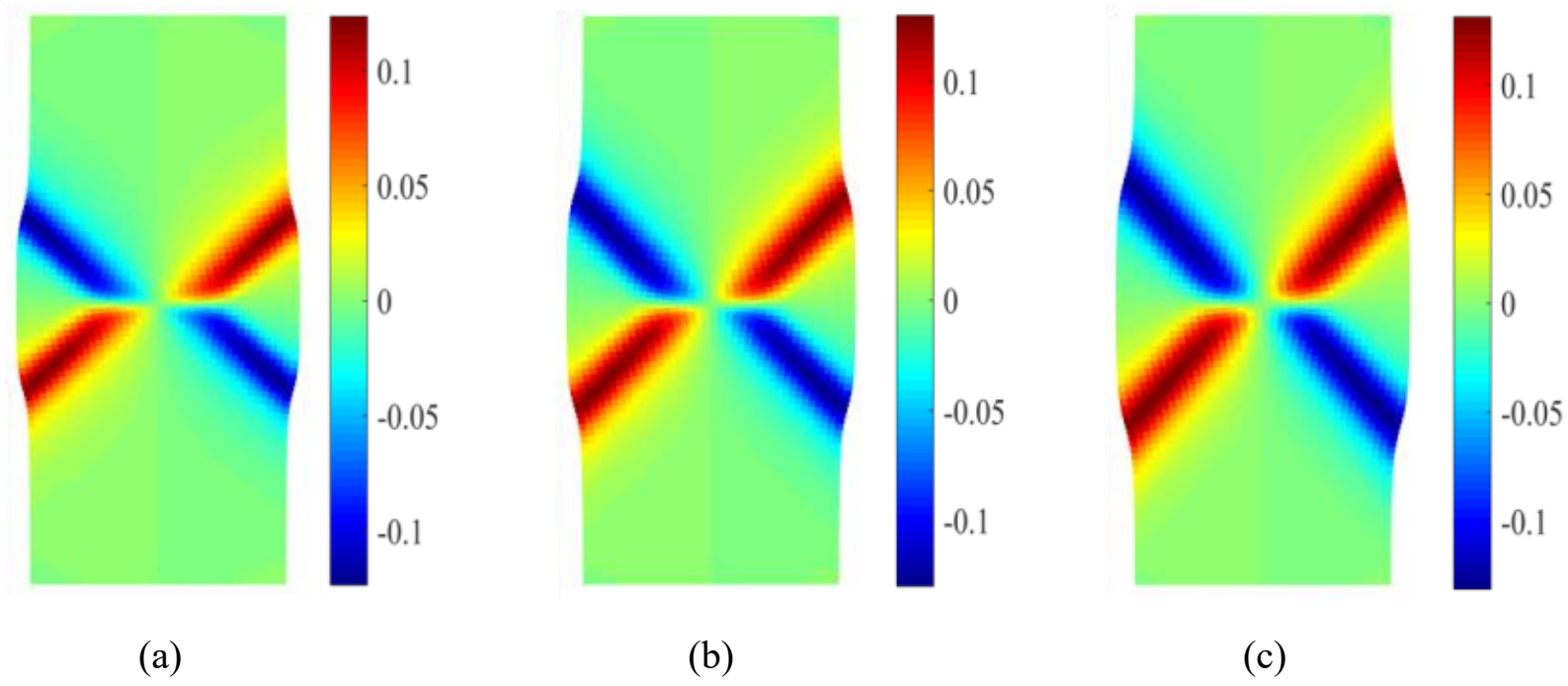

Figure 40: Contours of the micro rotation (rad) assuming three dilatation angles: (a) $\psi = 0^\circ$, (b) $\psi = 10^\circ$, and (c) $\psi = 20^\circ$ (the same frictional angle $\phi = 35^\circ$).

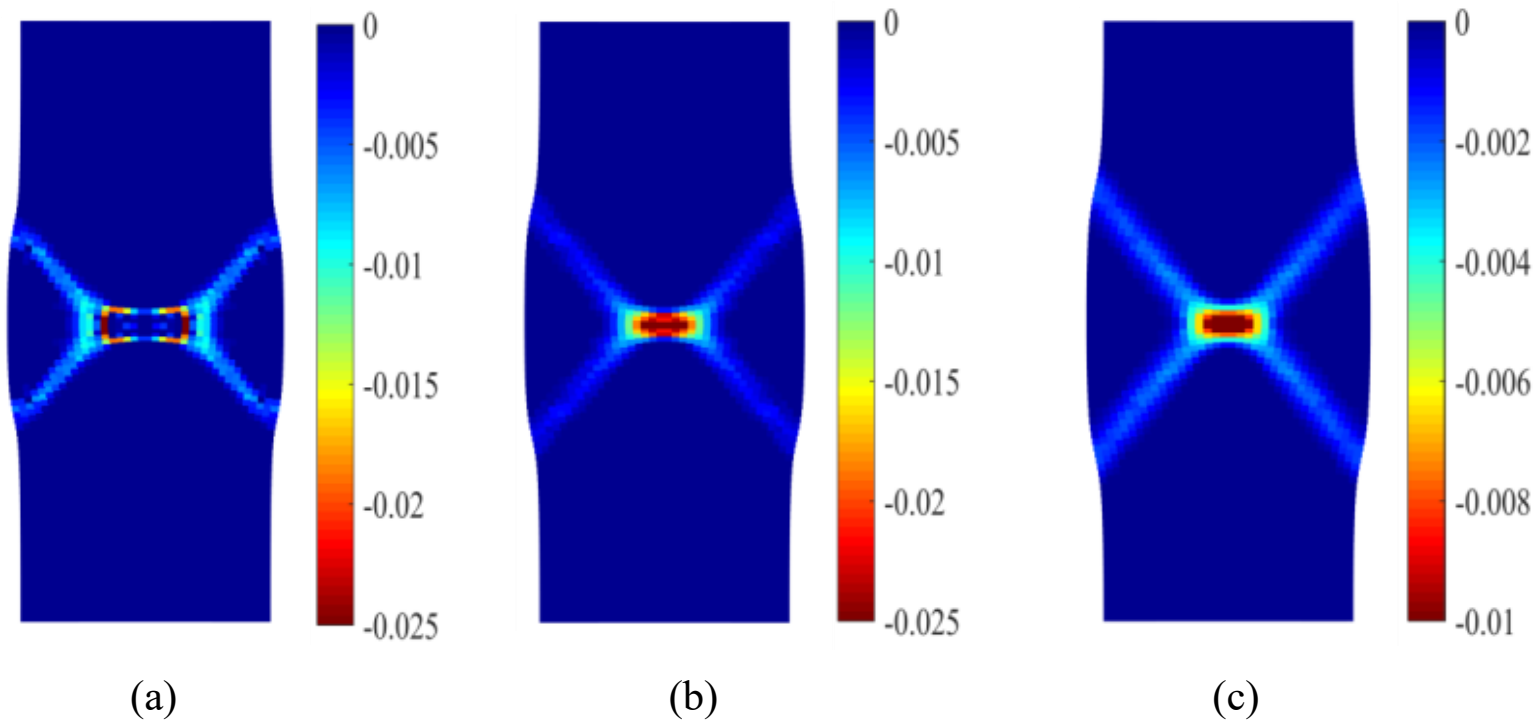

Figure 41: Contours of the second order work assuming three dilatation angle: (a) $\psi = 0^\circ$, (b) $\psi = 10^\circ$, and (c) $\psi = 20^\circ$ (the same frictional angle $\phi = 35^\circ$).

contours of equivalent plastic shear strain, micro rotation, and plastic volume strain at the same final loading stage, respectively. The results show that the loading rate affects the width of shear bands, e.g., decreasing loading rates decreases the width of shear bands. Furthermore, Figure 46 demonstrates that the zone of negative second-order work is consistent with the location of shear bands. It is implied that the augmented second-order work can be used to detect the formation of shear bands through the proposed micro-polar periporomechanics framework.

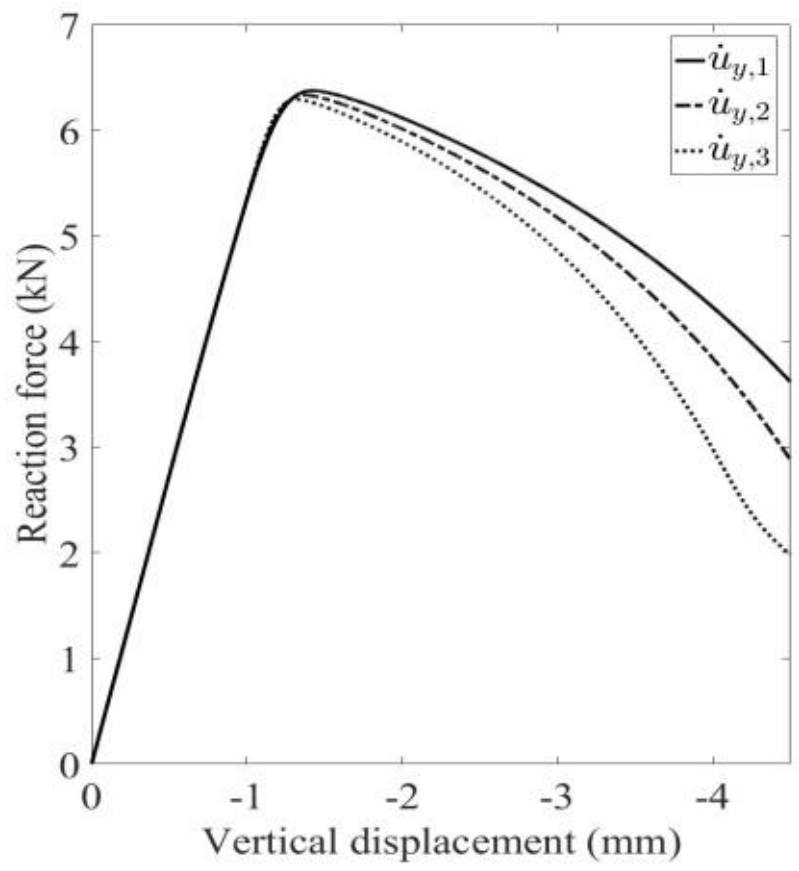

Figure 42: Loading curves from three loading rates, $\dot{u}_{y,1}$ = 0.083 m/s, $\dot{u}_{y,2}$ = 0.067 m/s, and $\dot{u}_{y,3}$ = 0.05 m/s.

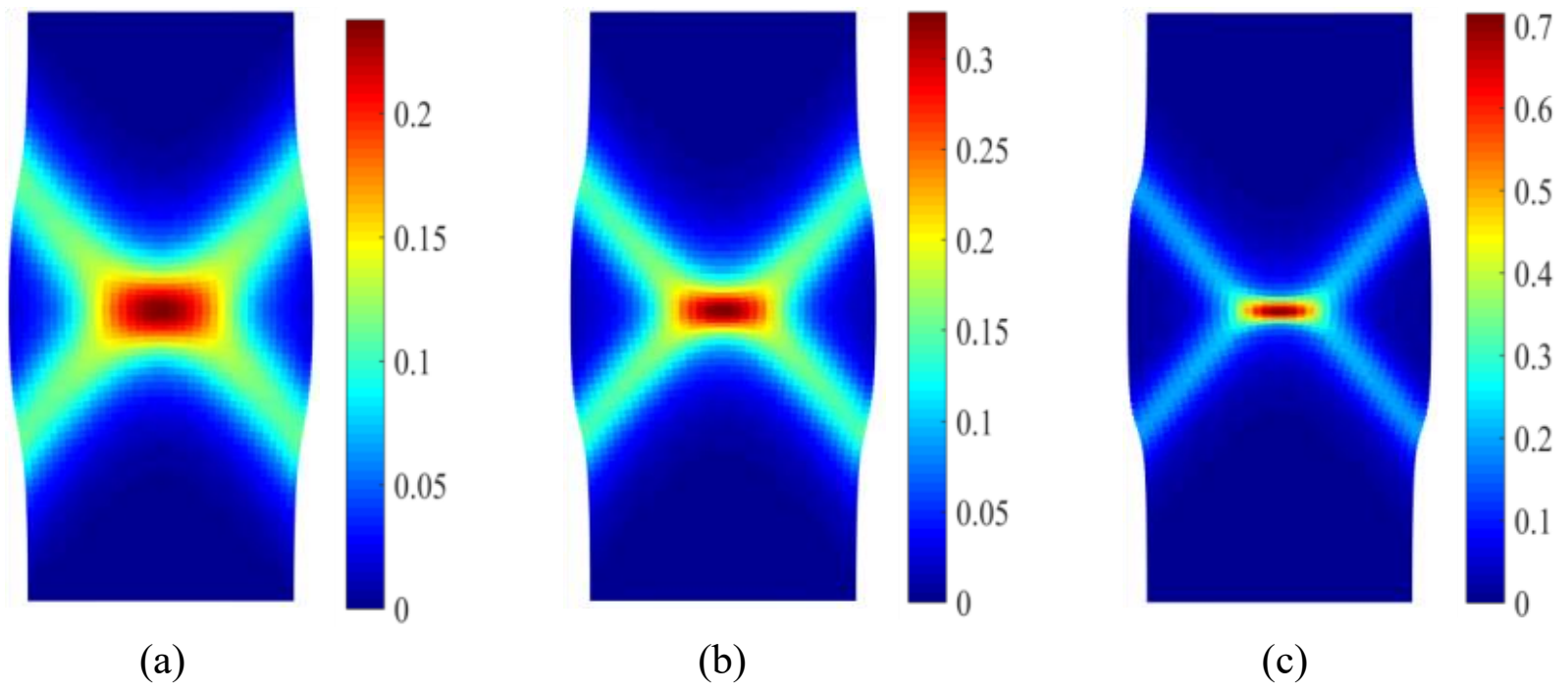

Figure 43: Contours of the plastic shear strain from three loading rates: (a) $\dot{u}_{y,1}$ = 0.083 m/s, (b) $\dot{u}_{y,2}$ = 0.067 m/s, and (c) $\dot{u}_{y,3}$ = 0.05 m/s.

### *4.4. Example 4: Crack branching under high loading rates*

This example deal with the crack branching in a dry visco-elastic porous material under high loading rates through the proposed energy-based cracking criterion considering the micro-rotation of material points at the crack tip. Crack branching in porous media such as clay can be related to the mechanical properties of the clay layers [61]. When the stress in a fracture zone is high in porous materials, the material cannot dissipate the energy, and the crack starts to branch due to a small critical energy release rate $G_c$. The crack branching is observed when the crack reaches the critical speed of propagation (e.g., [62]). The failure mode could change from mode I to mixed modes by increasing the loading rate. The inertia forces at the crack tip prevent crack propagation when the crack propagates fast, resulting in branching (e.g., [62]). The experimental studies have demonstrated that crack growth velocity in porous materials is lower than the Rayleigh wave velocity (e.g., [63, 64]). We note that the heterogeneity of materials can affect the crack branching in porous media, which is beyond the scope of this article [65]. In this example, we study the impact of loading rates, Cosserat length scales, and initial volume fractions on crack branching. We also investigate the micro rotation of material points along the crack path in crack branching.

Figure 47 shows the geometry, boundary conditions, and loading in this example. The specimen's dimensions are 100 mm × 40 mm × 10 mm. The initial crack is 50 mm long, as shown in

$$\sigma_y = \begin{cases} \dfrac{\sigma_1 t}{t_0} & \text{if } t < t_0, \\ \sigma_1 & \text{if } t \geq t_0, \end{cases} \tag{138}$$

where $t_0 = 6.25$ $\mu$s, and $\sigma_1 = 8$ MPa. A stable time step $\Delta t = 0.025$ $\mu$s. The specimen is discretized into 200 × 80 uniform material points with $\Delta x = 0.5$ mm. The horizon size is $\delta = 4.05\Delta x$ [66].

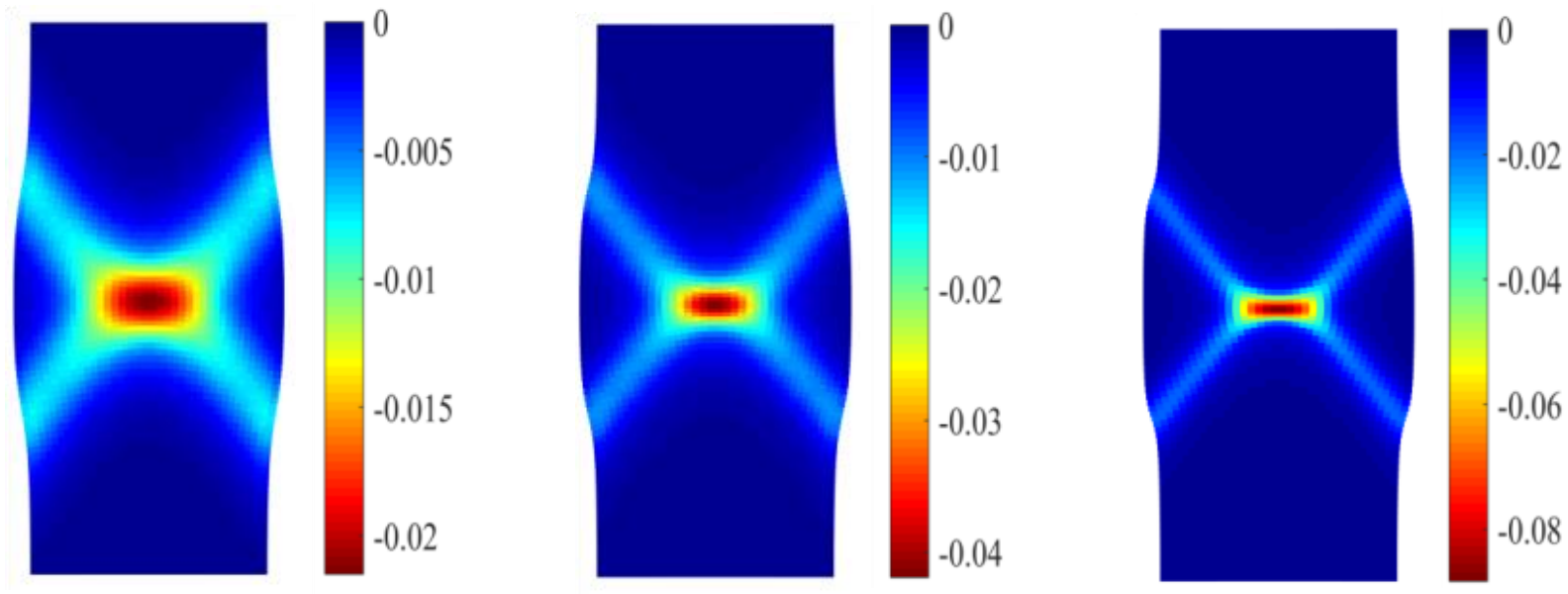

(a) (b) (c)

Figure 44: Contours of the plastic volume strain from three loading rates: (a) $\dot{u}_{y,1} = 0.083$ m/s, (b)
$\dot{u}_{y,2} = 0.067$
m/s, and (c) $\dot{u}_{y,3} = 0.05$ m/s.

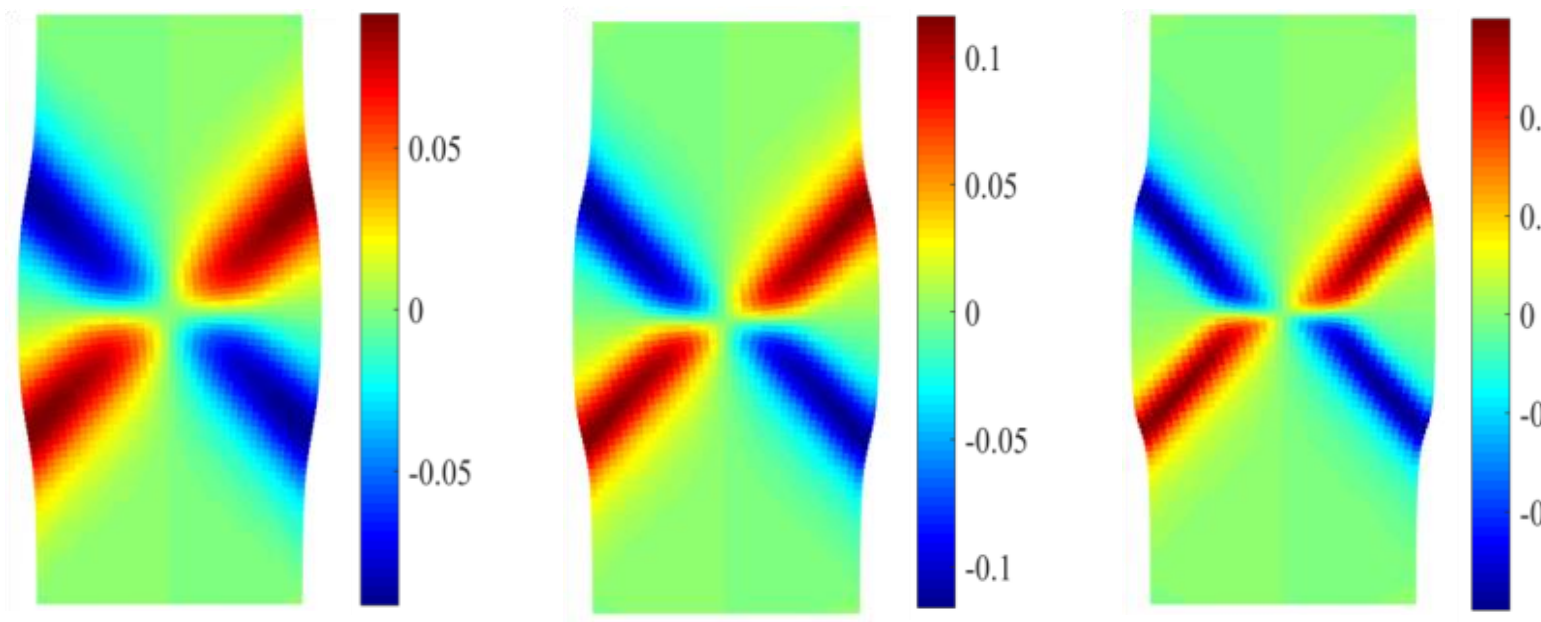

(a) (b) (c)

Figure 45: Contours of the micro rotation (rad) from three loading rates: (a) $\dot{u}_{y,1} = 0.083$ m/s, (b)
$\dot{u}_{y,2} = 0.067$
m/s, and (c) $\dot{u}_{y,3} = 0.05$ m/s.

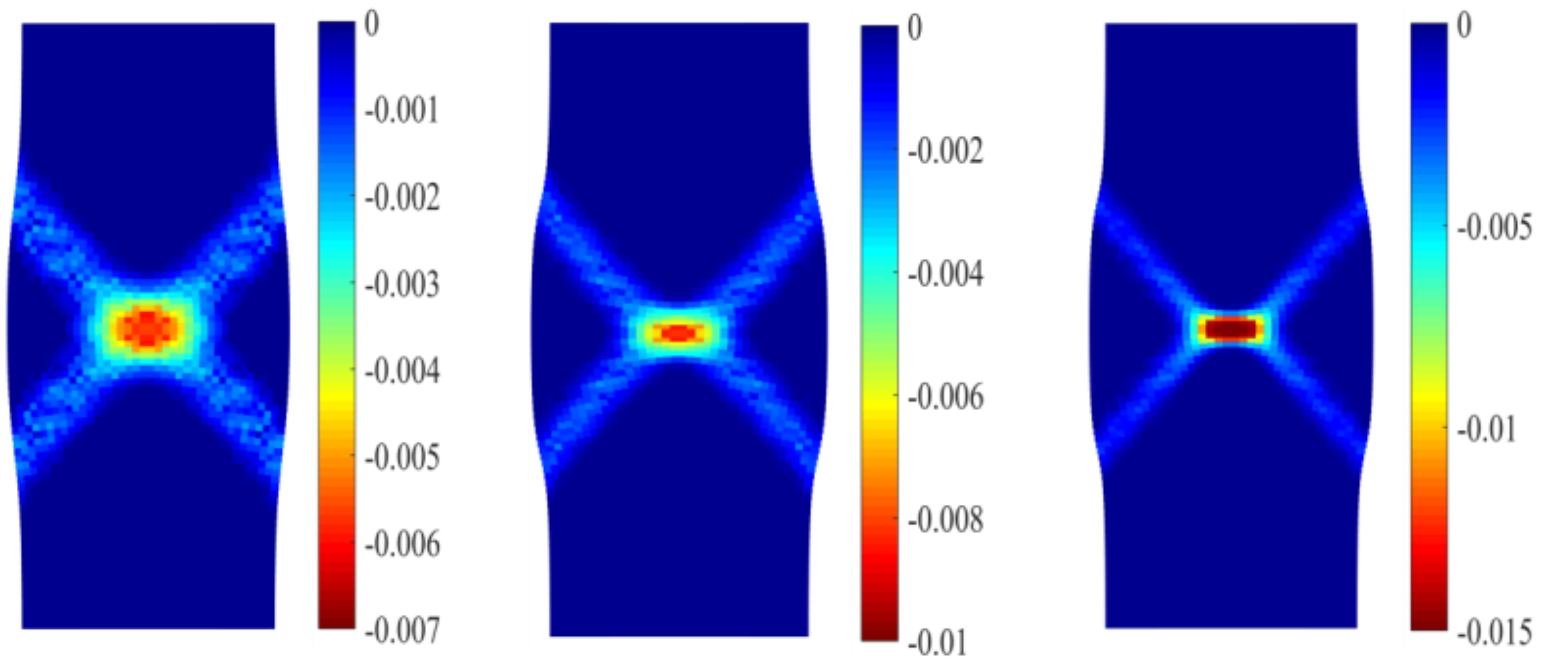

(a) (b) (c)

Figure 46: Contours of the second-order work from three loading rates: (a) $\dot{u}_{y,1} = 0.083$ m/s, (b) m/s, and (c) $\dot{u}_{y,3} = 0.05$ m/s.
$\dot{u}_{y,2} = 0.067$

Figure 47. The tensile stress in the vertical direction is applied on the top and bottom boundaries by the following equations.

The micro-polar visco-elastic material model in Section 2 is adopted for this example. The input material parameters are summarized as follows. The solid density $\rho^s = 2650$ Kg/m$^3$, initial volume fraction $\varphi_0 = 0.95$, Young's modulus $E = 35$ MPa, Poisson's ratio $\nu = 0.25$, Cosserat shear modulus $\mu_c = \mu/3$, Cosserat length scale $l = 2$ mm, relaxation time $\tau_r = 100$ $\mu$s.

For this example, $G_{cr} = 160$ N/m is assumed for the energy-based bond criteria. The stabilization parameters $G_1 = 0.1$ and $G_2 = 0.01$ are used.

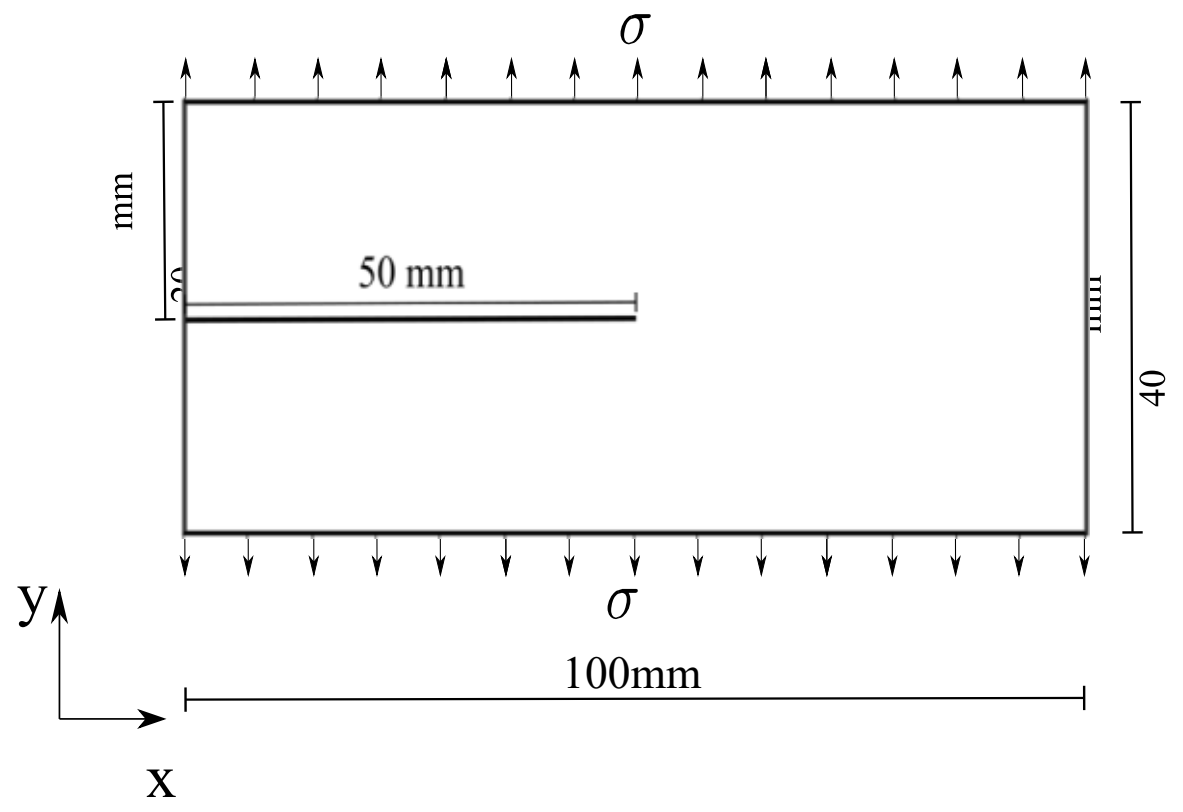

Figure 47: Model setup for example 4.

Figure 48 plots the applied loading-time curve. The four points shown in Figure 48 are at times $t_1 = 12.75$ $\mu$s, $t_2 = 29.75$ $\mu$s, $t_3 = 35$ $\mu$s, and $t_4 = 42.5$ $\mu$s. The results of the base simulations are presented in Figures 49 and 50. Figure 49 plots the snapshot of the crack propagation and branching on the deformed configuration at the four loading stages. The damage variable greater than 0.35 is in red [67]. The results show that at time $t_1 = 12.75$ $\mu$s the crack starts to grow and at time $t_2 = 29.75$ $\mu$s the crack start branching. Figure 50 plots the snapshots of the contours of micro-rotation of material points in the specimen. The results in Figure 50 show that the microrotation of material points is concentrated on the crack tip and crack propagation and branching paths, and its value increases with the crack growth. In what follows, we study the influence of the loading rate, the Cosserat internal length scale, and the initial volume fraction on the crack branching in this example.

*4.4.1. Influence of loading rates*

In this part, we analyze the results of crack branching under three loading rates assuming the same conditions. For the three loading rates, $t_{0,1} = 2.5$ $\mu$s, $t_{0,2} = 6.25$ $\mu$s, and $t_{0,3} = 12.5$ $\mu$s. For all three simulations, it is assumed that $\varphi_0 = 0.85$ and $l = 2$ mm. The other material parameters are the same as the base simulation.

Table 2 summarizes the timing of the crack propagation, the beginning of crack branching, and the end of crack branching. As the loading rate increases, the times for crack growth and the start of branching decrease. Figures 51 and 52 compare the contour of the damage variable at $t = 42.5$ $\mu$s and the contour of micro-rotation of material points at $t = 35$ $\mu$s for the three loading rates, respectively.

*4.4.2. Influence of Cosserat length scale*

We study the impact of the Cosserat length scale on the cracking branching. The three Cosserat length scales adopted are $l = 1$ mm, 2 mm and 3 mm. It is assumed that $\varphi_0 = 0.85$ and $t_0 = 6.25$ $\mu$s while the other parameters are the same as the base simulation. Table 3 summarizes the timing of crack propagation and the start and end of crack branching. The results show that the Cosserat length scale has less effect on crack propagation and branching. Figures 53 and 54 plot the contour of the damage variable at $t = 42.5\ \mu s$ and the contour of the micro-rotation of material points at $t = 35\ \mu s$ for three cases, respectively. As shown in 54, the micro rotation of material points decreases considerably as the Cosserat length scale increases. The results in Figure 53 show a similarity between the three cases. In contrast, the simulations with a larger Cosserat length scale show that the crack branching advances further from the initial branching point. This may be due to the same horizon size adopted for the three simulations.

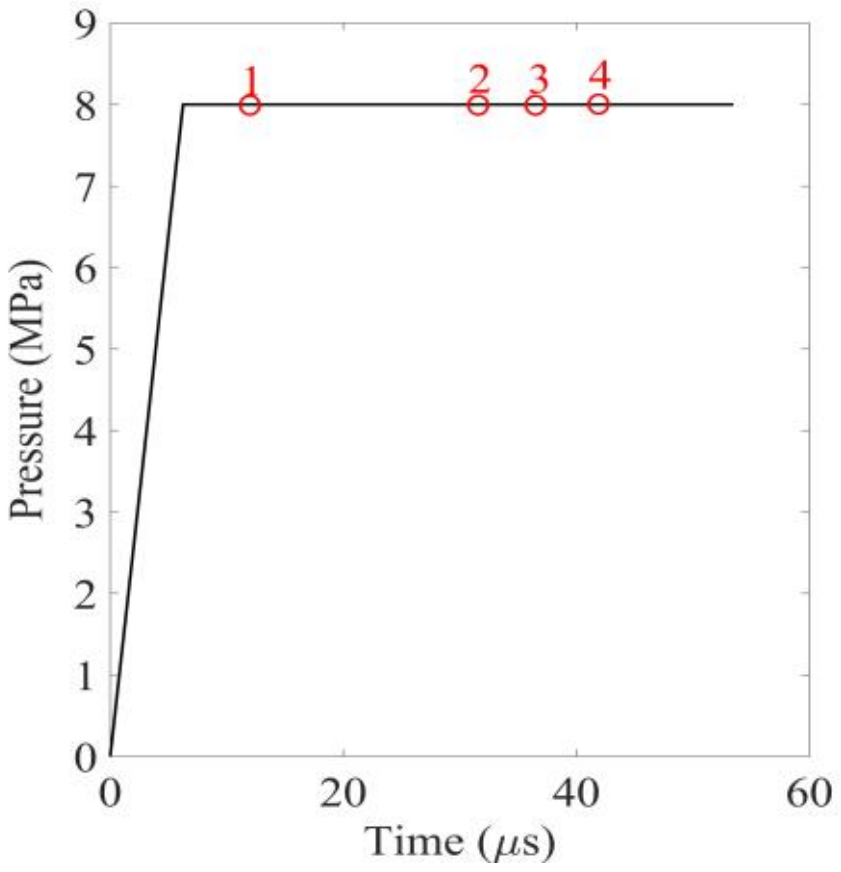

Figure 48: Applied load versus time.

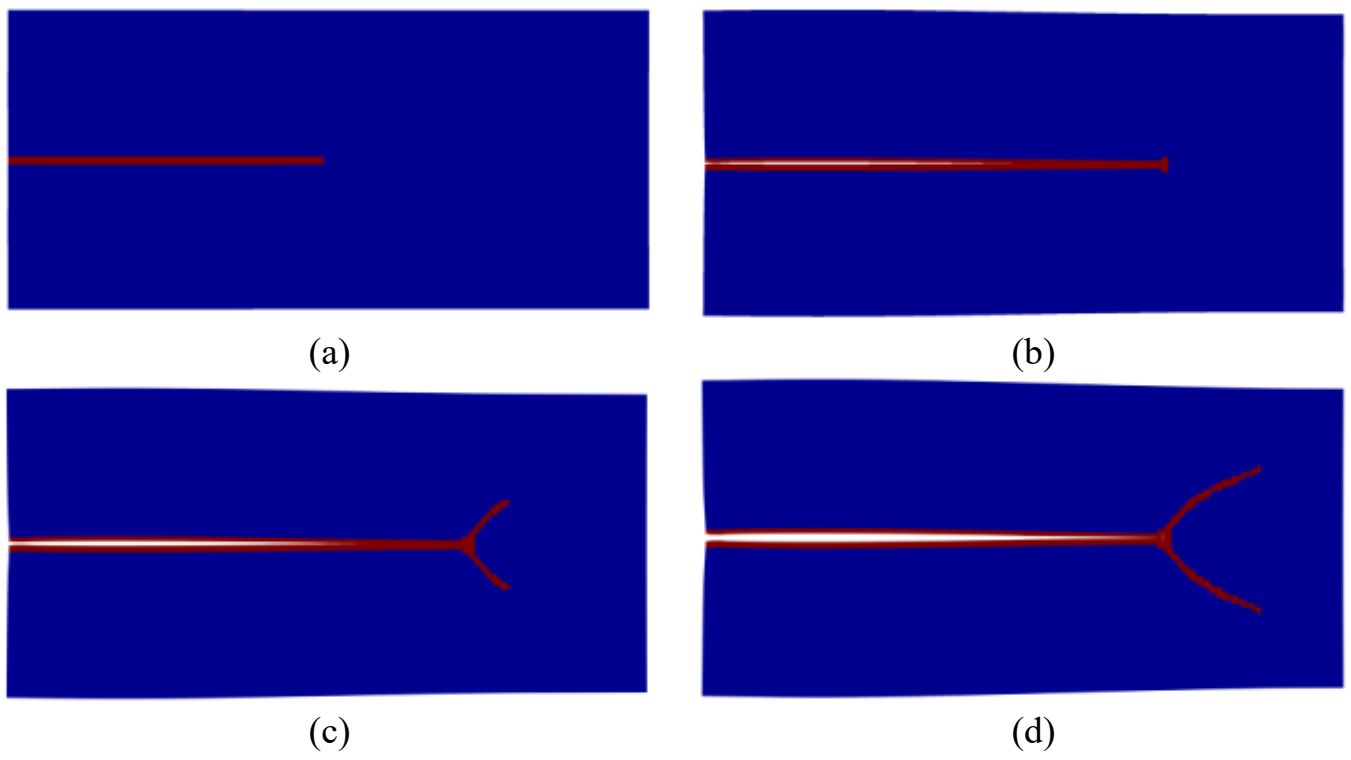

Figure 49: Contours of the crack path on the deformed configurations (magnification factor = 10) at (a) $t_1$ =12.75 $\mu$s, (b) $t_2$ =29.75 $\mu$s, (c) $t_3$ =35 $\mu$s, and (d) $t_4$ =42.5 $\mu$s.

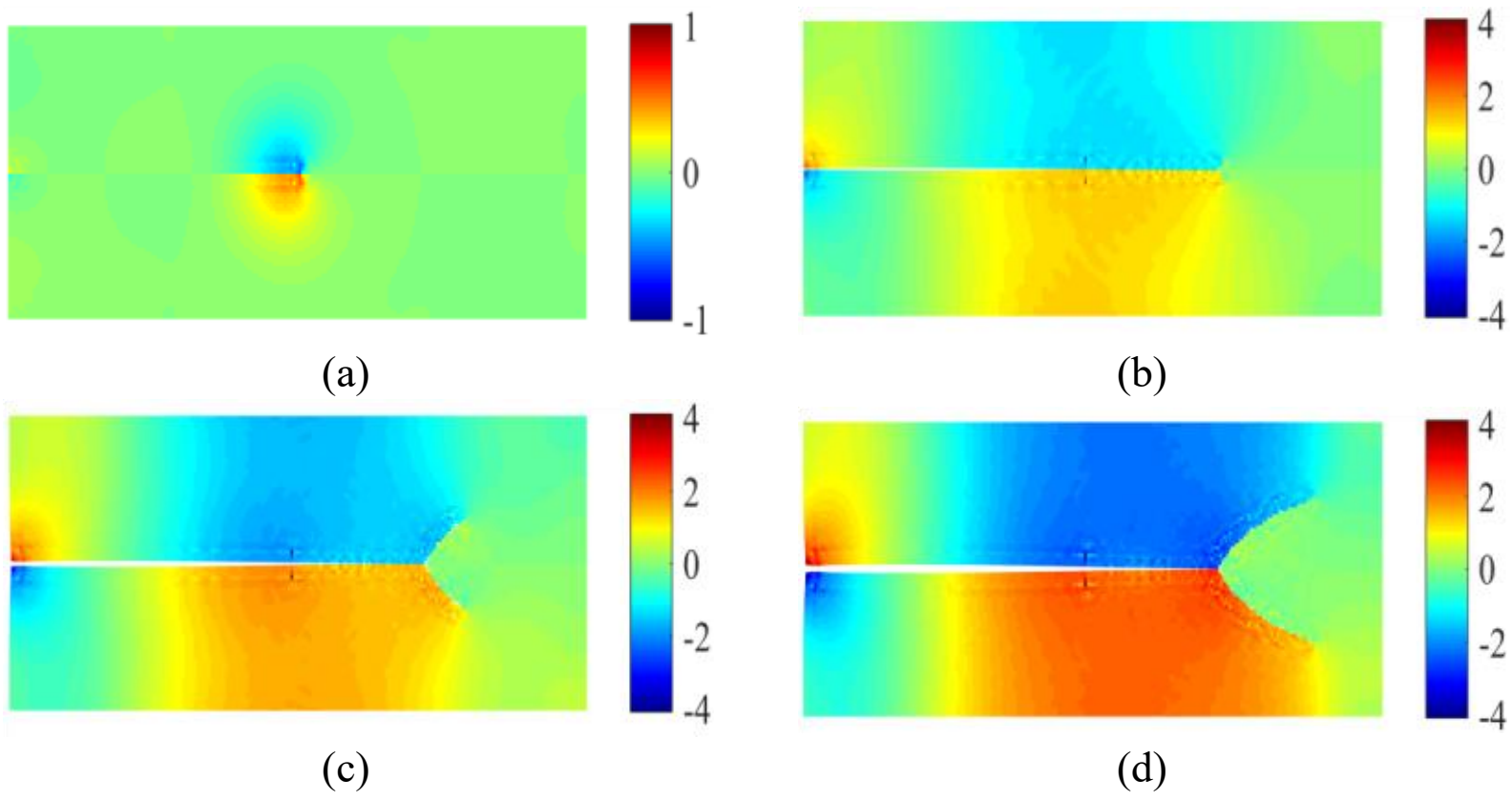

Figure 50: Contours of the micro rotation (0.001 rad) on the deformed configurations (magnification factor = 10) at (a) $t_1$ =12.75 $\mu$s, (b) $t_2$ =29.75 $\mu$s, (c) $t_3$ =35 $\mu$s, and (d) $t_4$ =42.5 $\mu$s.

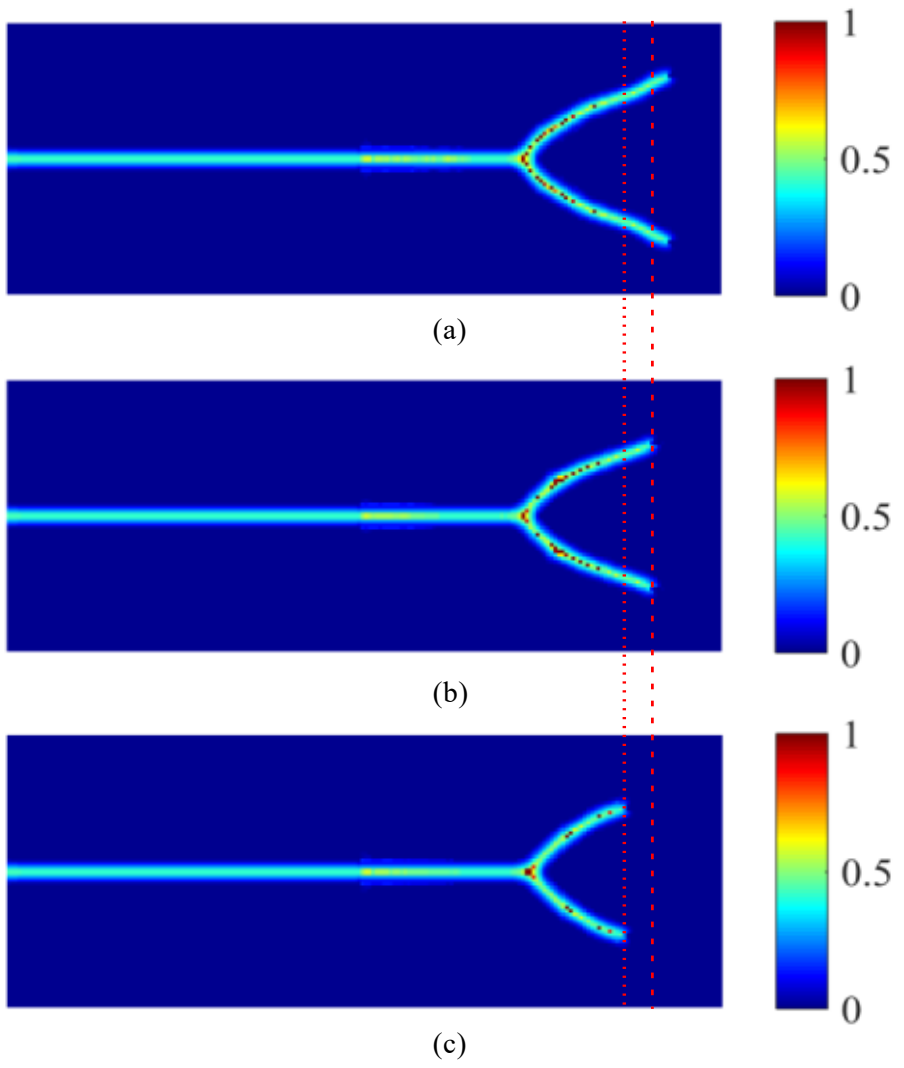

Figure 51: Contours of the damage variable at $t$ = 42.5 $\mu$s for three loading rates: (a) $t_{0,1}$ = 2.5 $\mu$s, (b) $t_{0,2}$ = 6.25 $\mu$s, and (c) $t_{0,3}$ = 12.5 $\mu$s. Note: the vertical dashed lines are plotted to compare the difference of crack branching.

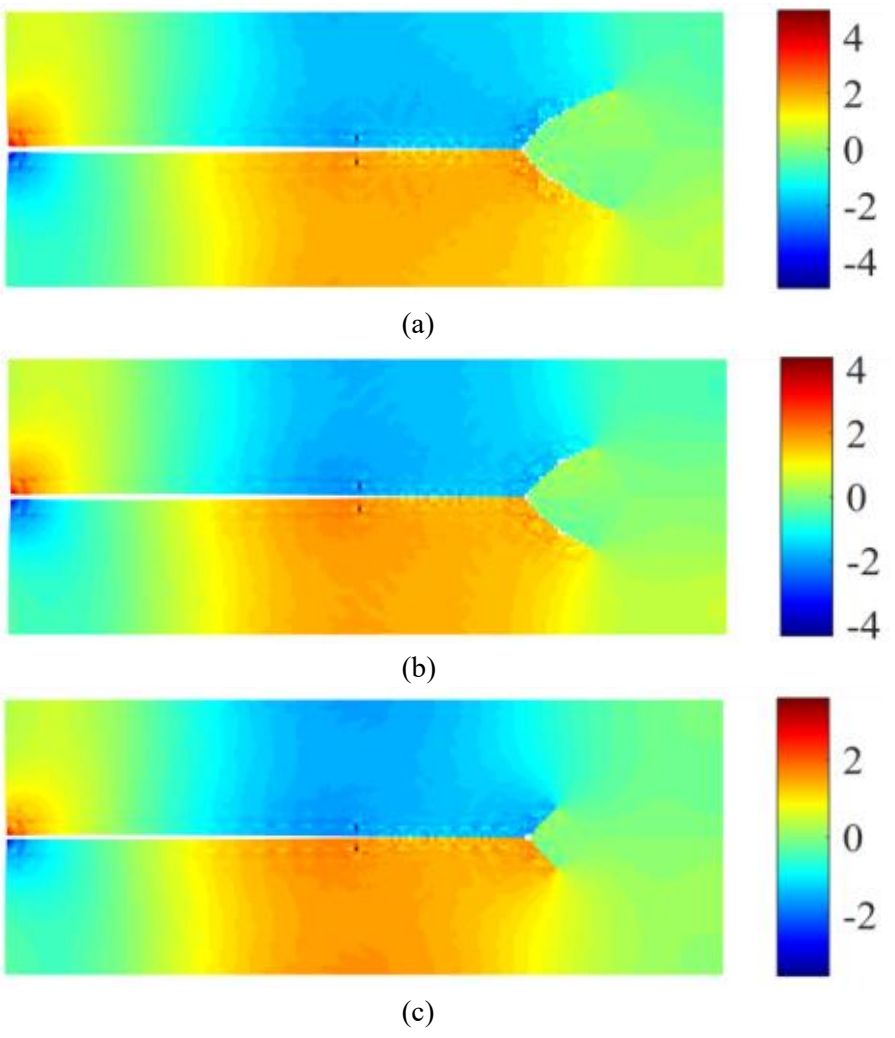

Figure 52: Contours of the micro rotation (×1$^{-3}$ rad) at $t$ = 35 $\mu$s for three loading rates: (a) $t_{0,1}$ = 2.5 $\mu$s, (b) $t_{0,2}$ = 6.25 $\mu$s and (c) $t_{0,3}$ = 12.5 $\mu$s.

Table 2: Summary of the timing of crack growth and branching for the three loading rates.

| $t_0$($\mu$s) | Start of crack growth ($\mu$s) | Start of branching ($\mu$s) | End of branching ($\mu$s) |
|---|---|---|---|
| 2.5 | 10 | 26.55 | 50.62 |
| 6.25 | 12.25 | 28.25 | 53.75 |
| 12.5 | 16 | 32.75 | 58 |

Table 3: Summary of the timing of crack growth and branching for different Cosserat length scales.

| $l$(mm) | Start of growth ($\mu$s) | Start of branching ($\mu$s) | End of branching ($\mu$s) |
|---|---|---|---|
| 1 | 12.25 | 28.5 | 55.5 |
| 2 | 12.25 | 28.25 | 52.75 |
| 3 | 12.12 | 28.12 | 48.75 |

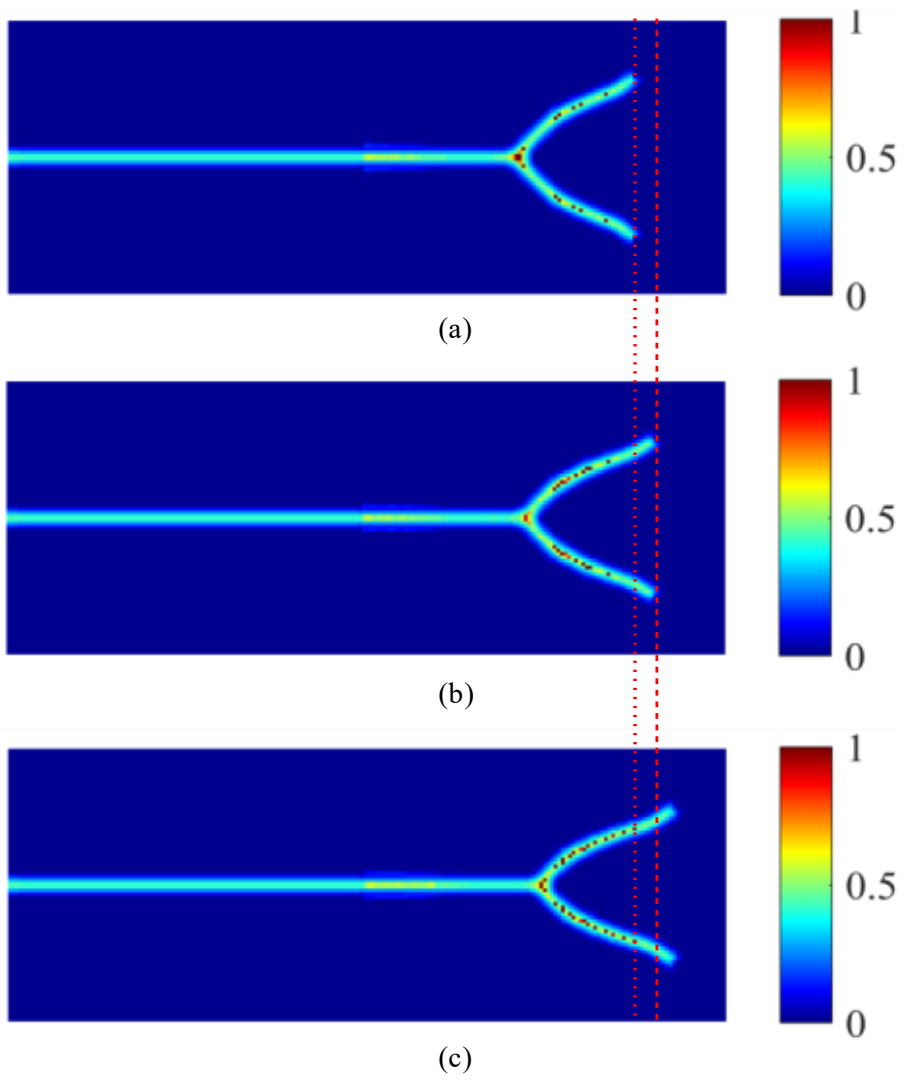

Figure 53: Contours of the damage variable at $t = 42.5\ \mu s$ for three Cosserat length scales: (a) $l_1 = 1$ mm, (b) $l_2 = 2$ mm, and (c) $l_3 = 3$ mm.

### *4.4.3. Influence of the initial volume fraction*

This part analyzes how the initial volume fraction can affect the crack branching. Three initial volume fractions considered are 0.75, 0.85, and 0.95. It is assumed that $l = 2$ mm and $t_0 = 6.25\ \mu s$ and the other parameters remain the same as the base simulation. Table 2 summarizes the timing of crack propagation and branching. It can be concluded from the results in Table 2 that the increase of the initial volume fraction of the porous material demands a larger load to initiate crack growth and branching. Figures 55 and 56 plot the contours of the damage variable at $t = 42.5\ \mu s$ and the contours of the micro rotation of material points at $t = 35\ \mu s$ for the three cases, respectively. As indicated by the results in Figures 55 and 56, the model with the largest initial volume fraction generates the smallest crack propagation and branching and the microrotation of material points around the crack.

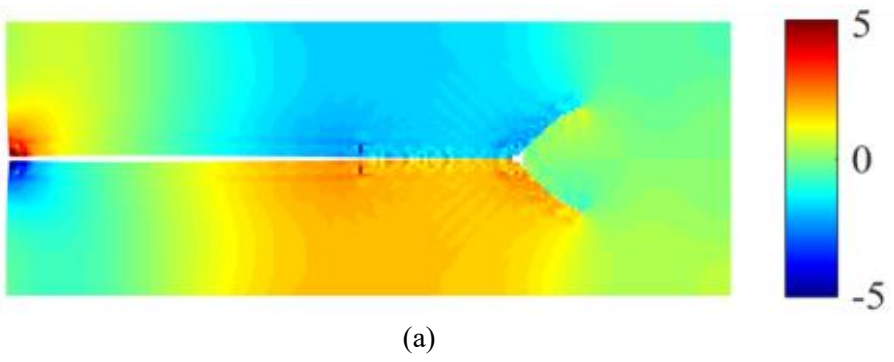

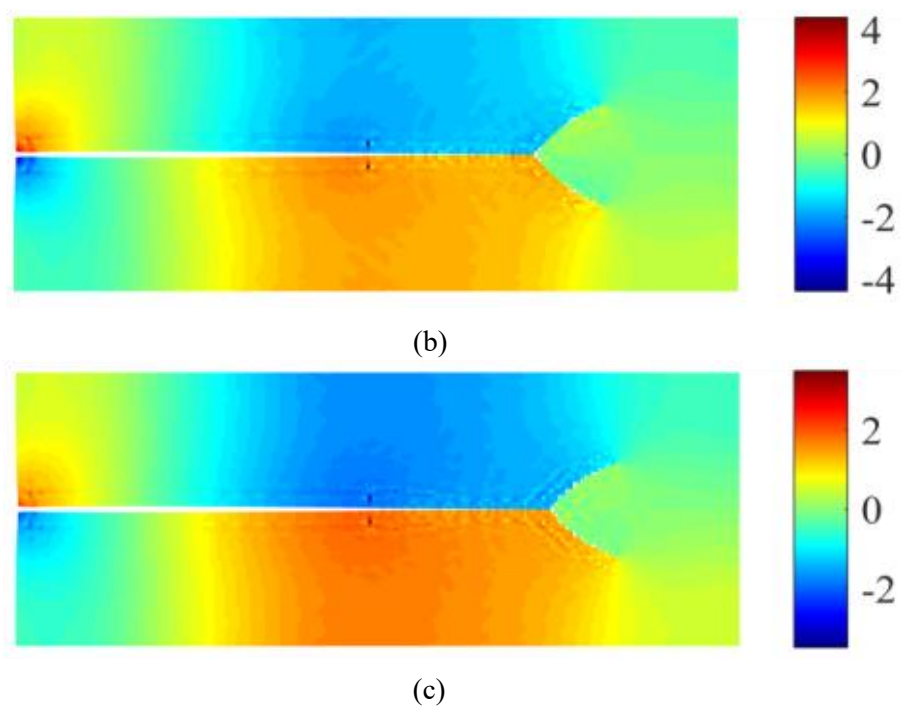

Figure 54: Contours of the micro rotation (×$10^{-3}$ rad) at $t$ = 35 $\mu$s for three Cosserat internal length scales: (a) $l_1$ = 1 mm, (b) $l_2$ = 2 mm, and (c) $l_3$ = 3 mm.

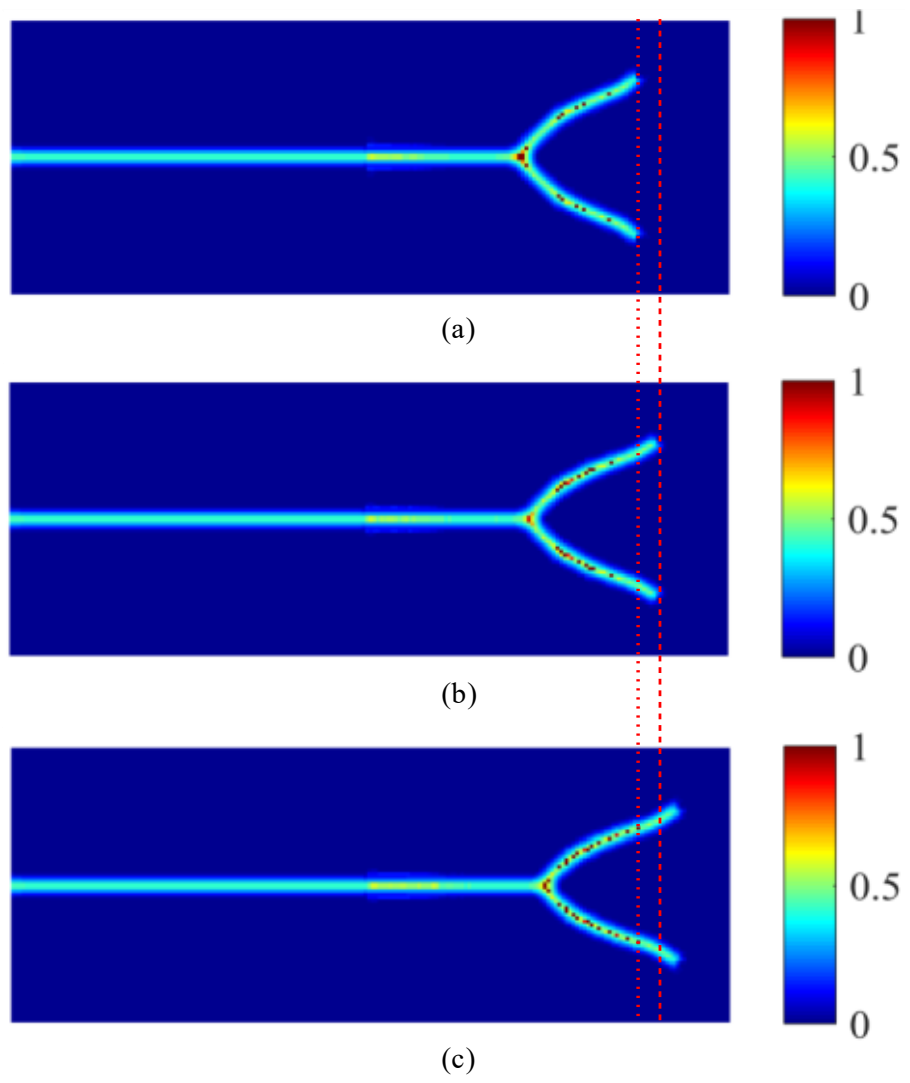

Figure 55: Contours of the damage variable at $t$ = 42.5 $\mu$s for three initial volume fractions: (a) $\varphi_{0,1}$ = 0.95, (b) $\varphi_{0,2}$ = 0.85, and (c) $\varphi_{0,3}$ = 0.75.

Table 4: Summary of the timing of crack growth and branching for three initial volume fractions.

| $\varphi_0$ | Start of growth ($\mu$s) | Start of branching ($\mu$s) | End of branching ($\mu$s) |
|---|---|---|---|
| 0.75 | 11.62 | 26.75 | 49.5 |
| 0.85 | 12.25 | 28.25 | 52.75 |
| 0.95 | 12.75 | 29.75 | 55 |

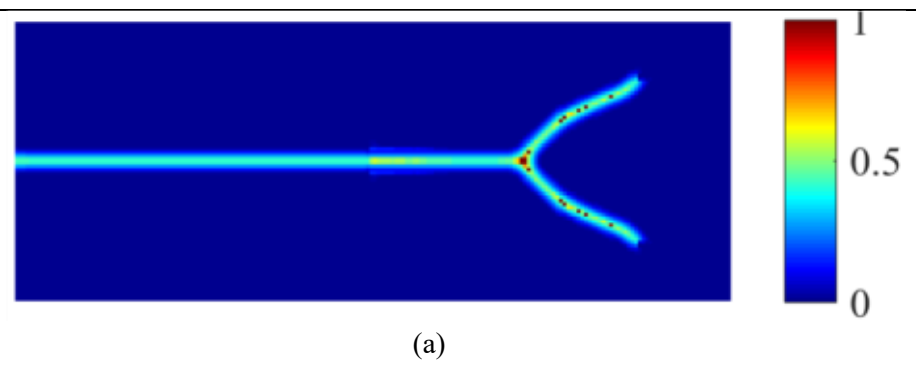

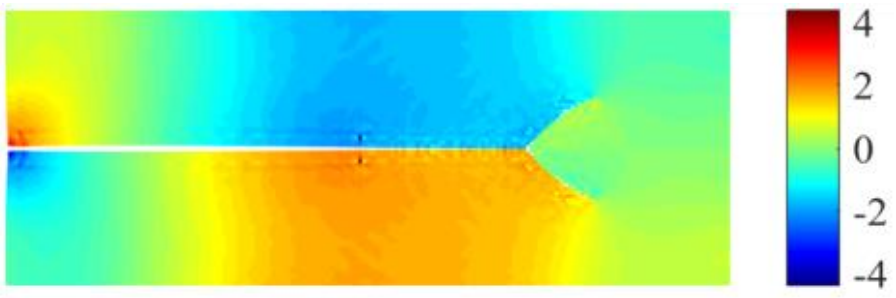

(b)

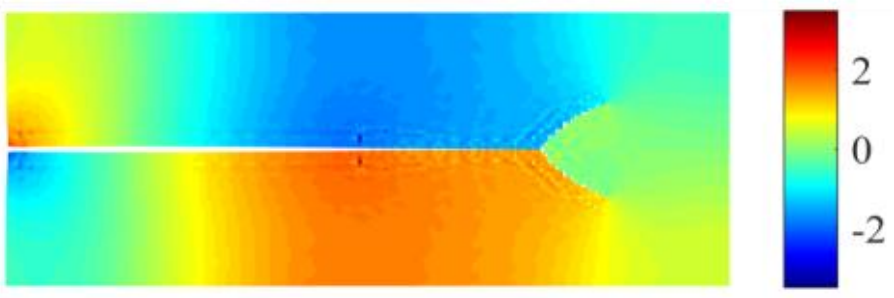

(c)

Figure 56: Contours of the micro rotation of material points ($\times 1^{-3}$ rad) at $t$ = 35 $\mu$s for three initial volume fractions: (a) $\varphi_{0,1}$ = 0.95, (b) $\varphi_{0,2}$ = 0.85, and (c) $\varphi_{0,3}$ = 0.75.

## 5. Closure

In this article, we formulate a viscous Cosserat periporomechanics paradigm for modeling dynamic shear banding and crack branching in dry porous media, in which a micro-structure based length scale, i.e., the Cosserat length scale, is incorporated. In this micro-periporomechanics paradigm, each material point has both translational and rotational degrees of freedom in line with the classical Cosserat continuum theory. The two field equations consisting of the force balance equation and moment balance equation are cast using the effective force state and moment state. The energy method is used to formulate the Cosserat periporomechanics correspondence principle for incorporating the classical micro-polar viscoplastic and viscoelastic constitutive models. We have also demonstrated that the Cosserat periporomechanics correspondence principle inherits zero-energy mode instability. To circumvent this numerical stability issue, we formulate a stabilized Cosserat constitutive correspondence principle through which classical viscous material models for porous media can be used in the proposed Cosserat periporomechanics. This new periporomechanics paradigm has been numerically implemented through an explicit Lagrangian meshfree algorithm for modeling dynamic failure in porous media. Numerical examples are presented to validate the computational Cosserat periporomechanics paradigm in modeling shear bands and mode-I cracks and demonstrate its efficacy and robustness in modeling dynamic shear banding and crack branching in dry porous media. We note that the classical second-order work criterion for detecting material instability incorporating micro-rotations of material points is utilized to validate the numerical results of shear banding. Meanwhile, through the numerical examples, we have analyzed the factors that could impact the dynamic shear banding and crack branching in dry porous media, such as loading rates, Cosserat length scales, and initial volume fractions.

## Acknowledgment

This work has been supported by the US National Science Foundation under contract number 1944009.